\documentclass[11pt,a4paper]{article}

\usepackage{amsmath,amssymb}
\newcommand{\e}{\varepsilon}
\newcommand{\va}{\varphi}
\newcommand{\D}{\Delta}

\newcommand{\n}{\nabla}

\newcommand{\N}{\frac{N}{2}}
\newcommand{\NN}{\frac{N}{p}}

\newcommand{\p}{\partial}

\newcommand{\R}{\mathbb{R}}
\newcommand{\h}{\hookrightarrow}

\def\hat{\widehat}

\newtheorem{theorem}{Theorem}
\newtheorem{proposition}{Proposition}
\newtheorem{prop}{Proposition}

\newtheorem{remarka}{Remark}

\newtheorem{lemme}{Lemma}

\newtheorem{lem}{Lemma}

\title{Global $\mbox{bmo}^{-1}(\R^N)$ radially symmetric solution for compressible Navier-Stokes equations with initial density in $L^\infty(\R^N)$}
\author{Boris Haspot  \thanks{Universit\'e Paris Dauphine, PSL Research University, Ceremade, Umr Cnrs 7534, Place du Mar\' echal De Lattre De Tassigny 75775 Paris cedex 16 (France), haspot@ceremade.dauphine.fr }
\thanks{ANGE project-team (Inria, Cerema, UPMC, CNRS), %
2 rue Simone Iff, CS 42112, 75589  Paris, France.}}
\date{}
\begin{document}
\maketitle

\begin{abstract}
In this paper we investigate the question of the existence of global weak solution 
for the compressible Navier Stokes equations provided that the initial momentum $\rho_0 u_0$ belongs to $\mbox{bmo}^{-1}(\R^N)$ with $N= 2,3$ and is radially symmetric. More precisely we deal with the so called viscous shallow water system when the viscosity coefficients verify $\mu(\rho)=\mu\rho$, $\lambda(\rho)=0$ with $\mu>0$. We prove then a equivalent of the so called Koch-Tataru theorem for the compressible Navier-Stokes equations. In addition we assume that the initial density $\rho_0$ is only bounded in $L^\infty(\R^N)$, it allows us in particular to consider initial density admitting shocks. Furthermore we show that if the coupling between the density and the velocity is sufficiently strong, then the initial density which admits initially shocks is instantaneously regularizing inasmuch as the density becomes Lipschitz. This coupling is expressed via the regularity of the so called effective velocity $v=u+2\mu\n\ln\rho$. In our case $v_0$ belongs to $L^2(\R^N)\cap L^\infty(\R^N)$, it is important to point out that this choice on the initial data implies that we work in a setting of infinite energy on the initial data $(\rho_0,u_0)$, it extends in particular the results of \cite{V}. In a similar way, we consider also the case of the dimension $N=1$ where the momentum $\rho_0 u_0$ belongs to $bmo^{-1}(\R)$ without any geometric restriction.\\
To finish we prove the global existence of strong solution for large initial data provided that the initial data are radially symmetric and sufficiently regular in dimension $N=2,3$ for $\gamma$ law pressure.
\end{abstract}
\section{Introduction}
We consider the compressible Navier Stokes system:
\begin{equation}
\begin{cases}
\begin{aligned}
&\p_{t}\rho+{\rm div}(\rho u)=0,\\
&\p_{t}(\rho u)+{\rm div}(\rho u\otimes u)-2{\rm div}(\mu(\rho) D(u))-\n (\lambda(\rho) {\rm div}u )+\n P(\rho)=0,\\
&(\rho,u)_{t=0}=(\rho_0,u_0).
\end{aligned}
\end{cases}
\label{0.1}
\end{equation}
Here $u=u(t,x)\in\R^{N}$ (with $N=1,2,3$) stands for the velocity field, $\rho=\rho(t,x)\in\R^{+}$ is the density, $D(u)=\frac{1}{2}(\n u+^{t}\n u)$ is the strain tensor and $P(\rho)$ is the pressure (we will only consider $\gamma$ law in the sequel which corresponds to $P(\rho)=a\rho^\gamma$ with $a>0$ and $\gamma\geq 1$). We denote by $\lambda$ and $\mu$ the two viscosity coefficients of the fluid. We supplement the problem with initial condition $(\rho_{0},u_{0})$.
Throughout the paper, we assume that the space variable $x$ belongs to $\R^{N}$ with $N\geq 1$. In the sequel we shall only consider the viscous shallow-water system, in this case the viscosity coefficients verify:
\begin{equation}
\mu(\rho)=\mu\rho\;\;\mbox{with}\;\mu>0\;\;\;\;\mbox{and}\;\;\;\;\lambda(\rho)=0.
\label{condition}
\end{equation}
It is motivated by the physical consideration
that in the derivation of the Navier-Stokes equations from the Boltzmann equation through the
Chapman-Enskog expansion to the second order (see \cite{CC70}), the viscosity coefficient is a function
of the temperature. If we consider the case of isentropic fluids, this dependence is reduced to the
dependence on the density function.  It is worth pointing out that the case $\mu(\rho)=\mu\rho$ corresponds to the so called viscous shallow water system. This system with friction has been derived by Gerbeau and Perthame in \cite{GP} from the Navier-Stokes system with a free moving boundary in the shallow water regime at the first order (it corresponds to a small shallowness parameter). This derivation relies on the hydrostatic approximation where the authors follow the role of viscosity and friction on the bottom. \\
We are now going to rewrite the system (\ref{0.1})
following the new formulation proposed in \cite{para} (see also \cite{BrDeZa,CPAM,CPAM1,glob}), indeed setting $v=u+2\mu\n\ln\rho$ we can rewrite the system (\ref{0.1}) as follows:
\begin{equation}
\begin{cases}
\begin{aligned}
&\p_t\rho-2\mu\D\rho+{\rm div}(\rho v)=0,\\
&\rho\p_t v+\rho u\cdot\n v-\mu{\rm div}(\rho{\rm curl}v)+\n P(\rho)=0,
\end{aligned}
\end{cases}
\label{0.1a}
\end{equation}
with ${\rm curl}v={\rm curl}u=\n u-^t\n u$. An other way is also to rewrite the system (\ref{0.1}) as follows by using the unknown $m$ corresponding to the momentum velocity $m=\rho u$:
\begin{equation}
\begin{cases}
\begin{aligned}
&\p_{t}\rho+{\rm div}m=0,\\
&\p_{t}m+\frac{1}{2}\big({\rm div}(m\otimes v)+{\rm div}(v\otimes m)\big)-\mu\D m-\mu\n{\rm div}m+\n P(\rho)=0,\\
&\p_t(\rho v)+{\rm div}(m\otimes v)-\mu{\rm div}(\rho{\rm curl}v)+\n P(\rho)=0,\\
&(\rho,m,v)_{t=0}=(\rho_0,m_0,v_0).
\end{aligned}
\end{cases}
\label{0.1m}
\end{equation}
This system is an augmented system since we deal with three unknowns $(\rho,m,v)$.
In the one dimensional case, we observe that the previous system (\ref{0.1a})
can be write simply under the following form with $v=u+2\mu\p_x\ln\rho$:
\begin{equation}
\begin{cases}
\begin{aligned}
&\p_t\rho-2\mu\p_{xx}\rho+\p_x(\rho v)=0,\\
&\rho\p_t v+\rho u\p_x v+\p_x P(\rho)=0,
\end{aligned}
\end{cases}
\label{0.1ad}
\end{equation}
It is important to point out that systematically in this case the effective velocity $v$ verifies essentially a  transport equation. We recall now the classical energy inequality of the system (\ref{0.1}) when the initial density $\rho_0$ is close from a constant state $\bar{\rho}>0$. Multiplying the momentum equation by $u$ we have at least heuristically: 
\begin{equation}
\begin{aligned}
&{\cal E}(\rho,u)(t)=\int_{\R^N}\big(\frac{1}{2}\rho(t,x)|u(t,x)|^{2}+(\Pi(\rho)(t,x)-\Pi(\bar{\rho}))\big)dx+\int^{t}_{0}\int_{\R^N}2\mu\,\rho|D u|^{2}(t,x) dt dx\\
&\hspace{3cm} \leq \int_{\R^N}\big(\rho_{0}(x)|u_{0}(x)|^{2}+(\Pi(\rho_{0})(x)-\Pi(\bar{\rho}))\big)dx={\cal E}(\rho_0,u_0).
\end{aligned}
\label{cenergie}
\end{equation}
with $\Pi(\rho)$ defined as follows $\Pi(s)=s\big(\int^s_{\bar{\rho}}\frac{P(z)}{z^2}dz-\frac{P(\bar{\rho})}{\bar{\rho}}\big)$.
Notice that $\Pi$ is convex in the case of the $\gamma$ law. In order to give a sense to (\ref{cenergie}), it is important to point out that the initial data must be of finite energy inasmuch as $\sqrt{\rho_0}u_0$ and $(\Pi(\rho_0)-\Pi(\bar{\rho}))$ must be respectively in $L^2(\R^N)$ and in $L^1(\R^N)$. Similarly we have also the BD entropy (see \cite{BD}), multiplying now momentum equation of (\ref{0.1a}) by $v$ we obtain at least heuristically:
\begin{equation}
\begin{aligned}
&{\cal E}_1(\rho, v)(t)=\int_{\R^N}\big(\frac{1}{2}\rho(t,x)|v(t,x)|^{2}+(\Pi(\rho)(t,x)-\Pi(\bar{\rho}))\big)dx\\
&+\int^{t}_{0}\int_{\R^N}\mu\,\rho(t,x)|{\rm curl} v|^{2}(t,x) dt dx+\frac{8\mu}{\gamma}\int^t_0\int_{\R^N}|\n\rho^{\frac{\gamma}{2}}|^2 dx dt\\
&\hspace{3cm} \leq \int_{\R^N}\big(\rho_{0}(x)|v_{0}(x)|^{2}+(\Pi(\rho_{0})(x)-\Pi(\bar{\rho})\big)dx={\cal E}_1(\rho_0,v_0).
\end{aligned}
\label{bcenergie}
\end{equation}
As previously, in order to obtain the BD entropy the initial data must verify $\sqrt{\rho_0}v_0\in L^2(\R^N)$ and $(\Pi(\rho_0)-\Pi(\bar{\rho}))\in L^1(\R^N)$ . In particular it is possible to choose initial data for which the BD entropy  (\ref{bcenergie}) is finite but the classical energy (\ref{cenergie}) infinite. In the sequel we will assume that $\bar{\rho}=1$ in order to simplify the notations.\\
When $N\geq 2$, we  now restrict our attention to the case of radially invariant initial velocity and radial initial density:
\begin{equation}
u_0(x)=\n\theta(|x|)\;\;\mbox{and}\;\;\rho_0(x)=\rho_1(|x|),
\label{initial}
\end{equation}
with $\theta$ a radial function. Under the assumption that this condition is preserved all along the time $t>0$:
\begin{equation}
u(t,x)=\n\theta(t,|x|)\;\;\;\mbox{and}\;\;\;\rho(t,x)=\rho_1(t,|x|),
\label{preserv}
\end{equation}
 it implies that ${\rm curl}v(t,\cdot)={\rm curl}u(t,\cdot)=0$ all along the time $t\geq 0$. We can then rewrite the system (\ref{0.1a}) as follows (provided that we show that (\ref{preserv}) is true all along the time):
\begin{equation}
\begin{cases}
\begin{aligned}
&\p_{t}\rho-2\mu\D\rho+{\rm div}(\rho v)=0,\\
&\rho \p_{t}v+\rho u\cdot\n v+\frac{a\gamma}{2\mu}\rho^\gamma v=\frac{a\gamma}{2\mu}\rho^\gamma u\\
&\p_t m+{\rm div}(v\otimes m)-2\mu\D m+\frac{a\gamma}{2\mu}\rho^\gamma v=\frac{a\gamma}{2\mu}\rho^\gamma u.
\end{aligned}
\end{cases}
\label{0.3}
\end{equation}
We observe that in this case the effective velocity verifies a damped transport equation as it is  systematically the case in one dimension. It will allow us to estimate simply the effective velocity $v$ in $L^\infty(\R^N)$ norm.\\
In this paper we address the question of the existence of global weak solution for initial density admitting some shocks. This is  motivated by the fact that discontinuous solutions are fundamental both in the physical theory of non-equilibrium thermodynamics as well as in the mathematical study of inviscid models for compressible flow. It is important, therefore, to understand in which functional setting we can deal with shocks on the initial density. Indeed an important question is to understand if it is possible to prove the existence of strong solution in a functional setting which includes this case. An other formulation is, can we prove the existence of global weak solution for initial density verifying only $(\rho_0-1)\in  L^\infty(\R^N)\cap L^2(\R^N)$ and is it possible to prove some regularizing effects on the velocity which are stronger than the energy estimates? 
We mention in particular that in \cite{Li,V}  the authors proved independently the existence of global weak solution provided that the initial data verify the classical energy estimate and the BD entropy. It implies in particular that $\n\sqrt{\rho_0}$ belongs to $L^2(\R^N)$. Unfortunately this assumption is too strong to deal with initial density admitting shocks (we can think for example to the case $\rho_0=a$ on $B(0,r)$ and $b$ on $^c B(0,r)$ with $a,b,r>0$).
\\
The second question is to understand the propagation of the shocks on the density all along the time? Indeed is it true that for $t>0$ the density $\rho(t,\cdot)$ admits some shocks (and how to describe them?) or is it possible to imagine that the density becomes regular in finite time? Similarly how can we describe the regularizing effects on the velocity? The results of this paper give some particular answers to these questions provided that the initial data are radially symmetric. Indeed we prove  the existence of global weak solution with initial density admitting shocks and with momentum $m=\rho u$ verifying all along the time $\log$ Lipschitz estimates. Roughly speaking, it is well known that Lipschitz estimate on the velocity is often a minimal condition to prove the uniqueness of the solution in fluids mechanics (we refer in particular to the so called Beale-Kato-Majda criterion \cite{Beale} or to the work of Danchin for compressible Navier Stokes equations \cite{Fourier}), it explains why it is not clear that our solution are also unique \footnote{We will give more indications on this question in the sequel.}.\\
 In order to obtain such results, a key point consists in working with initial data with minimal regularity assumption. Since we wish to deal with the largest functional space $X$ in which
local or global well-posedness may be proved, it is then important to understand in which sense a space is \textit{critical} for the existence of strong solution. In consequence the notion of invariance by scaling for the equations plays a crucial role, we will say that a functional space is critical if we can solve the system (\ref{0.1}) in functional spaces with norm
 invariant by the changes of scales which leave (\ref{0.1}) invariant. We can observe that the transformations:
\begin{equation}
\begin{aligned}
&(\rho_{0}(x),u_{0}(x))\rightarrow (\rho_{0}(lx),lu_{0}(lx)),\;\;\;(\rho(t,x),u(t,x))\rightarrow (\rho(l^{2}t,lx),lu(l^{2}t,lx))
\end{aligned}
\label{2}
\end{equation}
have that property for $l>0$, provided that the pressure term has been changed accordingly. A good candidate are the following homogeneous Besov spaces (see \cite{Danchin} for some definitions) $\dot{B}^0_{\infty,\infty}(\R^N)\times(\dot{B}^{-1}_{\infty,\infty}(\R^N))^N$, it is well known that it is the largest critical space ( see \cite{Lemarie}). In the case of the incompressible Navier-Stokes equations, it has been proved by Canonne et al \cite{Can} that the equations are globally well posed for small initial velocity $u_0$ in the following Besov space $\dot{B}^{\NN-1}_{p,\infty}(\R^N)$  \footnote{These spaces are embedded in $\dot{B}^{-1}_{\infty,\infty}(\R^N)$} with $p<+\infty$. In \cite{Koch}, Koch and Tataru haved proved the global existence of strong solution for small initial data $u_0$ in $\mbox{BMO}^{-1}(\R^N)$. It is important to point out that this space is the largest space in which the existence of global strong solution has been showed for incompressible Navier-Stokes equations. Let us recall precisely the theorem of Koch and Tataru (see also the results of Bourgain and Pavlovi\'c in \cite{BP} of ill-posedness in $\dot{B}^{-1}_{\infty,\infty}(\R^N)$). Let $e^{t\D}u_0$ be the solution to the heat equation with initial data $u_0$. We have then:
$$e^{t\D}u_0=u_0*\phi_{\sqrt{4t}}\;\;\mbox{with}\;\phi(x)=\pi^{-\N}e^{-|x|^2}\;\;\mbox{and}\;\;\phi_t (x)=t^{-N}\phi(\frac{x}{t}).$$
Therefore $BMO^{-1}(\R^N)$ is the set of temperated distribution $u_0$ for which the following norm is finite:
$$\|u_0\|_{BMO^{-1}(\R^N)}=\sup_{x\in\R^N,t>0}\big( t^{-\N}\int^t_0\int_{B(x,\sqrt{t})}|e^{s\D}u_0 (y)|^2 dy ds\big)^{\frac{1}{2}}.$$
As in \cite{Koch,Lemarie}, we define also $bmo^{-1}(\R^N)$ the set of temperated distribution $u_0$ for which for all $T>0$ we have $\sup_{x\in\R^N,0<t<T}\big( t^{-\N}\int^t_0\int_{B(x,\sqrt{t})}|e^{s\D}u_0 (y)|^2 dy ds\big)^{\frac{1}{2}}<+\infty.$
We define the norm of $bmo^{-1}(\R^N)$ by:
$$\|u_0\|_{bmo^{-1}(\R^N)}=\sup_{x\in\R^N,0<t<1}\big( t^{-\N}\int^t_0\int_{B(x,\sqrt{t})}|e^{s\D}u_0 (y)|^2 dy ds\big)^{\frac{1}{2}}.$$
In addition we have the following proposition (see \cite{Koch,Lemarie}).
\begin{prop}
\label{proptech}
Let $u_0$ be a tempered distribution. Then $u_0\in BMO^{-1}(\R^N)$ if and only if there exists $f_i\in BMO(\R^N)$ with $u_0=\sum_i\p_i f_i$. Similarly $u_0$ belongs to $bmo^{-1}(\R^N)$ if and only if there exists $N+1$ distributions $f_0,f_1,\cdots,f_N\in bmo(\R^N)$ such that $u_0=f_0+\sum_{1\leq i\leq N}\p_i f_i$.
\end{prop}
We refer to the Chapter 10 of \cite{Lemarie} for the definition of $BMO(\R^N)$ and $bmo(\R^N)$. Let us state now the theorem of Koch and Tataru \cite{Koch}. 
\begin{theorem}
[Koch-Tataru \cite{Koch}] Let $N\geq 2$. The incompressible Navier Stokes equations:
$$
\begin{cases}
\begin{aligned}
&\p_t u+{\rm div}(u\otimes u)-\D u+\n \Pi=0,\,{\rm div}u=0,
\end{aligned}
\end{cases}
$$
have a unique global small solution in $X$ provided that $u_0$ is small in $BMO^{-1}(\R^N)$  with ${\rm div} u_0=0$. Here $X$ is defined as follows:
\begin{equation}
\|u\|_{X}=\sup_{t\in\R^+}\sqrt{t}\|u(t,\cdot)\|_{L^\infty(\R^N)}+\big(\sup_{x\in\R^N,t>0}t^{-\N}\int^t_0\int_{ B(x,\sqrt{t})}|u|^2(s,y) dt dy\big)^{\frac{1}{2}},
\label{defXs}
\end{equation}
\label{theo1}
\end{theorem}
We mention that the previous result have been extended to general parabolic equations with quadratic nonlinearities including the particular case of the Navier Stokes equations when neither pointwise kernel bounds nor self-adjointness are available. We refer to the work of Auscher and Frey in \cite{Auscher} who develop a strategy making use of the tent spaces.
In the case of compressible Navier-Stokes equations with constant viscosity coefficients
, Danchin in \cite{DG} showed a result of existence of global strong solution with small initial data in critical space for the scaling of the system in dimension $N\geq 2$. More precisely the initial data are chosen as follows $(\rho_{0}-1,u_{0})\in (\dot{B}^{\N}_{2,1}(\R^N)\cap \dot{B}^{\N-1}_{2,1}(\R^N))\times \dot{B}^{\N-1}_{2,1}(\R^N)$. This result has been generalized in $\widetilde{B}^{\N-1,\NN}_{2,p,1}(\R^N)\times( \widetilde{B}^{\N-1,\NN-1}_{2,p,1}(\R^N))^N$
for $p$ suitably chosen in \cite{CD,CMZ,arma} with different methods ( see also \cite{FZZ} for some refinements in low frequencies). Let us make two remarks, the first one is that the initial density $(\rho_0-1)$ is systematically chosen belonging to a space $\widetilde{B}^{\N-1,\NN}_{2,p,1}(\R^N)$ which is embedded in $C^{0}(\R^N)$ the set of continuous functions decaying to $0$ at infinity.
In particular such results prevent the choice of initial density $\rho_0$ admitting shocks.
The second observation is that the velocity $u_0$ is generally chosen in a critical Besov space with third index equal to $1$. In particular it allows to prove that the velocity is Lipschitz all along the time or in other words $\n u$ belongs to $L^1(\R^+,L^\infty(\R))$.  It is important to point out that it is not the case for incompressible Navier Stokes equations since $u_0$ belongs only to $\dot{B}^{\NN-1}_{p,\infty}(\R^N)$ with $1\leq p<+\infty$ (see \cite{Can}).
For some extensions, we refer however to \cite{Geometrie}, where it is proven that there exists global strong solution with small initial velocity in $\dot{H}^{\N-1}$ when the initial velocity is assumed axisymmetric (it corresponds to the so called Fujita-Kato initial data \cite{Fujita}). More recently in \cite{DFP}, Danchin et al proved the existence in finite time for small initial density $\rho_0-1$ in $L^\infty(\R^N)$
with initial velocity slightly subcritical for constant viscosity coefficients. This result allows in particular to deal with density admitting shocks provided a small perturbation of a positive constant and at the condition to consider constant viscosity coefficients (indeed this last point is absolutely crucial in the analysis).
\\
In this paper we wish to establish the existence of global weak solution for the system (\ref{0.1}) for radially invariant initial data in a functional framework allowing to consider shocks on the initial data. 
In addition, we wish to exhibit some regularizing effects on the momentum $m=\rho u$  as it is proved in \cite{Koch} for the incompressible Navier Stokes equations. More precisely we are going to deal with initial density verifying $0\leq \rho_0\leq M<+\infty$ and with an effective velocity $v_0$ belonging to $L^\infty(\R)$. Furthermore $\rho_0$ and $u_0$ will verify the assumptions (\ref{initial}). From the definition of the effective velocity we have $m_0=\rho_0 u_0=\rho_0 v_0+2\mu\n\rho_0$, we deduce that $\rho_0 u_0$ is in $L^\infty(\R^N)+W^{-1,\infty}(\R^N)$. It implies using the proposition \ref{proptech} that 
$\rho_0 u_0\in bmo^{-1}(\R^N)$.
The proof of the result consists in combining $L^\infty$ estimates on the effective velocity $v$ and the Koch-Tataru estimates on the momentum $\rho u$ which verifies a heat equation with quadratic nonlinearities in $\rho u$ and $v$ \footnote{The fact that $v$ can be controlled in $L^\infty$ norm depends in a crucial way that $v$ verifies a damped transport equation (\ref{0.3}). It is not the case for general initial data since in (\ref{0.1a}) the evolution of $v$ depends also on the vorticity ${\rm curl}v$ which is generally not null.}. More precisely, we will prove in a first time the existence of weak solution in finite time, it means that there exists a finite time $T>0$ such that our solution $(\rho,v,m)$ belongs to the following space:
$$
\begin{aligned}
&Y_T=L^\infty_T(L^\infty(\R^N))\times (L^\infty_T(L^\infty(\R^N)))^N\times ({\cal E}_T\cap {\cal E}^1_T) , 
\end{aligned}
$$
with the norms of ${\cal E}_T$ and ${\cal E}^1_T$ defined as follows:
$$
\begin{aligned}
&\|m\|_{{\cal E}_T}=\sup_{0<t<T}\sqrt{t}\|m(t)\|_{L^\infty(\R^N)}+\big(\sup_{x\in\R^N,\,0<t<T}t^{-\N}\int^t_0 \int_{y\in B(x,\sqrt{t})}|m(t,y)|^2 dt \,dy\big)^{\frac{1}{2}}\\
&\|m\|_{{\cal E}_T^1}=\|m\|_{\widetilde{L}^\infty_T(\dot{B}^{-1}_{\infty,\infty}(\R^N))\cap\widetilde{L}_T^1(\dot{B}^1_{\infty,\infty}(\R^N))}.
\end{aligned}
$$
Concerning the definition of the Chemin-Lerner spaces for non homogeneous Besov spaces $\widetilde{L}^\infty_T(\dot{B}^{-1}_{\infty,\infty}(\R^N))\cap\widetilde{L}_T^1(\dot{B}^1_{\infty,\infty}(\R^N))$, we refer to \cite{Danchin}. \\
The previous choice allows us to choose initial density with shocks since we assume only that $\rho_0$ belongs to $L^\infty(\R^N)$. We observe also that instantaneously the density is regularized inasmuch as the density becomes Lipschitz. Indeed we have using the definition of the effective velocity $\n\rho=\frac{1}{2\mu}(\rho v-m). $
And by definition of $Y_T$, we observe that for $0<t\leq T$, $\n\rho(t,\cdot)$ belongs to $L^\infty(\R^N)$. Combining this information with the fact that $\rho\in L^\infty_T(L^\infty(\R^N))$, its shows that for $0<t<T$, $\rho(t,\cdot)$ belongs to $W^{1,\infty}(\R^N)$ and is in particular continuous and Lipschitz. This result is surprising if we consider that the density $\rho$ verifies an hyperbolic equation in (\ref{0.1}). This is due in fact to the regularity of the effective velocity $v$ which describes the coupling between the momentum $m$ and the density $\rho$. In other way if $v$ is sufficiently regular then we can expect regularizing effects on the density inasmuch as the density becomes instantaneously continuous. Next we will prove that we can extend beyond the time $T$ the solution in a global weak solution provided that the initial data verify the BD entropy. Since we have proved the existence of weak solution in finite time, it remains to deal with the long time and to do this we will show that the solution verifies classical energy estimate for long time and admits regularizing effects. We emphasize that our initial data are of infinite classical energy but instantaneously $\sqrt{\rho}u(t,\cdot)$ belongs to $L^2(\R^N)$ for $t>0$ if $\rho_0$ is far far away from the vacuum. Conversely we would like to recall some results due to  Hoff. In \cite{HJDE,H6,H7}
, Hoff consider the case of the compressible Navier-Stokes equations with constant viscosity. He proved the existence of global weak solution for small initial data provided that the initial density $(\rho_0-1)$ belongs only to $L^\infty(\R^N)\cap L^2(\R^N)$ with $\rho_0\geq c>0$ and that the velocity $u_0$ is sufficiently regular (the initial data are of finite energy in particular). In addition he showed regularizing effects on the vorticity ${\rm curl}u$ and the effective pressure $\mu {\rm div}u-P(\rho)+P(1)$. Compared with our result there is no regularizing effects on the density $\rho$ inasmuch as the density becomes continuous. This is probably due to the fact that the effective velocity $v$ which defines the coupling between the density and the velocity $u$ is not sufficiently regular and that there is no BD entropy for constant viscosity coefficients at least in dimension $N\geq 2$ (otherwise it is true). It means in particular that our results and the works of Hoff are opposite in the sense that the density admits or not regularizing effects and this is probably translated by the regularity of the effective velocity $v$.\\
We are going now to consider the one dimensional without any geometrical restriction. Similarly we prove the existence of global weak solution with initial density $(\rho_0-1)$ in $L^\infty(\R)\cap L^2(\R)$. We use the same ideas than previously and the fact that in this case the velocity $v$ verifies a damped transport equation (see (\ref{0.1ad})). In particular we extend the work \cite{glob} to the case of the viscous shallow water system. Indeed in \cite{glob}, we prove for small initial data the existence of global weak solution for initial density in $BV(\R)$ and for effective momentum $\rho_0 v_0$ in the set of finite measure ${\cal M}(\R)$. In the present work we do not assume any smallness assumption, however we deal with an effective velocity which is more integrable since $v_0$ is in $L^\infty(\R)$. We improve also the regularity assumptions on the density since it is only in $L^\infty(\R)$ (whereas in \cite{glob} the density is assumed to belong to $BV(\R)$). For the existence of global strong solution with large initial data in one dimension we refer to \cite{Glob,MV,Constantin}.\\
To finish we will prove the existence of global strong solution for the system (\ref{0.1}) (corresponding to the viscous shallow water system) for large initial data provided that the initial data are sufficiently regular and that the initial density and the initial velocity are respectively radial and  radially symmetric. We extend in particular the results of Choe and Kim \cite{CK1} to the case of euclidean space $\R^N$ with $N=2,3$. In \cite{CK1} the authors proved the existence of global strong solution for compressible Navier-Stokes equations with constant viscosity coefficients in an annular domain. The fact that the authors work in an annular domain allow to simplify the energy estimates on the velocity and the density since it avoids all the blow-up phenomena when $|x|$ goes to zero.
In some sense our result is an extension of \cite{Glob} to the multidimensional case when the data are radially symmetric. The main difficulty consists in proving that $\frac{1}{\rho}$ and $\rho$ remains in $L^\infty(\R^N)$ all along the time. To do this, we will use the fact that the effective velocity $v$ verifies a damped transport equation in order to obtain $L^\infty$ estimates on $v$. Using now the fact that the density $\rho$ is governed by a parabolic equation in (\ref{0.1a}), we will use the maximum principle and the $L^\infty$ estimate on $v$ to obtain $L^\infty$ estimate on $\rho$ and $\frac{1}{\rho}$. Combining classical blow up criterion and propagation of the regularity on the velocity $u$ when the density is bounded and far away from zero, we obtain the existence of global strong solution.
\section{Main results}
We are going now to state the main results of this paper. 
\begin{theorem} Let $N=2,3$, $\gamma\geq 2$, $\rho_0\in L^\infty(\R^N)$ a radial function and $v_0\in L^\infty(\R^N)$ which is radially symmetric. We define now the momentum velocity $m_0$ as follows:
\begin{equation}
\label{momentum1}
m_0=\rho_0 v_0+2\mu\n\rho_0.
\end{equation}
Then there exists a weak solution $(\rho,m,v)$ for the system (\ref{0.1m}) on a finite time interval $[0,T]$ with $T>0$ depending on $(\rho_0,v_0)$ provided that $v_0\in L^1(\R^N)$. In addition $(\rho,v,m)$ belongs to the space $Y_T$ defined above.\\
If in addition $\frac{1}{\rho_0}\in L^\infty(\R^N)$ then there exists a weak solution $(\rho,u)$ for the system (\ref{0.1}) on a finite time interval $[0,T]$ with $T>0$ depending on $(\rho_0,v_0)$. In addition $(\rho,v,m)$ belongs to the space $Y_T$ defined above and $\frac{1}{\rho}$ is in $L^\infty([0,T]\times\R^N)$. Furthermore $u\in L^2_{loc}([0,T]\times\R^N)$.\\
If in addition we assume that $N=2$, $v_0\in L^2(\R^2)$, $(\rho_0-1)\in L^2(\R^2)$ and 
$(u_0,v_0)$ are in $(L^1(\R^2))^{4}$ then there exists a global weak solution $(\rho,u)$ of (\ref{0.1}) for $\gamma\geq 2$. $(\rho,v,m)$ belongs to the space $Y_T$ with $T$ defined above and $\frac{1}{\rho}$ is in $L^\infty([0,T]\times\R^N)$. In addition we have for any $T'>0$:
$${\cal E}_1(\rho,v)(T')\leq C,$$
and for $C>0$ large enough.
\label{theo2}
\end{theorem}
\begin{remarka}
In the previous theorem, the regularity on the initial density $\rho_0$ seems to be minimal
, in particular our result allows to deal with initial density admitting shocks (what is not the case in the framework of strong solution except in \cite{DFP} under smallness assumption on the density and for constant viscosity coefficients). Furthermore if we compare with the results of global weak solution (see \cite{Li,V}), we do not ask any BD control on the initial density. Indeed it is not necessary to assume that $\n\sqrt{\rho_0}$ and $(\rho_0-1)$ are respectively in $L^2(\R^N)$ and $L^\gamma_2(\R^N)$ \footnote{We refer to \cite{Lio98} for the definition of the Orlicz spaces.}. Otherwise as we mentioned in the introduction, we can observe regularizing effects on the density even when the initial density admits shocks. Indeed since $\n\rho=\frac{1}{2\mu}(\rho v-m)$, it implies that $\sqrt{t} \n\rho(t,\cdot)$ belongs to $L^\infty(\R^N)$ for any $t\in(0,T)$ because $v$ is in $L^\infty([0,T]\times\R^N)$ and $m$ in ${\cal E}_T$. We deduce then that we have for $C>0$ large enough and any $t\in(0,T)$:
\begin{equation}
\|\rho(t,\cdot)\|_{W^{1,\infty}(\R^N)}\leq\frac{C}{\sqrt{t}}.
\end{equation}
It implies in particular that the density $\rho$ becomes instantaneously a Lipschitz function and in particular a continuous function.
\end{remarka}
\begin{remarka} It is important to mention that even if $\frac{1}{\rho_0}$ is not in $L^\infty(\R^N)$, we get a weak solution in finite time for the system (\ref{0.1m}) where the momentum $m$ is a unknown.\\
In opposite when $\frac{1}{\rho_0}$ belongs to $L^\infty(\R^N)$, we obtain a weak solution for the system  (\ref{0.1}) where the velocity $u$ is a unknown and is well defined on $(0,T)$. 
It is important to mention that in the first case we can deal with initial density admitting vacuum. In addition we are able to prove regularizing effects on the momentum since $m$ belongs to $\widetilde{L}^1_T(B^{1}_{\infty,\infty})$.
\end{remarka}
\begin{remarka}
It is also interesting to point out that we obtain weak solution in finite time without assuming that the initial density $\rho_0$ or the velocity $v_0$ are close from equilibrium when $|x|$ goes to $+\infty$.
Up to our knowledge, in order to use energy estimates all the results of weak solution concern initial density $\rho_0$ and initial velocity $u_0$ verifying roughly speaking:
$$\lim_{|x|\rightarrow +\infty}\rho_0(|x|)=\bar{\rho}\;\;\mbox{and}\;\;\lim_{|x|\rightarrow +\infty}u_0(|x|)=\bar{u},$$
with $\bar{\rho}\in\R^+$ and $\bar{u}\in\R^N$. It is also true when we consider strong solution (see \cite{Danchin}).
\end{remarka}
\begin{remarka}
We observe that if we assume in addition that $v_0$ and $(\rho_0-1)$ are in $L^2(\R^N)$ with $N=2$ then there exists a global weak solution. It implies in particular that the solution verifies the BD entropy but remains of infinite energy.
\end{remarka}
We have a similar result in dimension $N=1$ without any geometrical assumption.
\begin{theorem}
Let $N=1$, $\gamma\geq 2$, $\rho_0\in L^\infty(\R)$ and $v_0\in L^\infty(\R)$. We define now the momentum velocity $m_0$ as follows:
\begin{equation}
\label{momentum1}
m_0=\rho_0 v_0+2\mu\p_x\rho_0.
\end{equation}
Then there exists a weak solution $(\rho,m,v)$ for the system (\ref{0.1m}) on a finite time interval $[0,T]$ with $T>0$ depending on $(\rho_0,v_0)$. In addition $(\rho,v,m)$ belongs to the space $Y_T$ defined above.\\
If in addition $\frac{1}{\rho_0}\in L^\infty(\R)$ then there exists a weak solution $(\rho,u)$ for the system (\ref{0.1}) on a finite time interval $[0,T]$ with $T>0$ depending on $(\rho_0,v_0)$. In addition $(\rho,v,m)$ belongs to the space $Y_T$ defined above. If in addition we assume that $v_0\in L^2(\R)$ and $(\rho_0-1)\in L^2(\R)$ then there exists a global weak solution $(\rho,u)$ of (\ref{0.1}). \label{theo3}
\end{theorem}
\begin{remarka}
In the previous Theorem we improve the results of \cite{glob} inasmuch as we do not require that $\rho_0$ belongs to $BV(\R)$. In particular we can deal with initial density which have no limits in the left and on the right. In counterpart we need to ask more integrability on $v_0$ since in \cite{glob} $v_0$ is only in $L^2(\R)$. We can also mention that we do not need any smallness assumption on $u_0$ and $v_0$.
\end{remarka}
We finish by giving a result of global strong solution for the system (\ref{0.1}) when the initial velocity is radially symmetric and that the initial data are sufficiently regular.
\begin{theorem}
Let $N=2,3$. We assume that $(\rho_0-1,u_0)$ belongs to $(\dot{B}^{\NN}_{p,1}(\R^N)\cap \dot{B}^{\NN+\e}_{p,1}(\R^N))\times (\dot{B}^{\NN-1}_{p,1}(\R^N)\cap \dot{B}^{\NN-1+\e}_{p,1}(\R^N))$ with $\e>0$ and $\frac{N}{1-\e}<p<2N$. Furthermore we suppose that $0<c\leq \rho_0$, $\rho_0$ is a radial function, $u_0$ is radially symmetric and $v_0,u_0\in L^1(\R^N)\cap L^\infty(\R^N)$. We have also:
\begin{equation}
{\cal E}(\rho_0,u_0)<+\infty,\;\;{\cal E}_1(\rho_0,v_0)<+\infty,
\label{energietheo}
\end{equation}
then:
\begin{itemize}
\item if  $N=2$ and $\gamma=1$ or $\gamma\geq 2$ there exists a unique global solution $(\rho,u)$ of the system (\ref{0.1}) .
\item if $N=3$ and $\gamma=1$ there exists a unique global solution $(\rho,u)$ of the system (\ref{0.1}).
\end{itemize}
The solution $(\rho,u,v)$ verifies in addition:
\begin{equation}
\begin{cases}
\begin{aligned}
&\rho,\frac{1}{\rho}\in L^\infty_{loc}(\R^+,L^\infty(\R^N))\\
&u\in\widetilde{L}^\infty_{loc}(\R^+,\dot{B}^{\NN-1}_{p,1}\cap \dot{B}^{\NN-1+\e}_{p,1})\cap \widetilde{L}^1_{loc}(\R^+,\dot{B}^{\NN+1}_{p,1}\cap \dot{B}^{\NN+1+\e}_{p,1})\\
&v\in  L^\infty_{loc}(\R^+,L^\infty(\R^N)).
\end{aligned}
\end{cases}
\end{equation}
\label{theo4}
\end{theorem}
Section \ref{section3} deals with the proof of  the Theorem \ref{theo4}. In the section \ref{section4} and \ref{section5} we prove the Theorems \ref{theo2} and \ref{theo3}.
\section{Proof of the Theorem \ref{theo4}}
\label{section3}
We wish to prove the existence of global strong solution for large initial data provided that the initial density $\rho_0$ is radial and the initial velocity $u_0$ is radially symmetric ( in addition $(\rho_0,u_0)$ satisfy the regularity assumptions of the theorem \ref{theo4}). The proof is divided in different steps, in a first step we recall some generic results of existence of strong solution in finite time on a maximal time interval $(0,T^*)$. Next we prove that the solution remains radially symmetric on the time interval $(0,T^*)$. The last steps consist in proving that necessary $T^*=+\infty$, to do this it suffices to prove that:
$$\int^{T^*}_{0}\|\n u(t,\cdot)\|_{L^\infty} dt<+\infty.$$
In order to show the previous estimate, the main difficulty is to get the following estimate:
\begin{equation}
\|(\rho,\frac{1}{\rho})\|_{L^\infty((0,T^*),L^\infty(\R^N))}<+\infty.
\label{sudcru}
\end{equation}
To do this, we will observe that the $L^\infty$ norm of the effective velocity $v$ is bounded on $[0,T^*]\times\R^N$. Using the maximum principle for the mass equation of (\ref{0.3}) we can establish (\ref{sudcru}).
\subsection*{Existence of strong solution in finite time on a maximal time interval $(0,T^*)$}
Let us recall a result of existence of strong solution for the system (\ref{0.1}) in critical Besov spaces (we refer to \cite{Danchin} for the definition of the Besov spaces). In addition, we will give a blowup criterion which enables us to prove that the solution exists globally for initial data verifying the Theorem \ref{theo4}. In the sequel we will note $q=\rho-1$ and $q_0=\rho_0-1$.
\begin{theorem}
Let $N\geq 2$, the viscosity coefficients verify (\ref{condition}). Assume now that $(q_0,u_0)\in \dot{B}^{\NN}_{p,1}\times \dot{B}^{\NN-1}_{p,1}$ with $1\leq p<2N$ and that there exists $c$ such that $0<c\leq\rho_0$. Then there exists a time $T$ such that system (\ref{0.1}) has a unique solution on $[0,T]$ with:
\begin{equation}
q\in \widetilde{C}_{T}(\dot{B}^{\NN}_{p,1}),\;(\frac{1}{\rho},\rho)\in (L^{\infty}_{T}(L^\infty(\R^N)))^2\;\;\;\mbox{and}\;\;u\in \widetilde{C}_{T}(\dot{B}^{\NN-1}_{p,1})\cap L^1_T(\dot{B}^{\NN+1}_{p,1}).
\label{hypimportant}
\end{equation}
If in addition $(q_0,u_0)$ belongs to $\dot{B}^{s}_{p_1,1}\times \dot{B}^{s-1}_{p_1,1}$ for any $s\geq \frac{N}{p_1}$ with $1\leq p_1<2N$
then we have:
\begin{equation}
q\in \widetilde{C}_{T}(\dot{B}^{s}_{p_1,1})\;\;\mbox{and}\;\;u\in \widetilde{C}_{T}(\dot{B}^{s-1}_{p_1,1})\cap L^1_T
(\dot{B}^{s+1}_{p_1,1}).
\label{estimcru33}
\end{equation}
Assume now that (\ref{0.1}) has a strong solution $(q,u)\in C([0,T^*),\dot{B}^{\NN}_{p,1}\times \dot{B}^{\NN-1}_{p,1})$ on the time interval $[0,T^*)$ (with $T^*>0$) which satisfies the following three conditions:
\begin{enumerate}
\item the function $q$ belongs to $L^\infty([0,T^*],\dot{B}^{\NN}_{p,1})$,
\item $\frac{1}{\rho}\in L^\infty([0,T^*],L^\infty(\R^N))$,
\item we have $\int^{T^*}_{0}\|\n u(s)\|_{L^\infty} ds<+\infty$,
\end{enumerate}
Then $(q,u)$ may be continued beyond $T^*$.
\label{theo1a}
\end{theorem}
The theorem \ref{theo1} is a direct consequence of \cite{JDE,M3AS,DG,Fourier}. In \cite{JDE} the existence is proved for $(q_0,u_0)\in B^{\NN}_{p,1}\times B^{\NN-1}_{p,1}$ for $1\leq p<2N$ with uniqueness for $1\leq p\leq N$. The uniqueness is extended in \cite{Fourier} to the case $N\leq p<2N$. The blow up criteria is given in \cite{JDE} following an argument developed in \cite{DG}.
We can simply recall that the condition $p<2N$ is a restriction due to the paraproduct law , indeed we need to define properly the term of the form $\n\ln\rho\cdot Du$. For the sequel we recall the useful proposition (see \cite{Danchin}).
\begin{proposition}
\label{chaleur} Let $s\in\R$, $(p,r)\in[1,+\infty]^{2}$ and
$1\leq\rho_{2}\leq\rho_{1}\leq+\infty$. Assume that $u_{0}\in \dot{B}^{s}_{p,r}$ and $f\in\widetilde{L}^{\rho_{2}}_{T}
(\dot{B}^{s-2+2/\rho_{2}}_{p,r})$.
Let u be a solution of the heat equation:
$$
\begin{aligned}
&\p_{t}u-\mu\D u-(\lambda+\mu)\n{\rm div}u=f,\;u_{t=0}=u_{0}.
\end{aligned}
$$
with $\mu>0$, $\lambda+2\mu>0$.
Then there exists $C>0$ depending only on $N,\mu,\rho_{1}$ and
$\rho_{2}$ such that:
$$\|u\|_{\widetilde{L}^{\rho_{1}}_{T}(\dot{B}^{s+2/\rho_{1}}_{p,r})}\leq C\big(
 \|u_{0}\|_{\dot{B}^{s}_{p,r}}+\|f\|_{\widetilde{L}^{\rho_{2}}_{T}
 (\dot{B}^{s-2+2/\rho_{2}}_{p,r})}\big)\,.$$
 If in addition $r$ is finite then $u$ belongs to $\widetilde{C}([0,T],\dot{B}^{s}_{p,r})$. Similarly we have:
 $$\|u\|_{\widetilde{L}^{\rho_{1}}_{T}(B^{s+2/\rho_{1}}_{p,r})}\leq C(T)\big(
 \|u_{0}\|_{B^{s}_{p,r}}+\|f\|_{\widetilde{L}^{\rho_{2}}_{T}
 (B^{s-2+2/\rho_{2}}_{p,r})}\big)\,,$$
 with $C$ a continuous increasing function on $\R^+$. If in addition $r$ is finite then $u$ belongs to $\widetilde{C}([0,T],B^{s}_{p,r})$.
\end{proposition}
\begin{remarka}
We refer to \cite{Danchin} for the definition of the Besov spaces (homogeneous and non homogeneous), the Littlewood-Paley decomposition and the paraproduct laws. The spaces $\widetilde{L}^\rho(\dot{B}^{s}_{p,r})$,  $\widetilde{L}^\rho(B^{s}_{p,r})$ are the Chemin-Lerner spaces and are also defined in \cite{Danchin}. We will use in the sequel the same notations as \cite{Danchin}. 
\end{remarka}
\subsection*{Existence of strong solution in finite time which remains radially symmetric}
Let $(\rho_0,u_0,v_0)$ verifying the assumptions of the Theorem \ref{theo2}, 
we now smooth out the data as follows with $b_0=\rho_0-1$ :
\begin{equation}
(b_0)_{n}=S_{n}b_{0}\;\;(u_{0})_{n}=S_{n}u_0.
\label{minitial}
\end{equation}
We recall that for $w$ a temperated distribution, we have $S_nw=\chi(\frac{D}{2^n})w$ with $\chi\in C^\infty_0(\R^N)$, $\mbox{supp}\chi\subset B(0,\frac{4}{3})$, $\chi$ is equal to $1$ in a neighborhood of $0$ with values in $[0,1]$ and $\chi(\xi)$ radial.
We observe then that $(b_0)_n$ is always radial since the kernel associated to $S_n$ is radial and $b_0$ is radial. Indeed we recall that the convolution of two radial functions is radial. Let us show now that $(u_0)_n$ is radially symmetric for any $n\in\mathbb{N}$. Let us prove this claim, since $\hat{ u_{0}}$ is in $L^1_{loc}$ we deduce that $\hat{S_n u_0}$ belongs to $L^1$ and we can use the inverse Fourier theorem which yields (we have used the fact that $u_0\in L^1(\R^N)$):
$$
\begin{aligned}
& S_n u_{0}(x)=\int_{\R^N} e^{ix\cdot \xi} \hat{S_n u_0}(\xi) d\xi\\
\end{aligned}
$$
We recall that for $w$ a temperated distribution $S_n w=\chi(\frac{D}{2^n})w$ for all $n\in\mathbb{N}$. 
We set $(u_{0,k})(x)=(1-\chi(kx)) u_0 (x)$ and we denote in the sequel by $e_{r,x}=\frac{x}{|x|}$ for $x\in\R^N\backslash\{0\}$. Next we have since 
$(u_{0,k})(x)= x\, \frac{u_{0,k}(x)\cdot x}{|x|^2} $:
$$
\begin{aligned}
&\hat{S_n u_{0,k}}(\xi) =\chi(\frac{\xi}{2^n})\int_{\R^N} e^{-ix\cdot \xi}   x ( \frac{u_{0,k}(x)\cdot x}{|x|^2})d x\\
&=i \chi(\frac{\xi}{2^n})\n_\xi \int_{\R^N} e^{-ix\cdot \xi}  \,(\frac{u_{0,k}(x)\cdot x}{|x|^2} ) d x
\end{aligned}
$$
Since $\int_{\R^N} e^{-ix\cdot \xi}  \,(\frac{u_{0,k}(x)\cdot x}{|x|^2} ) d x$ is radial (indeed the Fourier transform of a radial function remains radial), we have if we denote $r=|\xi|$ and since $u_0\in L^1(\R^N)$.:
$$
\begin{aligned}
& \hat{S_n (u_{0,k})}(\xi)=i \chi(\frac{\xi}{2^n})\n_\xi \int_{\R^N} e^{-ix_1 |\xi|}  \,(\frac{u_{0,k} (x)\cdot x}{|x|^2} ) d x\\
&=i \chi(\frac{\xi}{2^n})\frac{\xi}{|\xi|}\p_r \int_{\R^N} e^{-ix_1 r}  \,(\frac{u_{0,k} (x)\cdot x}{|x|^2} ) d x\\
&=\chi(\frac{\xi}{2^n})\frac{\xi}{|\xi|} \int_{\R^N} x_1 e^{-ix_1 |\xi|}  \,(\frac{u_{0,k} (x)\cdot x}{|x|^2} ) d x\\
&=\chi(\frac{\xi}{2^n})\frac{\xi}{|\xi|}A_{k,n}(|\xi|),
\end{aligned}
$$
with:
$$A_{k,n}(|\xi|)= \int_{\R^N} x_1 e^{-ix_1 |\xi|}  \,(\frac{u_{0,k} (x)\cdot x}{|x|^2} ) d x.$$
Here $A_{k,n}$ belongs to $L^\infty$ since $u_0\in L^1(\R^N)$. 
Next we have if we define $ \hat{S_{n,k} u_{0,k}} (\xi)= \chi(\frac{\xi}{2^n})(1-\chi(k\xi))\hat{u_{0,k}} (\xi)=\chi_{k,n}(\xi)\hat{u_{0,k}} (\xi)$ , we deduce via the inverse Fourier transform that:
$$
\begin{aligned}
 S_{n,k} u_{0,k}(x)=&\int_{\R^N} e^{ix\cdot \xi} \chi_{k,n}(\xi) \frac{\xi}{|\xi|}A_{k,n}(|\xi|) d\xi\\
&=-i \n_x \int_{\R^N} e^{ix\cdot \xi}  \frac{1}{|\xi|}\chi_{k,n}(\xi) A_{k,n}(|\xi|) d\xi\\
&=-i \n_x \int_{\R^N} e^{i\xi_1 |x|}  \frac{1}{|\xi|}\chi_{k,n}(\xi)  A_{k,n}(|\xi|) d\xi\\
&=\frac{x}{|x|} \int_{\R^N} e^{i\xi_1 |x|}  \frac{\xi_1}{|\xi|}\chi_{k,n}(\xi)A_{k,n}(|\xi|) d\xi=B_{k,n}(x).
\end{aligned}
$$
Next we prove using  dominated convergence (because $|\chi_{k,n}|\leq |\chi_n|$, $|A_{k,n}|\leq |A_n|$ and $|u_{0,k}|\leq |u_0|$)  and the Plancherel theorem that:
\begin{equation}
\begin{aligned}
&S_{n,k} u_{0,k} \rightarrow_{k\rightarrow+\infty} S_n u_{0}\;\;\;\mbox{in}\;\;L^2\;\;\;\mbox{and}\;\;\;B_{k,n} \rightarrow_{k\rightarrow+\infty} B_n\;\;\;\mbox{a.e}.
\end{aligned}
\end{equation}
with:
$$
\begin{aligned}
&B_n(x)=\frac{x}{|x|} \int_{\R^N} e^{i\xi_1 |x|}  \frac{\xi_1}{|\xi|}\chi(\frac{\xi}{2n})A_n(|\xi|) d\xi\\
&A_n(|\xi|)= \int_{\R^N} x_1 e^{-ix_1 |\xi|}  \,(\frac{u_{0} (x)\cdot x}{|x|^2} ) d x.
\end{aligned}
$$
When we extract up to a subsequence in $k$ we get almost everywhere:
\begin{equation}
\begin{aligned}
& S_n u_{0} (x)=\frac{x}{|x|} \int_{\R^N} e^{i\xi_1 |x|}  \frac{\xi_1}{|\xi|}\chi(\frac{\xi}{2n})A_n(|\xi|) d\xi.
\end{aligned}
\end{equation}
In particular we deduce that $(u_0)_n$ conserves the rotational invariant structure $\frac{|\cdot|}{\cdot}P_n(\cdot)$ with $P_n$ radial. We can now observe that we have for $s\in\R$, $(p,r)\in[1,+\infty]^2$:
$$\forall l\in\mathbb{Z},\;\;\|\D_{l}(b_{0})_{n}\|_{L^{p}}\leq\|\D_{l}b_{0}\|_{L^{p}}\;\;\;\mbox{and}\;\;\;\|(b_{0})_{n}\|
_{\dot{B}^s_{p,r}}\leq \|b_{0}\|_{\dot{B}^s_{p,r}},$$
and similar properties for $(u_{0})_{n}$, a fact which will be used repeatedly during the next
steps. \\
We have seen from the Theorem \ref{theo1a} that there exists a sequence of strong solution $(\rho_n,u_n)_{n\in\mathbb{N}}$ on a maximal time interval $(0,T_n^*)$ of the system (\ref{0.1}) for the initial data $(1+(b_0)_n,(u_0)_n)_{n\in\mathbb{N}}$. It is important to point out that these solutions $(\rho_n,u_n)_{n\in\mathbb{N}}$ satisfy (\ref{hypimportant}), (\ref{estimcru33}) (with $p>N$ defined in the Theorem \ref{theo4}) and are in $C^\infty((0,T_n^*)\times\R^N)$, this is due to the fact that any $\dot{B}^s_{p,1}$ norm is propagated with $s>0$ large enough and we use in addition the Besov embeddings. In addition, $(\rho_n,u_n)_{n\in\mathbb{N}}$ verify for any $0<t<T_n^*$ the energy estimates (\ref{cenergie}) and (\ref{bcenergie}).
We are now interested in proving that the density $\rho_n$ remains radial on the time interval $[0,T_n^*)$ and similarly that the velocity $u_n$ is radially symmetric on  the interval $[0,T_n^*)$.
\begin{remarka}
 In the sequel, for simplicity in the notations we forget the subscript $n$. However we will need to use compactness arguments to prove that up to a subsequence $(\rho_n,u_n)_{n\in\mathbb{N}}$ converges to a global strong solution $(\rho,u)$, at this moment of the proof we will use again the subscript $n$.
 \end{remarka}
In other words we wish to prove that $(\rho,u)$ can be written under the following form for any $t\in(0,T^*)$ and $x\in\R^N\setminus \{0\}$:
\begin{equation}
\rho(t,x)=\rho_1(t,|x|)\;\;\mbox{and}\;\;u(t,x)=\frac{x}{|x|}u_1(t,|x|),
\label{invariance2}
\end{equation}
with $\rho_1$ and $u_1$ some functions with values in $\R$. We are going to deal with the case $N=3$. The case $N=2$ can be treated in a similar way.
In a first time, we want to check that for  any orthogonal matrix $A\in{\cal M}_{N}(\R^N)$ ( $A^t A=A A^t=Id_N$) we have for all $t\in(0,T^*)$:
\begin{equation}
\rho(t,x)=\rho(t,^t A x)\;\;\mbox{and}\;\;u(t,x)=Au(t,^t A x).
\label{invariance}
\end{equation}
We can observe since $(\rho,u)$ is a regular solution on $(0,T^*)$ that we have:
\begin{equation}
\begin{cases}
\begin{aligned}
&\p_t (\rho(t,^t A x))+{\rm div}\big( \rho(t,^t A x)A u (\rho(t,^t A x)\big)=0\\
&\rho(t,^t A x)\p_t (A u (t,^t A x))+\rho (t,^t A x)A u(t,^t A x)\cdot\n (A u (t,^t A x))\\
&\hspace{0,5cm}-\mu \D ( A u(t,^t A x))-\mu\n{\rm div}( A u(t,^t A x))  +\n[P(\rho)(t,^t A x)]=0. 
\end{aligned}
\end{cases}
\label{symcru}
\end{equation}
It implies that $(\rho_1(t,\cdot),u_1(t,\cdot))=(\rho(t,^t A \cdot),Au(t,^t A \cdot))$ is also a solution of the system (\ref{0.1}) on $(0,T^*)$ verifying the regularity assumption (\ref{hypimportant}). Indeed we have:
$$
\begin{aligned}
\D_l(\rho(t,^tA\cdot)-1)(x)&=2^{lN}\int_{\R^N}h(2^l y)(\rho(t,^t A(x-y))-1) dy\\
&=2^{lN}\int_{\R^N}h(2^l A(^t Ax-u))(\rho(t,u)-1) du.
\end{aligned}
$$
Now since $h={\cal F}^{-1}\va$ with $\va$ radial, we deduce that $h$ is radial and we get:
$$
\begin{aligned}
\D_l(\rho(t,^tA\cdot)-1)(x)=\D_l(\rho(t,\cdot)-1)(^t Ax)
\end{aligned}
$$
We deduce now that we have for any $p\in[1,+\infty]$ and any $t\in(0,T^*)$:
$$\|\D_l \rho(t,^tA\cdot)-1)\|_{L^p}=\|\D_l(\rho(t,\cdot)-1)\|_{L^p}.$$
It implies that $\rho_1$ has the same regularity than $\rho$, we proceed similarly for $u$ and $u_1$.
\\
We can now verify that $(\rho_0(^t A\cdot), A u_0(^t A\cdot))=(\rho_0,u_0)$ is in $\dot{B}^\NN_{p,1}(\R^N)\times \dot{B}^{\NN-1}_{p,1}(\R^N)$. Since $(\rho_1,u_1)(0,\cdot)=(\rho_0,u_0)$ (because $\rho_0$ is radial and $u_0$ is radially symmetric) we deduce by uniqueness (see  the Theorem \ref{theo1a}) that $\rho_1=\rho$ and $u_1=u$ which proves (\ref{invariance}). It implies in particular that for any rotation matrix $R_{h,\theta}$ ( with  $\theta$ is the angle of rotation and the axis of rotation which is given by the unitary vector $h$), we have for any $t\in(0,T^*)$ and any $x\in\R^3$:
\begin{equation}
\rho(t,R_ {h,\theta}x)=\rho(t,x).
\label{densradial}
\end{equation}
We deduce that the density $\rho$ remains radial on the time interval $(0,T^*)$ and we have for any $(t,x)\in(0,T^*)\times\R^3$:
\begin{equation}
\rho(t,x)=\rho_2(t,|x|),
\label{densradial1}
\end{equation}
with $\rho_2$ a function. We are going to prove in a similar way that the velocity remains radially symmetric on the time interval $(0,T^*)$.
We are interested now in working in spherical coordinates, more precisely we define the following orthonormal basis $(e_{r,x},e_{\theta,x},e_{\phi,x})$ depending on $x\in\R^3\backslash \{x_1=x_2=0\}$ with $x=^t (x_1,x_2,x_3)$ and $x'=^t(x_1,x_2)$. We note:
\begin{equation}
\begin{aligned}
&e_{r,x}=\frac{x}{|x|},\\
&e_{\theta,x}=^t(\frac{x_1 x_3}{\sqrt{(x_1^2+x_2^2)^2+x_3^2 (x_1^2+x_2^2)}},\frac{x_2 x_3}{\sqrt{(x_1^2+x_2^2)^2+x_3^2 (x_1^2+x_2^2)}},-\frac{\sqrt{x_1^2+x_2^2}}{\sqrt{x_1^2+x_2^2+x^2_3}}),\\
& e_{\phi,x}=\frac{1}{\sqrt{x_1^2+x_2^2}}^t(-x_2,x_1,0).
\end{aligned}
\label{formu1}
\end{equation}
We have in particular $e_{\theta,x}=^t(\frac{x_1 x_3}{|x'|\,|x|},\frac{x_2 x_3}{|x'|\,|x|},-\frac{|x'|}{|x|})$ for  $x\in\R^3\backslash \{x_1=x_2=0\}$. We observe that:
\begin{equation}
\begin{cases}
\begin{aligned}
&|x'|e_{r,x}+x_3 e_{\theta,x}=\frac{|x|}{|x'|}\,^t(x',0)\\
&\frac{x_3}{|x|} e_{r,x}-\frac{|x'| }{|x|}e_{\theta,x}=^t(0,0,1).
\end{aligned}
\end{cases}
\label{algebre}
\end{equation}
We wish now to prove that we can write the velocity $u$ under the following form for any $t\in(0,T^*)$ and $x\in\R^3\backslash\{x_1=x_2=0\}$:
\begin{equation}
u(t,x)=u_{1,1}(t,|x|) e_{r,x}+u_{2,1}(t,|x'|,x_3)e_{\theta,x}+u_{3,1}(t,|x'|,x_3)e_{\phi,x}.
\label{decomporad}
\end{equation}
For the moment since $(e_{r,x},e_{\theta,x},e_{\phi,x})$ is a orthonormal basis for $x\in\R^3\backslash \{x_1=x_2=0\}$ with $x=^t (x_1,x_2,x_3)$, we have:
\begin{equation}
u(t,x)=u_1(t,x) e_{r,x}+u_2(t,x)e_{\theta,x}+u_3(t,x)e_{\phi,x}.
\label{decomporada}
\end{equation}
Since $u_{1}(t,x)=u(t,x)\cdot e_{r,x}$, $u_{2}(t,x)=u(t,x)\cdot e_{\theta,x}$, $u_{3}(t,x)=u(t,x)\cdot e_{\phi,x}$, we deduce that $u_1$, $u_2$ and $u_3$ are regular on $(0,T^*)\times(\R^3\backslash \{x_1=x_2=0\})$.\\
For any rotation $R_{h,\theta}$ we have from (\ref{invariance}) and for any $t\in(0,T^*)$ and $x\in\R^3$ with $(x, ^t R_{h,\theta}x)\in (\R^3\backslash \{x_1=x_2=0\})^2$:
\begin{equation}
\begin{aligned}
&u_1(t,^t R_{h,\theta}x) R_{h,\theta} e_{r,^tR_{h,\theta}x}+u_2(t,^t R_{h,\theta}x) R_{h,\theta}e_{\theta,^t R_{h,\theta}x}+u_3(t,^t R_{h,\theta}x)R_{h,\theta} e_{\phi,^t R_{h,\theta}x}\\
&=u_1(t,x) e_{r,x}+u_2(t,x)e_{\theta,x}+u_3(t,x)e_{\phi,x}.
\end{aligned}
\label{decomporadb}
\end{equation}
Now using (\ref{formu1}) and the fact that $R_{h,\theta}$ is an isometry, we deduce that we have:
\begin{equation}
\begin{aligned}
&u_1(t,^t R_{h,\theta}x)e_{r,x}+u_2(t,^t R_{h,\theta}x) R_{h,\theta}e_{\theta,^t R_{h,\theta}x}+u_3(t,^t R_{h,\theta}x)R_{h,\theta} e_{\phi,^t R_{h,\theta}x}\\
&=u_1(t,x) e_{r,x}+u_2(t,x)e_{\theta,x}+u_3(t,x)e_{\phi,x}.
\end{aligned}
\label{decomporad1}
\end{equation}
From (\ref{decomporad1}), taking the scalar product with $e_{r,x}$, we have since  $(e_{r,x},e_{\theta,x},e_{\phi,x})$ is a orthonormal basis: 
\begin{equation}
\begin{aligned}
&u_1(t,^t R_{h,\theta}x)+u_2(t,^t R_{h,\theta}x) \langle R_{h,\theta}e_{\theta,^t R_{h,\theta}x},e_{r,x}\rangle+u_3(t,^t R_{h,\theta}x)\langle R_{h,\theta} e_{\phi,^t R_{h,\theta}x},e_{r,x}\rangle\\
&\hspace{11cm}=u_1(t,x).
\end{aligned}
\label{decomporad12}
\end{equation}
We have now since $R_{h,\theta}$ is an isometry and $(e_{r,x},e_{\theta,x},e_{\phi,x})$ is a orthonormal basis for any $x$ such that $(x,^t R_{h,\theta}x) \in(\R^3\backslash \{x_1=x_2=0\})^2$ with $x=^t (x_1,x_2,x_3)$:
\begin{equation}
\begin{aligned}
\langle R_{h,\theta} e_{\phi,^t R_{h,\theta} x},\frac{x}{|x|}\rangle &=\langle  e_{\phi,^t R_{h,\theta} x},\frac{^t R_{h,\theta} x}{|^tR_{h,\theta} x|} \rangle = \langle  e_{\phi,^t R_{h,\theta} x},e_{r,^t R_{h,\theta} x}\rangle\\
&=0.\\[2mm]
\langle R_{h,\theta} e_{\theta,^t R_{h,\theta} x},\frac{x}{|x|}\rangle &=\langle  e_{\theta,^t R_{h,\theta} x},\frac{^t R_{h,\theta} x}{|^t R_{h,\theta} x|} \rangle= \langle  e_{\theta,^t R_{h,\theta} x},e_{r,^t R_{h,\theta} x}\rangle \\
&=0.
\end{aligned}
\label{scalnul}
\end{equation}
Combining (\ref{decomporad1}) and (\ref{scalnul}), we deduce that for any $x$ such that $(x,^t R_{h,\theta}x) \in(\R^3\backslash \{x_1=x_2=0\})^2$  and any $t\in(0,T^*)$ we have for any rotation matrix $R_{h,\theta}$:
\begin{equation}
u_1(t,^t R_{h,\theta}x)=u_1(t,x).
\label{radial2}
\end{equation}
We deduce that $u_1(t,\cdot)$ is radial for any $t\in(0,T^*)$ on $\R^3\backslash \{x_1=x_2=0\}$. It proves in particular that for any $x\in \R^3\backslash \{x_1=x_2=0\}$ we have:
\begin{equation}
u_1(t,x)=u_{1,1}(t,|x|).
\label{radialu1}
\end{equation}
Similarly if we choose $h=^t(0,0,1)$, the matrix associated to $R_{h,\theta}$ is:
$$
\begin{pmatrix} 
\cos\theta &- \sin\theta &0\\
\sin\theta & \cos\theta &0\\
0&0&1\\
\end{pmatrix}
$$
We define also $A_\theta$ as follows \footnote{We observe that it is a rotation in $\R^2$ of angle $\theta$.}:
$
A_\theta=\begin{pmatrix} 
\cos\theta & -\sin\theta \\
\sin\theta & \cos\theta \\
\end{pmatrix}
$.
Noting $(x')^\perp=^t(-x_2,x_1)$, we have then since $A_\theta$ is a rotation and commute with $B=\begin{pmatrix} 
0 & -1 \\
1 & 0 \\
\end{pmatrix}$ :
\begin{equation}
\begin{aligned}
&R_{h,\theta} e_{\phi,^t R_{h,\theta}x}=R_{h,\theta}\begin{pmatrix} 
\frac{(^t A_\theta x')^{\perp}}{|x'|}\\
0\\
\end{pmatrix}
=\begin{pmatrix} 
A_\theta & 0 \\
0 &1 \\
\end{pmatrix}
\begin{pmatrix} 
\frac{ B\cdot ^t A_\theta x'}{|x'|}\\
0\\
\end{pmatrix}=e_{\phi,x}.
\end{aligned}
\label{techint}
\end{equation}
Proceeding as previously we deduce that for $h=^t (0,0,1)$ we have for any $x$ such that $x \in \R^3\backslash \{x_1=x_2=0\}$
and any $\theta$:
\begin{equation}
u_3(t,^t A_{\theta}x',x_3)=u_3(t,x',x_3).
\label{radial3}
\end{equation}
Indeed we have used the fact that from (\ref{techint}), we have $^t R_{h,\theta}e_{\phi,x}=e_{\phi,^t R_{h,\theta}x}$.
In other words we deduce that for any $x \in \R^3\backslash \{x_1=x_2=0\}$, we have:
\begin{equation}
u_3(t,x',x_3)=u_{3,1}(t,|x'|,x_3).
\label{radial4}
\end{equation}
Taking again $h=^t(0,0,1)$ and using the fact that $e_{\theta,x}=\frac{x_3}{|x|}^t(\frac{x'}{|x'|},-\frac{|x'|}{x_3})$ for $x \in \R^3\backslash \{x_1=x_2=0\}$, we have for $x_3\ne 0$:
\begin{equation}
\begin{aligned}
&e_{\theta,^t R_{h,\theta}x}=\frac{x_3}{|x|}\begin{pmatrix} 
\frac{^t A_\theta x'}{|x'|}\\
-\frac{|x'|}{x_3}\\
\end{pmatrix}\\
&R_{h,\theta} e_{\theta,^t R_{h,\theta}x}=\frac{x_3}{|x|}\begin{pmatrix} 
A_\theta & 0 \\
0 &1 \\
\end{pmatrix}
\begin{pmatrix} 
\frac{^t A_\theta x'}{|x'|}\\
-\frac{|x'|}{x_3}\\
\end{pmatrix}=e_{\theta,x}.
\end{aligned}
\label{ecriture}
\end{equation}
When $x_3=0$, we have $e_{\theta,x}=^t(0,0,-1)$. We have then easily:
\begin{equation}
\begin{aligned}
&R_{h,\theta} e_{\theta,^t R_{h,\theta}x}=\begin{pmatrix} 
A_\theta & 0 \\
0 &1 \\
\end{pmatrix}
\begin{pmatrix} 
0\\
-1\\
\end{pmatrix}=e_{\theta,x}.
\end{aligned}
\label{ecriture1}
\end{equation}
Proceeding as previously and using (\ref{ecriture}) and (\ref{ecriture1}), we deduce that for any $(t,x)\in (0,T^*)\times(\R^3  \R^3\backslash \{x_1=x_2=0\})$ we have:
\begin{equation}
u_2(t,x',x_3)=u_{2,1}(t,|x'|,x_3).
\label{radial5}
\end{equation}
Finally (\ref{radialu1}), (\ref{radial4}) and (\ref{radial5}) prove (\ref{decomporad}). In the sequel we note  $u^1(t,x)=u_1(t,|x|) e_{r,x}$, $u^2(t,x)=u_2(t,|x'|,x_3)e_{\theta,x}$ and $u^3(t,x)=u_3(t,|x'|,x_3)e_{\phi,x}$ with $u_{1,1}=u_1$, $u_{2,1}=u_2$ and $u_{3,1}=u_3$. \\
We would like now to understand what is the Besov regularity of each component $u^1$, $u^2$ and $u^3$. Let us consider now the vectors $u_i'$ in $\R^2$ with $i\in\{1,2,3\}$ defined as follows:
$$u_i'=^t(u^i_1,u^i_2). $$
We observe then (with ${\rm curl}u_i'=\p_1 (u_i')_2-\p_2 (u_i')_1$) that for any $x\in\R^3\backslash\{x_1=x_2=0\}$:
\begin{equation}
\begin{aligned}
&{\rm div}_{x'}u_3'(x)=0,\;{\rm curl}_{x'}u_3(x)'={\rm curl}_{x'}u'(x)\\[1mm]
&\hspace{2cm}\;\;\mbox{and}\;\;{\rm curl}_{x'}(u_1'+u_2')(x)=0,\;{\rm div}_{x'}(u_1'+u_2')(x)={\rm div}_{x'}u'(x).
\end{aligned}
\label{tech12}
\end{equation}
We recall now that for $w$ a temperated distribution we have ${\cal F}(\D_l u)(\xi)=\va(\frac{\xi}{2^l}){\cal F}u(\xi)$ with $\xi\in\R^3$ and $l\in\mathbb{Z}$.  Compared with \cite{Danchin}, we choose now a particular Litllewood-Paley decomposition and we assume here that $\va(\xi)=\va_1(\xi')\va(\xi_3)$ with $\xi'=^t (\xi_1,\xi_2)$ and $\va_1$ radial. Similarly we define $\D_l'$ in $\R^2$ and $\D_l^3$ in $\R$ as follows ${\cal F}'(\D_l' v)(\xi')=\va_1(\frac{\xi'}{2^l}){\cal F}'v(\xi')$,  ${\cal F}^3(\D_l^3 w)(\xi_3)=\va_2(\frac{\xi_3}{2^l}){\cal F}w(\xi_3)$  with ${\cal F}'$, ${\cal F}^3$ respectively  the Fourier transform in $\R^2$ and $\R$.\\
From (\ref{tech12}), we deduce that:
$$
{\cal F}(\D_l u'_3)(\xi)= \va_2(\frac{\xi_3}{2^l})\va_1(\frac{\xi'}{2^l})A(\xi'){\cal F}u(\xi),$$
with $A$ the Fourier multiplier associated to $(\D)_{x'}^{-1}\n^\perp_{x'}{\rm curl}_{x'}$.
We deduce using the inverse Fourier transform that:
$$
\begin{aligned}
&(\D_l u'_3)(t,x)=\int_{\R^2} e^{ix'\cdot\xi'}\va_1(\frac{\xi'}{2^l})A(\xi') \int_{\R}e^{i x_3 \xi_3}
\va_2(\frac{\xi_3}{2^l}) \hat{u'}(t,\xi) d\xi_3 d\xi'
\end{aligned}
$$
Next we have:
$$
\begin{aligned}
 \int_{\R}e^{i x_3 \xi_3}
\va_2(\frac{\xi_3}{2^l}) \hat{u'}(t,\xi) d\xi_3 &=\int_{\R^2}e^{-i z'\cdot\xi'}\big(\int_{\R}e^{i x_3 \xi_3}\va_2(\frac{\xi_3}{2^l})(\int_{\R}e^{-i z_3\xi_3}u'(t,z',z_3) dz_3)d\xi_3\big) dz'\\
&=\int_{\R^2}e^{-i z'\cdot\xi'}\big(\int_{\R}e^{i x_3 \xi_3}\va_2(\frac{\xi_3}{2^l}){\cal F}^3u'(t,z',\cdot) (\xi_3) d\xi_3\big) dz'\\
&=\int_{\R^2}e^{-i z'\cdot\xi'}\D_l^3 u'(t,z',\cdot)[x_3]dz'\\
&={\cal F}'[\D_l^3[u(t,',^3)](x_3)](\xi').
\end{aligned}
$$
The notation $(',^3)$ means that the unknown for $\xi'$ is $'$ and for $x_3$ is $^3$. We obtain then that:
$$
\begin{aligned}
(\D_l u'_3)(t,x)&=\int_{\R^2} e^{ix'\cdot\xi'}\va_1(\frac{\xi'}{2^l})A(\xi') {\cal F}'[\D_l^3[u(t,',^3)](x_3)](\xi') d\xi'\\
&=A(D')\D_l'[\D_l^3[u(t,',^3)](x_3)](x')
\end{aligned}
$$
We deduce then that for any $l\in\mathbb{Z}$, $p\in(1,+\infty)$ and any $x_3\in\R$, we have:
\begin{equation}
\begin{aligned}
&\|\D_l u'_3(t,\cdot,x_3)\|_{L^p(\R^2)}\leq C_p \|\D_l u'(t,\cdot,x_3)\|_{L^p(\R^2)}
\end{aligned}
\end{equation}
Integrating in $\R^3$ we obtain that:
\begin{equation}
\begin{aligned}
&\|\D_l u'_3(t,\cdot)\|_{L^p(\R^3)}\leq C_p \|\D_l u'(t,\cdot)\|_{L^p(\R^3)}
\end{aligned}
\end{equation}
We deduce in particular that $u'_3$ and then $u^3$ has the same regularity as $u$ in terms of Besov spaces $\dot{B}^s_{p,r}$ \footnote{It is important to mention that Besov space does not depend on our particular choice of $\va$ in the Littlewood-Paley decomposition.} with $p\in(1,+\infty)$. The same result is also true for $u'_1+u'_2$ and then $u^1+u^2$. 
We are now going to prove that for any $t\in(0,T^*)$ and any $x\in\R^3$ we have:
\begin{equation}
u(t,x)=u_1(t,|x|)\frac{x}{|x|}=u^1(t,x).
\label{decompfinal}
\end{equation}
We deduce then 
From Appendix, (\ref{densradial1}), (\ref{decomporad}), (\ref{calcu1}), (\ref{Lapla}), (\ref{Lapla1}) and (\ref{geometrie5}) we can rewrite the system (\ref{0.1}) as follows for any $(t,x)\in(0,T^*)\times (\R^3\backslash \{x_1=x_2=0\})$:
$$
\begin{aligned}
&\rho \p_t u^3+\rho (u\cdot\n u^3+u^3\cdot\n u^1+u^3\cdot\n u^2)-2\mu {\rm div}(\rho D u^3)+B_3(t,x) e_{r,x}+B_4(t,x)e_{\theta,x}=0,
\end{aligned}
$$
with $B_3$ and $B_4$ regular terms on $\R^3\backslash\{x_1=x_2=0\}$. 
Taking now the scalar product of the previous equation with $u^3$ and integrating over $(0,t)\times\R^3$ with $t\in(0,T^*)$:
\begin{equation}
\begin{aligned}
&\int^t_0\int_{\R^N}\rho(t,x)\p_t u^3\cdot u^3 (s,x) ds dx+\int^t_0\int_{\R^N}\rho(t,x) (u\cdot \n u^3)\cdot u^3 (s,x) ds dx\\
&+\int^t_0\int_{\R^N}(\rho( u_3\cdot\n u)\cdot u_3-\rho (u_3\cdot \n u^3)\cdot u_3)(s,x) ds dx\leq 0.
\end{aligned}
\end{equation}
It is important to point out that the previous integrals have a sense since $u^3$ has the same regularity as $u$. Indeed from (\ref{cenergie}), we know that $u$ is in $L^\infty((0,T^*),L^2)$, this is then also true for $u^3$. In addition $\n u$, $\n u^3$ are in $L^\infty((0,T^*),L^\infty)$.
We can observe now that $\p_t u^3=\p_t u\cdot e_{\phi,x}$ and then $\p_t u^3$ belongs to $L^\infty((0,T^*),L^2)$ since $\p_t u$ belongs to $L^\infty((0,T^*),L^2)$ (it suffices to use the fact that the regularity $\dot{B}^s_{2,1}$ is conserved on $(0,T^*)$ for $u$ with $s\geq 0$) . After integration by parts, we obtain:
\begin{equation}
\begin{aligned}
&\frac{1}{2}\int_{\R^N}\rho(t,x)|u^3|^2(t,x) dx+\int^t_0\int_{\R^N}\rho( u_3\cdot\n u)\cdot u_3 ds dx\\
&+\frac{1}{2}\int^t_0\int_{\R^N}\n\rho\cdot u^3 |u^3|^2 ds dx+\frac{1}{2}\int^t_0\int_{\R^N}\rho\,{\rm div}u^3 |u^3|^2 ds dx \leq 0.
\end{aligned}
\label{dperimp}
\end{equation}
We observe now using (\ref{densradial1}) that $\n\rho(s,x)=\p_{r}\rho_1(t,|x|) e_{r,x}$ for any $(s,x)\in(0,T^*)\times(\R^3\backslash\{0\})$, it implies in particular using (\ref{decomporad}) that $\n\rho\cdot u^3=0$ almost everywhere on $(0,T^*)\times\R^3$. Using again (\ref{decomporad}), we can verify that ${\rm div}u^3=0$ on $(0,T^*)\times(\R^3\backslash\{x_1=x_2=0\})$. It implies from (\ref{dperimp}) that we have for any $t\in(0,T^*)$:
\begin{equation}
\begin{aligned}
&\frac{1}{2}\int_{\R^N}\rho(t,x)|u^3|^2(t,x) dx+\int^t_0\int_{\R^N}\rho( u_3\cdot\n u)\cdot u_3 ds dx\leq 0.
\end{aligned}
\label{dperimp1}
\end{equation}
Using the Gronwall lemma (this is possible because the term $\n u$ is in $L^1((0,T^*),L^\infty(\R^N))$ indeed $L^1((0,T^*),\dot{B}^{\NN}_{p,1}(\R^N))$ is embedded in $L^1((0,T^*),L^\infty(\R^N)$). We deduce that for any $t\in(0,T^*)$ we have $u^3(t,\cdot)=0$ on $\R^3$.\\
Similarly from Appendix, (\ref{densradial1}), (\ref{decomporad}), (\ref{Lapla}), (\ref{Lapla1}), (\ref{geometrie5}) and the fact that $u^3=0$, we can rewrite the system (\ref{0.1}) as follows for any $(t,x)\in(0,T^*)\times (\R^3\backslash \{x_1=x_2=0\})$:
$$
\begin{aligned}
&\rho \p_t u^2+\rho (u\cdot\n u^2+u^2\cdot\n u^1)-2\mu {\rm div}(\rho D u^2)+B_5(t,x) e_{r,x}=0,
\end{aligned}
$$
with $B_5$ a regular term on $(0,T^*)\times (\R^3\backslash \{x_1=x_2=0\})$. 
Taking now the scalar product of the previous equation with $u^2$ and integrating over $(0,t)\times\R^3$ with $t\in(0,T^*)$, we have:
\begin{equation}
\begin{aligned}
&\int^t_0\int_{\R^N}\rho(\p_t u^2\cdot u^2+( u\cdot\n u^2)\cdot u^2)(s,x)ds dx+\int^t_0\int_{\R^N}\rho (u^2\cdot\n u^1)\cdot u^2  (s,x) ds dx\\
&-2\mu \int^t_0\int_{\R^N} {\rm div}(\rho D u^2)\cdot u^2(s,x) ds dx=0.
\end{aligned}
\label{eu2}
\end{equation}
We are briefly justifying the sense of the previous integrals. First since $u^1(s,x)=u(t,x)\cdot e_{r,x}\,e_{r,x}$, we can verify that it exists $C>0$ large enough such that for all $x\in\R^3\backslash\{x_1=x_2=0\}$ and $s\in(0,t]$ we have:
\begin{equation}\begin{aligned}
&|\n u^1(s,x)|\leq C(\frac{|u(s,x)|}{|x|}+|\n u(s,x)|)\\
&|\n^2 u^1(s,x)|\leq C(\frac{|u(s,x)|}{|x|^2}+\frac{|\n u (s,x)|}{|x|}+|\n^2 u(s,x)|).
\end{aligned}
\label{ertimtech}
\end{equation}
 It implies that we have for $C>0$ large enough and since $\frac{1}{|\cdot|}$ and $\frac{1}{|\cdot|}^2$ are respectively in $L^2(B(0,1))\cap L^\infty(\R^3\backslash B(0,1))$ and $L^{\frac{3}{2}-\e}(B(0,1))\cap (L^\infty\cap L^2 )(\R^3\backslash B(0,1))$ with $\e>0$ small enough:
\begin{equation}
\begin{aligned}
&\|\n u^1(s,\cdot)\|_{L^2}\leq C(\|\n u(s,\cdot)\|_{L^2}+\| u(s,\cdot)\|_{L^\infty}+\|u(s,\cdot)\|_{L^2})\\
&\|\n^2 u^1(s,\cdot)\|_{L^{\frac{3}{2}-\e}}\leq C(\| u(s,\cdot)\|_{L^\infty}+\|u(s,\cdot)\|_{L^2}+\|\n u(s,\cdot)\|_{L^\infty}\\
&\hspace{7cm}+\|\n u(s,\cdot)\|_{L^2}+\|\n^2 u(s,\cdot)\|_{L^{\frac{3}{2}-\e}}).
\end{aligned}
\label{ertimtech1}
\end{equation}
Now using the fact that $u, \n u$ are respectively in $L^\infty([0,t],L^2\cap L^\infty)$, $\n u\in L^\infty([0,t],L^2\cap L^\infty)$ and $\n^2 u\in L^\infty(|0,t],L^{\frac{3}{2}})$ (we use here in particular the fact that the norm $\dot{B}^2_{\frac{3}{2},1}$ is conserved on $(0,T^*)$), we deduce that
$\n u^1$ and $\n^2 u^1$ are respectively bounded in $L^\infty([0,t],L^2)$ and $L^\infty([0,t],L^{\frac{3}{2}-\e})$ for $\e>0$ small enough. Similarly $\n u$ and $\n^2 u$ have the same regularity. Since $u^1+u^2$ has the same regularity than $u$ we deduce then that $\n u^2$ and $\n^2 u^2$ are respectively bounded in $L^\infty([0,t],L^2)$ and $L^\infty([0,t],L^{\frac{3}{2}-\e})$ for $\e>0$ small enough. Using now the fact that $u^2$ is easily bounded in $L^\infty([0,t],L^2\cap L^\infty)$, we deduce that we can bound each integral in (\ref{eu2}).
\\
Using now integrations by parts,
we have since $u^2(0,\cdot)=0$ (indeed $u_0$ is radially symmetric):
\begin{equation}
\begin{aligned}
&\frac{1}{2}\int_{\R^N}\rho(t,x)|u^2|^2(t,x) dx+2\mu  \int^t_0\int_{\R^N}\rho |Du^2|(s,x) ds dx\\
&\hspace{2cm}+
 \int^t_0\int_{\R^N}\rho (u^2\cdot\n u^1)\cdot u^2  (s,x) ds dx\leq 0.
\end{aligned}
\end{equation}
Since $\frac{1}{\rho}\in L^\infty([0,t]\times \R^N)$, we deduce that there exists $c_t>0$ such that:
\begin{equation}
\begin{aligned}
&\frac{1}{2}\int_{\R^N}\rho(t,x)|u^2|^2(t,x) dx+2\mu c_t \int^t_0\int_{\R^N} |Du^2|(s,x) ds dx\\
&\hspace{2cm}+
 \int^t_0\int_{\R^N}\rho (u^2\cdot\n u^1)\cdot u^2  (s,x) ds dx\leq 0.
\end{aligned}
\end{equation}
Since $ \int^t_0\int_{\R^N} |Du^2|(s,x) ds dx=\frac{1}{2} \int^t_0\int_{\R^N} (|\n u^2|+({\rm div}u^2)^2)(s,x) ds dx$, we get then:
\begin{equation}
\begin{aligned}
&\frac{1}{2}\int_{\R^N}\rho(t,x)|u^2|^2(t,x) dx+\mu c_t \int^t_0\int_{\R^N} |\n u^2|(s,x) ds dx\\
&\hspace{2cm}+
 \int^t_0\int_{\R^N}\rho (u^2\cdot\n u^1)\cdot u^2  (s,x) ds dx\leq 0.
\end{aligned}
\label{mro}
\end{equation}
Let us bound now the last term on the left hand side, it yields by Gagliardo-Niremberg inequality and Young inequality that for $\e>0$ and $C_\e>0$ large enough:
\begin{equation}
\begin{aligned}
&|\int^t_0\int_{\R^N}\rho (u^2\cdot\n u^1)\cdot u^2  (s,x) ds dx|\leq \int^t_0\|\rho(s,\cdot)\|_{L^\infty}
\|u_2(s,\cdot)\|_{L^4}^2\|\n u^1(s,\cdot)\|_{L^2} ds\\
&\leq \int^t_0\|\rho(s,\cdot)\|_{L^\infty}
\|u_2(s,\cdot)\|_{L^2}^{\frac{1}{2}}\|\n u^2(s,\cdot)\|_{L^2}^{\frac{3}{2}}\|\n u^1(s,\cdot)\|_{L^2} ds\\
&\leq \int^t_0(\e \|\n u^2(s,\cdot)\|_{L^2}+C_\e \|u^2(s,\cdot)\|_{L^2}^2 \|\rho(s,\cdot)\|^4_{L^\infty}
\|\n u^1(s,\cdot)\|^4_{L^2} ) ds
\end{aligned}
\label{mro1}
\end{equation}
Plugging (\ref{mro1}) in (\ref{mro}), we obtain for $\e>0$ sufficiently small:
\begin{equation}
\begin{aligned}
&\frac{1}{2}\int_{\R^N}\rho(t,x)|u^2|^2(t,x) dx+\frac{\mu c_t}{2} \int^t_0\int_{\R^N} |\n u^2|(s,x) ds dx\\
&\leq C_\e \int^t_0\|\frac{1}{\rho}(s,\cdot)\|_{L^\infty}\|\sqrt{\rho}u^2(s,\cdot)\|_{L^2}^2 \|\rho(s,\cdot)\|^4_{L^\infty}
\|\n u^1(s,\cdot)\|^4_{L^2}  ds.
\end{aligned}
\label{mro2}
\end{equation}
Now applying Gronwall lemma we deduce that $u^2(t,\cdot)=0$ almost everywhere. Indeed we have seen that $\|\n u^1(\cdot)\|_{L^2}^4$ is in $L^1(0,t)$.  
We have then proved that our sequence of solutions $u_n$ satisfied on $(0,T_n^*)\times\R^3$ $u_n=u_n^1$ and then remains radially symmetric all along the time interval $(0,T_n^*)$.
\subsection*{Gain of integrability on the velocity $u$ and $v$ when $\gamma=1$ and $N=2,3$}
We forget again the subscript $n$ in the next sections. We have seen that there exists a strong solution $(\rho,u)$ on the maximal time interval $(0,T^*)$, we wish now to prove that $T^*=+\infty$. In addition we have seen that for any $t\in(0,T^*)$ and any $x\in\R^N\backslash\{0\}$ we have:
\begin{equation}
\rho(t,x)=\rho_1(t,|x|)\;\;\;\mbox{and}\;\;\;u(t,x)=u_1(t,|x|)\frac{x}{|x|}=\n\theta(t,|x|).
\end{equation}
It implies in particular that $(\n u-^t\n u)(t,x)=0$ for any $t\in(0,T^*)$ and any $x\in\R^N\backslash\{0\}$. In the case when $P(\rho)=a\rho$, we recall that $u$ and $v$ verify since $\n u=^t\n u$ on $(0,T^*)$:
$$
\begin{cases}
\begin{aligned}
&\rho\p_t u+\rho u\cdot\n u-2\mu {\rm div}(\rho\n u)+\frac{a}{2\mu}\rho(v-u)=0,\\
&\rho\p_t v+\rho u\cdot\n v+\frac{a}{2\mu}\rho(v-u)=0.
\end{aligned}
\end{cases}
$$
Multiplying the first equation by $u|u|^p$ and the second one by $v|v|^p$ and integrating over $(0,t)\times\R^N$ with $t\in(0,T^*)$ (we recall that this procedure is possible since $(\rho,u,v)$ is sufficiently integrable and regular), we obtain:
\begin{equation}
\begin{aligned}
&\frac{1}{p+2}\int_{\R^N}\rho |u|^{p+2}(t,x)dx+2\mu \int^t_0\int_{\R^N} \rho|\n u|^2|u|^p (s,x)dx ds\\
&+\frac{\mu p}{2} \int^t_0 \int_{\R^{N}}\rho|\n(|u|^{2})|^{2}|u|^{p-2}(s,x)dx ds\\
&\hspace{2cm}=\frac{1}{p+2}\int_{\R^N}\rho_0 |u_0|^{p+2}(x)dx-\frac{a}{2\mu}\int^t_0\int_{\R^N}\rho(v-u)\cdot u|u|^p dx ds,
\end{aligned}
\label{gainpdu}
\end{equation}
and:
\begin{equation}
\begin{aligned}
&\frac{1}{p+2}\int_{\R^N}\rho |v|^{p+2}(t,x)dx+\frac{a}{2\mu}\int^t_0\int_{\R^N}\rho |v|^{p+2}(s,x) ds dx\\
&\hspace{2cm}=\frac{1}{p+2}\int_{\R^N}\rho_0 |v_0|^{p+2}(x)dx+\frac{a}{2\mu}\int^t_0\int_{\R^N}\rho u\cdot v|v|^p dx ds.
\end{aligned}
\label{gainpduo}
\end{equation}
Combining the two previous equation and applying Young inequality, we get for:
\begin{equation}
\begin{aligned}
&\big(\int_{\R^N}\rho |u|^{p+2}(t,x)dx+\int_{\R^N}\rho |v|^{p+2}(t,x)dx\big)+2\mu(p+2) \int^t_0\int_{\R^N} \rho|\n u|^2|u|^p (s,x)dx ds\\
&+\frac{\mu p}{2}(p+2) \int^t_0 \int_{\R^{N}}\rho|\n(|u|^{2})|^{2}|u|^{p-2}(s,x)dx ds+\frac{a(p+2)}{2\mu}\int^t_0\int_{\R^N}\rho |v|^{p+2}(s,x) ds dx\\
&\leq \big(\int_{\R^N}\rho_0 |u_0|^{p+2}(x)dx+\frac{1}{p+2}\int_{\R^N}\rho_0 |v_0|^{p+2}(x)dx\big)+\frac{a(p+2)}{2\mu}\int^t_0\int_{\R^N}\rho|v|^{p+2}dx ds\\[3mm]
&+\frac{a(p+2)}{\mu}\int^t_0\int_{\R^N}\rho|u|^{p+2}dx ds.
\end{aligned}
\label{gainpdu1}
\end{equation}
Using Gronwall lemma, we deduce that it exists $C>0$ large enough and independent on $p$ such that for any $p\geq 0$ and  for any $t\in(0,T^*)$ we have:
\begin{equation}
\|\rho^{\frac{1}{p+2}}u(t,\cdot)\|_{L^{p+2}}^{p+2}+\|\rho^{\frac{1}{p+2}}v(t,\cdot)\|_{L^{p+2}}^{p+2}\leq
(\|\rho_0^{\frac{1}{p+2}}u_0\|_{L^{p+2}}^{p+2}+\|\rho_0^{\frac{1}{p+2}}v_0\|_{L^{p+2}}^{p+2})e^{C(p+2)t}.
\label{finmp}
\end{equation}
It implies in particular that it exists $C>0$ large enough and independent on $p$ such that for any $p\geq 0$ and  for any $t\in(0,T^*)$ we have:
\begin{equation}
\|\rho^{\frac{1}{p+2}}u(t,\cdot)\|_{L^{p+2}}+\|\rho^{\frac{1}{p+2}}v(t,\cdot)\|_{L^{p+2}}\leq
C(\|\rho_0^{\frac{1}{p+2}}u_0\|_{L^{p+2}}+\|\rho_0^{\frac{1}{p+2}}v_0\|_{L^{p+2}})e^{Ct}.
\label{finsuper}
\end{equation}
\begin{remarka}
We can observe that we can extend this result to the case where the pressure $P(\rho)$ verifies $P'\in L^\infty(\R)$. 
\label{apression}
\end{remarka}
\subsection*{Gain of integrability on the velocity $u$ when $N=2$, $\gamma\geq 2$}
As previously, we wish to prove a gain of integrability on $u$, multiplying the momentum equation of (\ref{0.1}) \footnote{At this level, it is important to use the fact that ${\rm curl}u=0$ on $(0,T^*)\times(\R^N\backslash\{0\})$.} by $u|u|^p$ with $p\geq 2$ and integrating over  $(0,t)\times\R^N$ with $0<t<T^*$ we get:
\begin{equation}
\begin{aligned}
&\frac{1}{p+2}\int_{\R^N}\rho |u|^{p+2}(t,x)dx+2\mu \int^t_0\int_{\R^N} \rho|\n u|^2|u|^p (s,x)dx ds\\
&+\frac{\mu p}{2} \int^t_0 \int_{\R^{N}}\rho|\n(|u|^{2})|^{2}|u|^{p-2}(s,x)dx ds\leq\frac{1}{p+2}\int_{\R^N}\rho_0 |u_0|^{p+2}(x)dx\\
&\hspace{7cm}+|\int^t_0\int_{\R^N}a\n\rho^\gamma\cdot u|u|^{p}(t,x) dx ds |.
\end{aligned}
\label{gainp}
\end{equation}
We are now going to estimate the right hand side and by integration by parts, we have:
\begin{equation}
\begin{aligned}
&|\int_{\R^N}a\n\rho^\gamma\cdot u|u|^{p}(s,x) dx |=a|-\int_{\R^N}\rho^\gamma {\rm div}u \,|u|^p dx
-\frac{p}{2}\int_{\R^N}\rho^\gamma |u|^{p-2} u\cdot\n |u|^2\, dx|
\end{aligned}
\end{equation}
Using Young inequality we deduce that there exists $C>0$ large enough independent on $p$ such that:
\begin{equation}
\begin{aligned}
|\int_{\R^N}a\n\rho^\gamma\cdot u|u|^{p}(s,x) dx |\leq &
\frac{\mu p}{4}\int_{\R^N}\rho|\n(|u|^{2})|^{2}|u|^{p-2}(t,x)dx+\mu \int_{\R^N} \rho|\n u|^2|u|^p (s,x)dx\\
&+
(p+1) C\int_{\R^N}\rho^{2\gamma-1}|u|^p(s,x) dx.
\end{aligned}
\label{3.38}
\end{equation}
We are going now to estimate the right hand side as follows with $M\geq 4$ and using the fact that $p\geq 2$:
\begin{equation}
\begin{aligned}
&\int_{\R^N}\rho^{2\gamma-1}|u|^p(s,x) dx\leq\int_{\R^N}|\rho^{2\gamma-1-\frac{p}{p+2}}-1|\, 1_{\{\rho\geq M\}} \rho^{\frac{p}{p+2}}|u|^p(s,x) dx\\
&\hspace{2cm}+\int_{\R^N}1_{\{\rho\geq M\}}\rho^{\frac{p}{p+2}}|u|^p(s,x) dx+\int_{\R^N}\, 1_{\{\rho< M\}}\rho^{2\gamma-1}|u|^p(s,x) dx\\[2mm]
&\leq \big((\int_{\R^N}|\rho^{2\gamma-1-\frac{p}{p+2}}-1|^{\frac{p+2}{2}}1_{\{\rho \geq M\}} (s,x) dx)^{\frac{2}{p+2}}+|\{\rho(s,\cdot)\geq M\}|^{\frac{2}{p+2}}\big)
\big(\int_{\R^N} \rho |u|^{p+2}(s,x)  dx\big)^{\frac{p}{p+2}}\\
&+M^{2\gamma-2}(\int_{\R^N}\rho |u|^{p+2}(s,x) dx)^{\frac{p-2}{p}}(\int_{\R^N}\rho |u|^2 (s,x)dx)^{\frac{2}{p}}.
\end{aligned}
\label{gainp1}
\end{equation}
We wish now to estimate $(\int_{\R^N}|\rho^{2\gamma-1-\frac{p}{p+2}}-1|^{\frac{p+2}{2}}(t,x) 1_{\{\rho\geq M\}} dx)$. Before we are going to prove the following Lemma.
\begin{lemme}
\label{lemL2}
For $N=2,3$ there exists $C>0$ large enough such that for any $s\in(0,T^*)$ we have:
\begin{equation}
\begin{aligned}
&\|(\sqrt{\rho}-1)(s,\cdot)\|^2_{H^1(\R^N)}\leq C,\\
&\|(\rho-1)(s,\cdot)\|^2_{L^2(\R^N)}\leq C
%
\end{aligned} 
\label{controle}
\end{equation}
FAIRE AUSSI POUR $\gamma=1$!
\end{lemme}
{\bf Proof:} We recall now that $(\rho,u)$ verify for any $0<t<T^*$ (\ref{cenergie}) and (\ref{bcenergie}), we are going to start with dealing with the case $\gamma>1$. It implies that $\n\sqrt{\rho}$ and $\rho-1$ are respectively in $L^\infty([0,T^*];L^2(\R^N))$ and $L^\infty([0,T^*],L^\gamma_2(\R^N))$ \footnote{ We refer to \cite{Lio98} for the definition of the Orlicz space $L^\gamma_2(\R^N)$.}. Let us prove now that $\sqrt{\rho}-1$ belongs to $L^\infty([0,T^*],H^1(\R^N))$. Indeed we have for any $t\in(0,T^*)$, $M\geq 3$ and $C>0$ large enough:
\begin{equation}
\begin{aligned}
&\|(\sqrt{\rho}-1)(s,\cdot)\|^2_{L^2(\R^N)}\leq C\|(\rho-1)1_{\{|\rho-1|\leq M\}}(s,\cdot)\|_{L^2(\R^N)}^2+C \|(\rho-1)1_{\{|\rho-1|\geq M\}}(s,\cdot)\|_{L^1(\R^N)}\\
&\leq C\big(\|(\rho-1)1_{\{|\rho-1|\leq M\}}(s,\cdot)\|_{L^2(\R^N)}^2\\
&\hspace{3cm}+ \|(\rho-1)1_{\{|\rho-1|\geq M\}}(s,\cdot)\|_{L^\gamma(\R^N)}|\{|\rho(s,\cdot)-1|\geq M\}|^{\frac{1}{\gamma'}}\big),
%
\end{aligned} 
\label{H1a}
\end{equation}
with $\frac{1}{\gamma}+\frac{1}{\gamma'}=1$.
Next we have for $C>0$ large enough:
\begin{equation}
|\{|\rho(s,\cdot)-1|\geq M\}|\leq \frac{ \|(\rho-1)1_{\{|\rho-1|\geq M\}}(s,\cdot)\|_{L^\gamma(\R^N)} ^\gamma}{M^\gamma}.
\label{mesure}
\end{equation}
Combining  (\ref{cenergie}), (\ref{bcenergie}), (\ref{H1a}) and (\ref{mesure}) it implies then that $\sqrt{\rho}-1$ is bounded in $L^\infty([0,T^*],H^1(\R^N))$. We have then for  $s\in(0,T^*)$ and $C>0$ large enough and independent on $T^*$:
\begin{equation}
\begin{aligned}
&\|(\sqrt{\rho}-1)(s,\cdot)\|^2_{H^1(\R^N)}\leq C.
%
\end{aligned} 
\label{controle1u}
\end{equation}
Since we have $\rho-1=(\sqrt{\rho}-1)^2+2(\sqrt{\rho}-1)$, using (\ref{controle1u}) and Sobolev embedding, we obtain that $(\rho-1)$ is bounded in $L^\infty([0,T^*],L^2(\R^N))$. It concludes the proof of the Lemma. \\
Let us treat now the case $\gamma=1$, we recall that in this case we have:
$$\Pi(\rho)-\Pi(1)=a\big(\rho\ln(\rho)+1-\rho\big).$$
For $\delta>0$ we observe that it exists $C>0$ such that:
$$\frac{1}{C}|\rho-1|1_{\{ |\rho-1|\geq \delta\} }\leq (\Pi(\rho)-\Pi(1))1_{\{ |\rho-1|\geq \delta\} }.$$
Next since it exists $C_1>0$ such that:
$$\frac{1}{C_1}|\rho-\bar{\rho}|^21_{\{ |\rho-1|\leq \delta\} }\leq (\Pi(\rho)-\Pi(1))1_{\{ |\rho-1|\leq \delta\} }.$$
It implies in particular that $(\rho-1)$ is in $L^\infty_T(L^2_1(\R^N))$ for any $T>0$ since $(\Pi(\rho)-\Pi(1))$ is in $L^\infty_T(L^1(\R^N))$ according to the energy estimate (\ref{cenergie}). By Young inequality we deduce that it exists $C>0$ large enough:
$$\frac{1}{2}\|(\sqrt{\rho}-1)1_{\{|\rho-1|\geq\delta\}}\|^2_{L^2}\leq \|(\rho-1)   1_{\{|\rho-1  |\geq\delta\}}\|_{L^1}+C |  \{|\rho-1 |\geq\delta\}|<+\infty,$$
because $(\rho-1)\in L^\infty(L^2_1(\R^N))$. Next there exists $C>0$ such that:
$$\frac{1}{C}|\sqrt{\rho}-1|1_{\{|\rho-1|\leq\delta\}}\leq 
|\rho-1|1_{\{|\rho-1|\leq\delta\}} ,$$
it yields that $(\sqrt{\rho}-1)1_{\{|\rho-1|\leq\delta\}}$ is in $L_T^\infty(L^2)$ because $(\rho-1)$ is in $L^\infty_T(L^2_1(\R^N))$. It shows that $(\sqrt{\rho}-1)$ is in $L^\infty_T(L^2(\R^N))$. Since $\n\sqrt{\rho}$ is also in $L^\infty_T(L^2(\R^N))$ thanks to the BD entropy, it concludes the proof of the proposition for $\gamma=1$. 
%
{\hfill $\Box$}\\
\\
Next we have for $s\in(0,t)$, $p\geq 2$, $M\geq 4$, and $C_p>0$ large enough:
\begin{equation}
\begin{aligned}
&\int_{\R^N}|\rho^{2\gamma-1-\frac{p}{p+2}}-1|^{\frac{p+2}{2}}1_{\{\rho\geq M\}}(s,x)  dx\leq 
\int_{\R^N}|2(\sqrt{\rho}-1)|^{(2\gamma-1-\frac{p}{p+2})(p+2)} 1_{\{\rho\geq M\}} (s,x) dx\\
&\leq  C_p
\int_{\R^N}|\sqrt{\rho}-1|^{(2\gamma-1-\frac{p}{p+2})(p+2)}1_{\{\rho\geq M\}}  (s,x) dx.
\end{aligned}
\label{moput}
\end{equation}
From (\ref{moput}), (\ref{controle}) and Sobolev embedding in dimension $N=2$ we deduce that there exists $C^1_p>0$ such that for any $s\in(0,T^*)$ we have:
\begin{equation}
\begin{aligned}
&\int_{\R^N}|\rho^{2\gamma-1-\frac{p}{p+2}}-1|^{\frac{p+2}{2}}1_{\{\rho\geq M\}} (s,x)  dx\leq C^1_p.
\end{aligned}
\label{controle1}
\end{equation}
Combining (\ref{gainp}), (\ref{3.38}), (\ref{gainp1}), 
(\ref{mesure}), (\ref{cenergie}) and (\ref{bcenergie}) we have for any $t\in(0,T^*)$, $p\geq 2$, $M\geq 4$ and $C_p>0$ depending on $p$ large enough:
\begin{equation}
\begin{aligned}
&\frac{1}{p+2}\int_{\R^N}\rho |u|^{p+2}(t,x)dx+\mu\int^t_0 \int_{\R^N} \rho|\n u|^2|u|^p (s,x)dx ds+\frac{\mu p}{4}\int^t_0 \int_{\R^{N}}\rho|\n(|u|^{2})|^{2}|u|^{p-2}(s,x)dx ds\\
&\leq \frac{1}{p+2}\int_{\R^N}\rho_0 |u_0|^{p+2}(x)dx+C_p \int^t_0(\int_{\R^N} \rho |u|^{p+2}(s,x)  dx\big)^{\frac{p}{p+2}}ds\\
&\hspace{6cm}+M^{2\gamma-2}C_p \int^t_0(\int_{\R^N}\rho |u|^{p+2}(s,x) dx)^{\frac{p-2}{p}}ds.
\end{aligned}
\label{gainpfin}
\end{equation}
Using the Gronwall lemma, we deduce that for any $t\in(0,T^*)$ and $p\geq 2$ there exists $C_p>0$ sufficiently large depending on $p $ such that:
\begin{equation}
\|\rho^{\frac{1}{p+2}}u(t,\cdot)\|_{L^{p+2}(\R^N)}\leq C_p e^{C_p t}(1+ \|\rho_0^{\frac{1}{p+2}}u_0\|_{L^{p+2}(\R^N)}).
\label{estimfin}
\end{equation}
\subsection*{Gain of integrability on the effective velocity $v$ when $N=2$ and $\gamma\geq 2$}
Since we have seen that $(\n u-^t\n u)(t,x)=0$ for any $t\in(0,T^*)$ and any $x\in\R^N\backslash\{0\}$, we deduce from (\ref{0.1a}) that the effective velocity $v$ verifies on $(0,T^*)$:
\begin{equation}
\rho\p_t v+\rho u\cdot\n v+\frac{a\gamma}{2\mu}\rho^\gamma v=\frac{a\gamma}{2\mu}\rho^\gamma u.
\label{equasurv}
\end{equation}
Multiplying the previous equation by $v |v|^p$ with $p>0$ and integrating over $(0,t)\times\R^N$ for any $t\in(0,T^*)$ we obtain as previously:
\begin{equation}
\begin{aligned}
&\frac{1}{p+2}\int_{\R^N}\rho |v|^{p+2}(t,x)dx+\frac{a\gamma}{2\mu}\int^t_0\int_{\R^N}\rho^\gamma |v|^{p+2} dx ds\\
&\hspace{2cm}\leq \frac{1}{p+2}\int_{\R^N}\rho_0 |v_0|^{p+2}(x)dx+ \frac{a\gamma}{2\mu}|\int^t_0\int_{\R^N}\rho^\gamma u \cdot v |v|^p dx ds|.
\end{aligned}
\label{vitv}
\end{equation}
Next we have:
\begin{equation}
|\int^t_0\int_{\R^N}\rho^\gamma u \cdot v |v|^p dx ds|\leq (\int^t_0\int_{\R^N}\rho^\gamma |v|^{p+2} dx ds)^{\frac{p+1}{p+2}}(\int^t_0\int_{\R^N}\rho^\gamma |u|^{p+2} dx ds)^{\frac{1}{p+2}}.
\end{equation}
From (\ref{vitv}) and using Young inequality we have for $C_p>0$ sufficiently large:
\begin{equation}
\begin{aligned}
&\frac{1}{p+2}\int_{\R^N}\rho |v|^{p+2}(t,x)dx+\frac{a\gamma}{4\mu}\int^t_0\int_{\R^N}\rho^\gamma |v|^{p+2} dx ds\\
&\hspace{2cm}\leq \frac{1}{p+2}\int_{\R^N}\rho_0 |v_0|^{p+2}(x)dx+C_p \int^t_0\int_{\R^N}\rho^\gamma |u|^{p+2} dx ds.
\end{aligned}
\label{vitv3}
\end{equation}
We have now to deal with the term $\int^t_0\int_{\R^N}\rho^\gamma |u|^{p+2} dx ds$. We have then for $M\geq 4$ and using (\ref{controle}) and (\ref{estimfin}):
\begin{equation}
\begin{aligned}
&|\int^t_0\int_{\R^N}\rho^\gamma |u|^{p+2} dx ds|\leq \int^t_0\int_{\R^N}|\rho^{\gamma-\frac{1}{2}}-1|1_{\{|\rho-1|\geq M\}}\sqrt{\rho} |u|^{p+2} dx ds\\
&+\int^t_0\int_{\R^N}1_{\{|\rho-1|\geq M\}}\sqrt{\rho} |u|^{p+2} dx ds +\int^t_0\int_{\R^N}1_{\{|\rho-1|\leq M\}} \rho^\gamma |u|^{p+2} dx ds\\
&\leq  (\int^t_0\int_{\R^N}|\rho^{\gamma-\frac{1}{2}}-1 |^2 1_{\{|\rho-1|\geq M\}} dx ds)^{\frac{1}{2}}
( \int^t_0\int_{\R^N}\rho  |u|^{2p+4} dx ds)^{\frac{1}{2}}\\
&+( \int^t_0\int_{\R^N}\rho  |u|^{2p+4} dx ds)^{\frac{1}{2}}( \int^t_0 |\{|\rho(s,\cdot)-1|\geq M\}| ds)^{\frac{1}{2}}+(M+1)^{\gamma-1}\int^t_0\int_{\R^N}\rho |u|^{p+2} dx ds\\[2mm]
&\leq
F_p(t),
\end{aligned}
\label{estimfin3}
\end{equation}
with $F_p$ a continuous function depending on $p$ and $t$. Plugging (\ref{estimfin3}) in (\ref{vitv3}) we obtain for any $t\in(0,T^*)$, $p\geq 2$ and $C>0$ large enough that:
:
\begin{equation}
\|\rho^{\frac{1}{p+2}}v(t,\cdot)\|_{L^{p+2}(\R^N)}\leq (\|\rho_0^{\frac{1}{p+2}}v_0\|_{L^{p+2}(\R^N)} +(p+2)CF_p(t))^{\frac{1}{p+2}}.
\label{estimfin4}
\end{equation}
\subsection*{Estimates on the density $(\rho-1)$ when $N=2$ $\gamma\geq 2$ or $N=2,3$, $\gamma=1$}
From (\ref{finsuper}) and (\ref{estimfin4}) we have proved in particular that for $p=N+\e$ with $\e>0$ sufficiently small that we have for any $t\in[0,T^*)$:
\begin{equation}
\begin{aligned}
&\|\rho^{\frac{1}{p}}v(t,\cdot)\|_{L^p(\R^N)}\leq C(t),
\end{aligned}
\label{superimpo}
\end{equation}
with $C$ a continuous function on $[0,T^*]$. We recall now that the density verifies on $(0,T^*)$ the following equation:
$$\p_t\rho-2\mu\D\rho=-{\rm div}(\rho v).$$
Using the proposition \ref{chaleur}, there exists $C>0$ such that for any $T\in(0,T^*)$ we have (with $q=\rho-1$):
\begin{equation}
\begin{aligned}
\|q\|_{\widetilde{L}^\infty_T(\dot{B}^{\NN+\e'}_{p,1})}&\leq C (\|q_0\|_{\dot{B}^{\NN+\e'}_{p,1}}+\|\rho v\|_{\widetilde{L}^\infty_T(\dot{B}^{\NN-1+\e'}_{p,1})}),
\end{aligned}
\label{impoq2}
\end{equation}
with $\e'$ small enough such that $\NN-1+\e'<0$. Using Besov embedding we have:
\begin{equation}
\begin{aligned}
&\|\rho v\|_{\widetilde{L}^\infty_T(\dot{B}^{0}_{p,\infty})}=\|\rho v\|_{L^\infty_T(\dot{B}^{0}_{p,\infty})}\leq \|\rho v\|_{L^\infty_T(L^p(\R^N))}\leq \|\rho\|^{1-\frac{1}{p}}_{L^\infty_T(L^\infty)} \|\rho^{\frac{1}{p}} v\|_{L^\infty_T(L^p)},\\
&\|\rho v\|_{\widetilde{L}^\infty_T(\dot{B}^{-N(\frac{1}{2}-\frac{1}{p})}_{p,\infty})}=\|\rho v\|_{L^\infty_T(\dot{B}^{-N(\frac{1}{2}-\frac{1}{p})}_{p,\infty})}\leq \|\rho\|^{\frac{1}{2}}_{L^\infty_T(L^\infty)} \|\rho^{\frac{1}{2}} v\|_{L^\infty_T(L^2)}.
\end{aligned}
\label{impov2}
\end{equation}
Since the inequality $-\N+\NN<\NN-1+\e'<0$ holds, we get by interpolation with $\NN-1+\e'=\theta(-\N+\NN)$ and $C>0$ large enough:
\begin{equation}
\begin{aligned}
\|\rho v\|_{\widetilde{L}^\infty_T(\dot{B}^{\NN-1+\e'}_{p,1})}&\leq C \|\rho v\|_{\widetilde{L}^\infty_T(\dot{B}^{0}_{p,\infty})}^{1-\theta}\|\rho v\|_{\widetilde{L}^\infty_T(\dot{B}^{-N(\frac{1}{2}-\frac{1}{p})}_{p,\infty})}^{\theta}\\
&\leq C\|\rho\|^{\frac{\theta}{2}}_{L^\infty_T(L^\infty)} \|\rho^{\frac{1}{2}} v\|^\theta_{L^\infty_T(L^2)} \|\rho\|^{(1-\theta)(1-\frac{1}{p})}_{L^\infty_T(L^\infty)} \|\rho^{\frac{1}{p}} v\|^{1-\theta}_{L^\infty_T(L^p)}\\
&\leq  \|\rho\|^{1-\frac{1-\e'}{N}}_{L^\infty_T(L^\infty)} \|\rho^{\frac{1}{2}} v\|^\theta_{L^\infty_T(L^2)} \|\rho^{\frac{1}{p}} v\|^{1-\theta}_{L^\infty_T(L^p)}.
\end{aligned}
\label{cuestimv1}
\end{equation}
According to (\ref{superimpo})and (\ref{bcenergie}) there exists a continuous function $M_1$ on $[0,T^*]$ such that:
\begin{equation}
\begin{aligned}
\|\rho v\|_{\widetilde{L}^\infty_T(\dot{B}^{\NN-1+\e'}_{p,1})} \leq \|\rho\|^{1-\frac{1-\e'}{N}}_{L^\infty_T(L^\infty)}M_1(T).
\end{aligned}
\label{estimv1}
\end{equation}
Plugging the previous estimate in (\ref{impoq2}) it gives:
\begin{equation}
\begin{aligned}
\|q\|_{\widetilde{L}^\infty_T(\dot{B}^{\NN+\e'}_{p,1})}
&\leq C (\|q_0\|_{\dot{B}^{\NN+\e'}_{p,1}}+(\|q\|_{L^\infty_T(L^\infty)} +1)^{1-\frac{1-\e'}{N}}M_1(T)).
\end{aligned}
\label{impoq21}
\end{equation}
From Besov embedding and interpolation, we know that there exists $C,C',C_1>0$ such that:
\begin{equation}
\begin{aligned}
\|q\|_{L^\infty_T(L^\infty)}&\leq C\|q\|_{\widetilde{L}^\infty_T(\dot{B}^{\NN}_{p,1})}\leq   C'\|q\|^{\theta_1}_{\widetilde{L}^\infty_T(\dot{B}^{-N(\frac{1}{2}-\frac{1}{p})}_{p,\infty})} \|q\|^{1-\theta_1}_{\widetilde{L}^\infty_T(\dot{B}^{\NN+\e'}_{p,\infty})},\\
&\leq C_1 '\|\rho-1\|^{\theta_1}_{L^\infty_T(L^2)} \|q\|^{1-\theta_1}_{\widetilde{L}^\infty_T(\dot{B}^{\NN+\e'}_{p,\infty})}
\end{aligned}
\label{controleL}
\end{equation}
with $\NN=-\theta_1 N(\frac{1}{2}-\frac{1}{p})+(1-\theta_1)(\NN+\e')$. From the Lemma \ref{lemL2}, we have seen that $\rho-1$ is bounded
 in $L^\infty([0,T^*],L^2(\R^N))$, by Young inequality and (\ref{controleL}) we show that there exists $M'>0$ large enough such that:
\begin{equation}
\begin{aligned}
\|q\|_{L^\infty_T(L^\infty)}\leq &M' +\|q\|_{\widetilde{L}^\infty_T(\dot{B}^{\NN+\e'}_{p,1})}.
\end{aligned}
\label{controleL1}
\end{equation}
From (\ref{impoq2}), (\ref{controleL1}) and by Young inequality  we have for $C>0$ large enough:
\begin{equation}
\begin{aligned}
&\|q\|_{\widetilde{L}^\infty_T(\dot{B}^{\NN+\e'}_{p,1})}
\leq C (\|q_0\|_{\dot{B}^{\NN+\e'}_{p,1}}+(\|q\|_{\widetilde{L}^\infty_T(\dot{B}^{\NN+\e'}_{p,1})} +M'+1)^{1-\frac{1-\e'}{N}}M_1(T))\\
&\leq C \big(\|q_0\|_{\dot{B}^{\NN+\e'}_{p,1}}+\frac{(N-1+\e')}{N}(\|q\|_{\widetilde{L}^\infty_T(\dot{B}^{\NN+\e'}_{p,1})} +M'+1)+C M_1(T)^{\frac{N}{1-\e'}}\big).
\end{aligned}
\label{impoq3}
\end{equation}
It implies that for any $T\in(0,T^*)$, we have for $C>0$ large enough:
\begin{equation}
\begin{aligned}
\|q\|_{\widetilde{L}^\infty_T(\dot{B}^{\NN+\e'}_{p,1})}
&\leq C (\|q_0\|_{\dot{B}^{\NN+\e'}_{p,1}}+C_1(T)),
\end{aligned}
\label{impoq4b}
\end{equation}
with $C_1$ a continuous function on $[0,T^*]$.
From (\ref{controleL1}) it implies in particular that for any $t\in(0,T^*)$ we have:
\begin{equation}
\|\rho(t,\cdot)\|_{L^\infty(\R^N)}\leq C(t),
\label{Norminf}
\end{equation}
with $C$ a continuous function on $[0,T^*]$.
\subsection*{Gain of regularity on $\rho u$ when $N=2$ and $\gamma\geq 2$}
We recall that $\rho u$ verifies from (\ref{0.1m}) the following equation on $(0,T^*)$:
$$\p_{t}(\rho u)+\frac{1}{2}\big({\rm div}(\rho u\otimes v)+{\rm div}(v\otimes \rho u)\big)-\mu\D (\rho u)-\mu\n{\rm div}(\rho u)+\n P(\rho)=0.$$
Using proposition \ref{chaleur}, the maximum principle and the fact that  ${\rm curl}(\rho u)=0$ on $(0,T^*)$, we have for any $t\in(0,T^*)$, for $p=N+\e$ with $\e>0$ sufficiently small and $C>0$ large enough:
\begin{equation}
\begin{aligned}
&\|\rho u\|_{\widetilde{L}^\infty_t(\dot{B}^{1}_{p,\infty})+L^\infty_t(L^\infty\cap L^2)}\leq C( \|\rho_0 u_0\|_{L^\infty\cap L^2}+\|{\rm div}(\rho u\otimes v)\|_{\widetilde{L}^\infty_t(\dot{B}^{-1}_{p,\infty})}+\|\n P(\rho)\|_{\widetilde{L}^\infty_t(\dot{B}^{-1}_{p,\infty})})\\
&\leq C( \|\rho_0 u_0\|_{L^\infty\cap L^2}+\|\rho u\otimes v\|_{L^\infty_t(L^p)}+\|P(\rho)-P(1)\|_{\widetilde{L}^\infty_t(\dot{B}^{0}_{p,\infty})})\\
&\leq C( \|\rho_0 u_0\|_{L^\infty\cap L^2}+\|\rho\|_{L^\infty_t(L^\infty)}^{1-\frac{1}{p}} \|\rho^{\frac{1}{2p}}u\|_{L^\infty_t(L^{2p})} \|\rho^{\frac{1}{2p}}v\|_{L^\infty_t(L^{2p})} \\
&\hspace{8cm}+\|P(\rho)-P(1)\|_{L^\infty_t(L^p(\R^N))}).
\end{aligned}
\label{rumpi}
\end{equation}
We have now using (\ref{cenergie}) for any $t\in(0,T^*)$, $M=\|\rho\|_{L^\infty([0,T^*],L^\infty(\R^N))}+2$ and $C,C_1>0$ sufficiently large:
\begin{equation}
\begin{aligned}
&\|P(\rho(t,\cdot))-P(1)\|_{L^p(\R^N)}^p\leq \int_{\R^N}|P(\rho(t,x))-P(1)|^{p}1_{\{|\rho(t,x)-1|\leq M\}} dx\\
&\leq C
\int_{\R^N}|\rho(t,x)-1|^{\gamma p} 1_{\{|\rho(t,x)-1|\leq M\}} dx,\\
&\leq C M^{\gamma N+\gamma \e-2}\int_{\R^N}|\rho(t,x)-1|^{2} 1_{\{|\rho(t,x)-1|\leq M\}} dx\leq  C_1 M^{\gamma N+\gamma \e-2}.
\end{aligned}
\label{pression}
\end{equation}
 From (\ref{pression}) and (\ref{Norminf}) we deduce that for any $t\in(0,T^*)$, we have:
\begin{equation}
\begin{aligned}
&\|P(\rho(t,\cdot))-P(1)\|_{L^p(\R^N)}\leq C(p,t)
\end{aligned}
\label{Pimp}
\end{equation}
with $C$ a continuous function on $[0,T^*]$.
From (\ref{rumpi}), (\ref{Pimp}), (\ref{Norminf}), (\ref{estimfin}) and (\ref{estimfin4}) we deduce that for any $t\in(0,T^*)$, we have:
\begin{equation}
\begin{aligned}
&\|\rho u\|_{\widetilde{L}^\infty_t(\dot{B}^{1}_{p,\infty})+L^\infty_t(L^\infty\cap L^2)}\leq C_1(p,t),
\end{aligned}
\label{rumpi1}
\end{equation}
with $C_1$ a continuous function on $[0,T^*]$. From (\ref{cenergie}) and (\ref{Norminf}) we deduce that for any $t\in(0,T^*)$ we have:
\begin{equation}
\|\rho u(t,\cdot)\|_{L^2(\R^N)}\leq C(t),
\label{rumpi3}
\end{equation}
with $C$ continuous on $[0,T^*]$. From (\ref{rumpi1}), (\ref{rumpi2}) and by interpolation in Besov space we deduce that for $p=N+\e$ and any $t\in(0,T^*)$:
\begin{equation}
\|\rho u(t,\cdot)\|_{L^\infty_t(\dot{B}^{\NN}_{p,1})+L^\infty_t(L^\infty)}\leq C(t),
\label{rumpi2}
\end{equation}
with $C$ continuous on $[0,T^*]$. By Besov embedding we deduce that
for any $t\in(0,T^*)$:
\begin{equation}
\|\rho u(t,\cdot)\|_{L^\infty_t(L^\infty)}\leq C(t),
\label{rumpi4}
\end{equation}
with $C$ continuous on $[0,T^*]$. 
\subsection*{$L^\infty$ estimate on $v$ when $N=2$, $\gamma\geq 2$ or $N=2,3$, $\gamma=1$}
Let us start with the case $N=2$ and $\gamma\geq 2$, we wish to obtain similar estimate as (\ref{finsuper}) which is true for the case $N=2,3$, $\gamma=1$.\\
Multiplying again the equation (\ref{equasurv}) by $v |v|^p$ for any $p\geq 2$ and integrate over $(0,t)\times\R^N$ for any $t\in(0,T^*)$ we obtain as previously:
\begin{equation}
\begin{aligned}
&\frac{1}{p+2}\int_{\R^N}\rho |v|^{p+2}(t,x)dx+\frac{a\gamma}{2\mu}\int^t_0\int_{\R^N}\rho^\gamma |v|^{p+2} dx ds\\
&\hspace{2cm}\leq \frac{1}{p+2}\int_{\R^N}\rho_0 |v_0|^{p+2}(x)dx+ \frac{a\gamma}{2\mu}\int^t_0\int_{\R^N}\rho^\gamma u \cdot v |v|^p dx ds.
\end{aligned}
\label{vitvb}
\end{equation}
We have now since $\gamma\geq 2$ and using (\ref{Norminf}), (\ref{rumpi4}):
\begin{equation}
\begin{aligned}
&\int_{\R^N}\rho |v|^{p+2}(t,x)dx+\frac{a\gamma}{2\mu}(p+2)\int^t_0\int_{\R^N}\rho^\gamma |v|^{p+2} dx ds\\
&\leq \int_{\R^N}\rho_0 |v_0|^{p+2}(x)dx+ \frac{a\gamma (p+2)}{2\mu}\|\rho u\|_{L^\infty([0,T^*],L^\infty)}\|\rho\|_{L^\infty([0,T^*],L^\infty)}^{\gamma-2}  \int^t_0\int_{\R^N} \rho   |v|^{p+1} dx ds\\
&\leq \int_{\R^N}\rho_0 |v_0|^{p+2}(x)dx\\
&+ \frac{a\gamma (p+2)}{2\mu}\|\rho u\|_{L^\infty([0,T^*],L^\infty)}\|\rho\|_{L^\infty([0,T^*],L^\infty)}^{\gamma-2}  \int^t_0\|\rho^{\frac{1}{p+2}}v(s,\cdot)\|_{L^{p+2}}^{\frac{(p-1)(p+2)}{p}} \|\sqrt{\rho}v(s,\cdot)\|_{L^2}^{\frac{2}{p}} ds.
\end{aligned}
\label{vitvbb}
\end{equation}
We deduce from  (\ref{rumpi4}), (\ref{Norminf}), (\ref{cenergie}) and (\ref{bcenergie}) that for any $t\in(0,T^*)$ there exists $C>0$ independent on $t$ such that:
\begin{equation}
\begin{aligned}
&\|\rho^{\frac{1}{p+2}}v(t,\cdot)\|_{L^{p+2}}^{p+2}\\
&\leq \|\rho_0^{\frac{1}{p+2}}v_0\|_{L^{p+2}}^{p+2}+C(p+2)\big(t+\int^t_0\|\rho^{\frac{1}{p+2}}v(s,\cdot)\|^{p+2}ds\big).
\end{aligned}
\label{vitvbbb}
\end{equation}
Using the Gronwall lemma, it yields that there exists a function $C$ continuous such that for any $p\geq 2$:
\begin{equation}
\begin{aligned}
&\|\rho^{\frac{1}{p+2}}v(t,\cdot)\|_{L^{p+2}}\leq C(t).
\end{aligned}
\label{vitvbbbfin}
\end{equation}
It is important to observe that the previous estimate is uniform in $p\geq 2$. This estimate is true for the case $N=2$, $\gamma\geq 2$ or $N=2,3$, $\gamma=1$.
\begin{proposition}\label{prop5}
We have then when  $N=2$, $\gamma\geq 2$ or $N=2,3$, $\gamma=1$:
\begin{equation}
	\forall T\in(0,T^*), \;\|v\|_{ L^\infty_T(L^\infty)}\leq C(T),
\label{Lfini}
\end{equation}
with $C$ a continuous function on $[0,T^*]$.
\end{proposition}
{\bf Proof:} We recall that for all $t\in (0,T^*)$ the solution $(\rho,v)$ verifies \footnote{This is true since $(\rho,v)$ is a strong regular solution on the maximal time interval $(0,T^*)$. It implies in particular that for any $t\in(0,T^*)$, we have $\frac{1}{\rho}$ and $v$ wich belong to $L^\infty([0,t],L^\infty(\R^N))$.}:
\begin{equation}
\rho(t,x)\geq  B_t>0\;\;\forall x\in\R^N\;\;\mbox{and}\;\;\|v(t,\cdot)\|_{L^\infty}\leq C^1_t<+\infty,
\label{3.76} 
\end{equation}
with possibly $B_t\rightarrow_{t\rightarrow T^*} 0$ and $C^1_t\rightarrow_{t\rightarrow T^*} +\infty$.
We observe that $\forall \e>0$ sufficiently small (such that $\e<\frac{\|v(t,\cdot)\|_{L^\infty}}{2}$ for $\|v(t,\cdot)\|_{L^\infty}\ne 0$), 
we have for any $p\geq 2$ and $t\in(0,T^*)$:
\begin{equation}
\begin{aligned}
\|\rho^{\frac{1}{p}}v(t)\|_{L^p}
&\geq \big(\int_{\{x,\;|v(t,x)|\geq \|v(t,\cdot)\|_{L^\infty}-\e\}}
    \rho(t,x) |v|^p (t,x) dx \big)^{\frac{1}{p}}\\  
&\geq \big(\|v(t,\cdot)\|_{L^\infty}-\e\big) 
    \big(\int_{\{x,\;|v(t,x)|\geq \|v(t,\cdot)\|_{L^\infty}-\e\}}
   \rho(t,x)  dx\big)^{\frac{1}{p}}\\
	&\geq \bigl(\|v(t,\cdot)\|_{L^\infty}-\e\bigr) B_{t}^{\frac{1}{p}}
    \,\big|\{x,\;|v(t,x)|\geq\|v(t,\cdot)\|_{L^\infty}-\e\} \big|^{\frac{1}{p}}.
\end{aligned}
\label{norm}
\end{equation}
Since we have $B_{t}>0$ and 
$0<\big| \{x,\;|v(t,x)|\geq \|v(t,\cdot)\|_{L^\infty}-\e \} \big|<+\infty$ (indeed we recall that $v(t,\cdot)$ via (\ref{bcenergie}) and (\ref{3.76}) is in $L^2(\R^N)$), 
we can pass to the limit when $p$ goes to $+\infty$ in \eqref{norm}. 
It implies that for any $\e>0$ small enough, we get using (\ref{vitvbbbfin}):
$$\|v(t,\cdot)\|_{L^\infty}-\e\leq C(t),\;\;\forall t\in (0,T^*).$$
It concludes the proof of the proposition \ref{prop5}.  {\hfill $\Box$}
\\
\\
It implies in particular that $v$ belongs to $L^\infty([0,T^*],L^\infty(\R^N))$.
\subsection*{Control of the $L^\infty$ norm of $\frac{1}{\rho}$ when $N=2$, $\gamma\geq 2$ or $N=2,3$, $\gamma=1$}
We recall now that the density verifies on $(0,T^*)$ the following equation:
$$\p_t\rho-2\mu\D\rho=-{\rm div}(\rho v).$$
Using (\ref{Lfini}) and the maximum principle (see \cite{La}), we deduce that for any $t\in(0,T^*)$:
\begin{equation}
\|\frac{1}{\rho}(t,\cdot)\|_{L^\infty(\R^N)}\leq C(t),
\label{superimp}
\end{equation}
with $C$ a continuous function on $[0,T^*]$. 
\subsection*{Global strong solution for the approximate solutions when $N=2$, $\gamma\geq 2$ and $N=2,3$, $\gamma=1$}
We recall that we have constructed a sequence of strong solution $(\rho_n,u_n,v_n)_{n\in\mathbb{N}}$ for the system (\ref{0.1}) in finite time on the time interval $[0,T_n^*)$ with $T_n^*$ the maximal lifespan time. These solutions have for initial data an approximate sequence $(\rho_n^0,u^0_n,v^0_n)_{n\in\mathbb{N}}$ (see (\ref{minitial})).\\
We would like now to prove that for any $n\in\mathbb{N}$, we have $T_n^*=+\infty$. To do this we wish to use the blowup criterion of the Theorem \ref{theo1a} with $N<p<2N$ and $p$ defined as in the Theorem \ref{theo4}. It suffices then to verify the three assertions are satisfied.\\
 First from (\ref{superimp}), we have seen that $\frac{1}{\rho_n}$ belongs to $L^\infty([0,T_n^*],L^\infty(\R^N))$. It remains now to prove the two following claims:
 \begin{equation}
 \int_{0}^{T_n^*}\|\n u_n(t,\cdot)\|_{L^\infty} dt<+\infty
 \label{ghyp1}
 \end{equation}
 and:
  \begin{equation}
 \|\rho_n-1\|_{L^\infty([0,T_n^*],\dot{B}^{\NN}_{p,1})}<+\infty. 
 \label{ghyp2}
 \end{equation} 
We point out that (\ref{ghyp2}) is a direct consequence of (\ref{impoq4b}) and (\ref{controle}) by interpolation in Besov spaces. Let us prove now
(\ref{ghyp1}), we have then on $(0,T_n^*)$:
$$\p_t u_n-2\mu\D u_n=-\frac{1}{\rho_n}\n P(\rho_n)-2 u_n\cdot\n u_n+v_n\cdot\n u_n=F_n.$$
Using proposition \ref{chaleur}, we have for any $T\in(0,T_n^*)$ and $C>0$ large enough with $p$ defined as in the Theorem \ref{theo4}:
\begin{equation}
\begin{aligned}
&\|u_n\|_{\widetilde{L}^\infty_T(\dot{B}^{\NN-1}_{p,1})}+\|u_n\|_{\widetilde{L}^1_T(\dot{B}^{\NN+1}_{p,1})}
\leq C(\|u_n(0,\cdot)\|_{\dot{B}^{\NN-1}_{p,1}}+\|F_n\|_{\widetilde{L}^1_T(\dot{B}^{\NN-1}_{p,1})}).
\end{aligned}
\label{particum}
\end{equation}
Let us estimate now $\|F_n\|_{\widetilde{L}^1_T(B^{\NN-1}_{p,1})}$, we have then using classical paraproduct laws and interpolation for Besov spaces with $\e'>0$ such that $-\N+\NN<\NN-1+\e'<0$ and $0<\theta<1$:
\begin{equation}
\begin{aligned}
&\|u_n\cdot\n u_n\|_{\widetilde{L}^1_T(\dot{B}^{\NN-1}_{p,1})}\leq C \|u_n\|_{L^\infty_T(\dot{B}^{\NN-1+\e'}_{p,1})} \|\n u_n\|_{L^1_T(\dot{B}^{\NN-\e'}_{p,1})}\\
&\leq C_1 \|u_n\|_{L^\infty_T(\dot{B}^{\NN-1+\e'}_{p,1})} \|\n u_n\|^\theta_{L^1_T(\dot{B}^{\NN}_{p,\infty})}\|\n u_n\|^{1-\theta}_{L_T^1(\dot{B}^{-N(\frac{1}{2}-\frac{1}{p})}_{p,\infty})}\\
&\leq C_2\|u_n\|_{L^\infty_T(\dot{B}^{\NN-1+\e'}_{p,1})} \|\n u_n\|^\theta_{\widetilde{L}^1_T(\dot{B}^{\NN}_{p,1})}\|\n u_n\|^{1-\theta}_{L_T^1(L^2)},
\end{aligned}
\label{3.86}
\end{equation}
with $C, C_1, C_2>0$ large enough.
 Since the inequality $-\N+\NN<\NN-1+\e'<0$ holds, we get by interpolation with $\NN-1+\e'=\theta_1(-\N+\NN)$ and $C>0$ large enough:
\begin{equation}
\begin{aligned}
\|u_n\|_{L^\infty_T(\dot{B}^{\NN-1+\e'}_{p,1})}&\leq C \|u_n\|_{L^\infty_T(\dot{B}^{0}_{p,\infty})}^{1-\theta_1}\|u_n\|_{L^\infty_T(\dot{B}^{-N(\frac{1}{2}-\frac{1}{p})}_{p,\infty})}^{\theta_1}\\
&\leq \|u_n\|^{\theta_1}_{L^\infty_T(L^2)} \|u_n\|^{1-\theta_1}_{L^\infty_T(L^p)}.
\end{aligned}
\label{cuestimv1g}
\end{equation}
Combining (\ref{3.86}), (\ref{cuestimv1g}), (\ref{cenergie}), (\ref{bcenergie}), (\ref{superimp}), (\ref{estimfin}), (\ref{finsuper}) we have using Young inequality with $\e_1>0$ arbitrary small:
\begin{equation}
\begin{aligned}
&\|u_n\cdot\n u_n\|_{\widetilde{L}^1_T(\dot{B}^{\NN-1}_{p,1})}\leq \e_1 \|\n u_n\|^\theta_{\widetilde{L}^1_T(\dot{B}^{\NN}_{p,1})}+C_{\e_1}(T),
\end{aligned}
\label{3.86a}
\end{equation}
with $C_{\e_1}$ a continuous function on $[0,T_n^*]$. In a similar way, we can prove that:
\begin{equation}
\begin{aligned}
&\|v_n\cdot\n u_n\|_{\widetilde{L}^1_T(\dot{B}^{\NN-1}_{p,1})}\leq \e_1 \|\n u_n\|^\theta_{\widetilde{L}^1_T(\dot{B}^{\NN}_{p,1})}+C'_{\e_1}(T),
\end{aligned}
\label{3.86b}
\end{equation}
with $C'_{\e_1}$ a continuous function on $[0,T_n^*]$. 
Let us deal now with the term $\frac{1}{\rho_n}\n P(\rho_n)$ when $N=2$, $\gamma\geq 2$, we recall that $\frac{1}{\rho_n}\n P(\rho_n)=\frac{a\gamma}{2\mu}\rho_n^{\gamma-1}(v_n-u_n)$ and $\frac{1}{\rho_n}\n P(\rho_n)=\frac{a\gamma}{\gamma-1}\n\rho_n^{\gamma-1}$, it implies by interpolation that for $C>0$ large enough and $\theta\in(0,1)$ we get:
\begin{equation}
\begin{aligned}
\|\frac{1}{\rho_n}\n P(\rho_n)\|_{\widetilde{L}^1_T(\dot{B}^{\NN-1}_{p,1})}&\leq C  \|\rho_n^{\gamma-1}(v_n-u_n)\|_{L^\infty_T(\dot{B}^0_{p,\infty})}^{\theta} \|\n(\rho_n^{\gamma-1}-1)\|^{1-\theta}_{L^\infty_T(\dot{B}^{-N(\frac{1}{2}-\frac{1}{p})-1}_{p,\infty})}\\
&\leq C_1  \|\rho_n^{\gamma-1}(v_n-u_n)\|_{L^\infty_T(L^p)}^{\theta} \|(\rho_n^{\gamma-1}-1)\|_{L^\infty_T(L^2)}^{1-\theta}
\end{aligned}
\label{3.86c}
\end{equation}
with $C,C_1>0$ large enough.
Combining (\ref{controle}), (\ref{cenergie}), (\ref{bcenergie}), (\ref{superimp}), (\ref{estimfin}), (\ref{vitvbbbfin}), (\ref{Norminf}) we deduce that there exists a continuous function $C$ on $[0,T_n^*]$ such that:
\begin{equation}
\begin{aligned}
&\|\frac{1}{\rho_n}\n P(\rho_n)\|_{\widetilde{L}^1_T(\dot{B}^{\NN-1}_{p,1})}\leq  C(T).
\end{aligned}
\label{3.86d}
\end{equation}
For the case $N=2,3$, $\gamma=1$, we have $\frac{1}{\rho_n}\n P(\rho_n)=\frac{a}{2\mu}(v_n-u_n)=a\n\ln\rho_n$, and by interpolation with $\theta\in(0,1)$ we obtain then by interpolation:
\begin{equation}
\begin{aligned}
\|\frac{1}{\rho_n}\n P(\rho_n)\|_{\widetilde{L}^1_T(\dot{B}^{\NN-1}_{p,1})}&\leq C  \|(v_n-u_n)\|_{L^\infty_T(\dot{B}^0_{p,\infty})}^{\theta} \|\n\ln\rho_n\|^{1-\theta}_{L^\infty_T(\dot{B}^{-N(\frac{1}{2}-\frac{1}{p})-1}_{p,\infty})}\\
&\leq C_1  \|(v_n-u_n)\|_{L^\infty_T(L^p)}^{\theta} \|\ln\rho_n\|_{L^\infty_T(L^2)}^{1-\theta}
\end{aligned}
\label{3.86e}
\end{equation}
Combining (\ref{cenergie}), (\ref{bcenergie}), (\ref{superimp}), (\ref{finsuper}), (\ref{Norminf}) we obtain:
\begin{equation}
\begin{aligned}
&\|\frac{1}{\rho_n}\n P(\rho_n)\|_{\widetilde{L}^1_T(\dot{B}^{\NN-1}_{p,1})}\leq  C(T),
\end{aligned}
\label{3.86da}
\end{equation}
with $C$ a continuous function on $[0,T_n^*]$. Using (\ref{particum}), (\ref{3.86a}), (\ref{3.86b}), (\ref{3.86d}) and (\ref{3.86da}) and
taking $\e_1>0$ sufficiently small, there exists a continuous function $C$ on $[0,T_n^*]$ such that for any $T\in(0,T_n^*)$:
\begin{equation}
\begin{aligned}
&\|u_n\|_{\widetilde{L}^\infty_T(\dot{B}^{\NN-1}_{p,1})}+\|u_n\|_{\widetilde{L}^1_T(\dot{B}^{\NN+1}_{p,1})}
\leq C(T).
\end{aligned}
\label{particumv}
\end{equation}
It yields in particular from Besov embedding that there exists a continuous function $C$ on $[0,T_n^*]$ such that for any $T\in(0,T_n^*)$:
\begin{equation}
\begin{aligned}
&\|\n u_n\|_{L^1_T(L^\infty)} \leq C(T).
\end{aligned}
\label{particumv}
\end{equation}
It proves in particular (\ref{ghyp1}) since $C$ is continuous on $[0,T_n^*]$ and we have  $\lim_{t\rightarrow T_n^*} C(t)<+\infty$.
We have then proved that $(\rho_n,v_n,m_n)_{n\in\mathbb{N}}$ exists globally in time.\\
From the Theorem \ref{theo1a}, (\ref{cenergie}), (\ref{bcenergie}), (\ref{Norminf}), (\ref{Lfini}) and (\ref{superimp}), we deduce that for any $T\in\R$ and any $n\in\mathbb{N}$, we have for $C_1>0$ large enough and $p,\e$ defined in the Theorem \ref{theo4}:
\begin{equation}
\begin{cases}
\begin{aligned}
&\|\rho_n-1\|_{\widetilde{L}^\infty_T(\dot{B}^{\NN}_{p,1}\times \dot{B}^{\NN+\e}_{p,1})}+\|u_n\|_{\widetilde{L}^\infty_T(\dot{B}^{\NN-1}_{p,1}\times \dot{B}^{\NN-1+\e}_{p,1})}+\|u_n\|_{\widetilde{L}^1_T(\dot{B}^{\NN+1}_{p,1}\times \dot{B}^{\NN+1+\e}_{p,1})}\leq C(T),\\
&\|(\frac{1}{\rho_n},\rho_n)\|_{L^\infty_T(L^\infty)}\leq C(T),\\
&\|v_n\|_{L^\infty_T(L^\infty)}\leq C(T),\\
&{\cal E}(\rho_n,u_n)(T)\leq C_1,\,{\cal E}_1(\rho_n,v_n)(T)\leq C_1
\end{aligned}
\end{cases}
\label{estimcrucial}
\end{equation}
with $C$ locally bounded. It is important to point out that the function $C$ is independent on $n$. 
\subsection*{Compactness and global strong solution}
From (\ref{estimcrucial}) it is now classical to prove that there exists a subsequence $(\rho_{\va(n)},v_{\va(n)},u_{\va(n)})_{n\in\mathbb{N}}$ such that $(\rho_{\va(n)},v_{\va(n)},u_{\va(n)})_{n\in\mathbb{N}}$ converges in the sense of the distribution to a global weak solution $(\rho,v,u)$ of the system (\ref{0.1}) with initial data $(\rho_0,v_0,u_0)$. We refer for a proof to \cite{Danchin}. In addition $(\rho,v,u)$ verifies for any $T>0$ and for $C_1>0$ large enough:
\begin{equation}
\begin{cases}
\begin{aligned}
&\|\rho-1\|_{\widetilde{L}^\infty_T(\dot{B}^{\NN}_{p,1}\times B^{\NN+\e}_{p,1})}+\|u\|_{\widetilde{L}^\infty_T(\dot{B}^{\NN-1}_{p,1}\times \dot{B}^{\NN-1+\e}_{p,1})}+\|u\|_{\widetilde{L}^1_T(\dot{B}^{\NN+1}_{p,1}\times \dot{B}^{\NN+1+\e}_{p,1})}\leq C(T),\\
&\|(\frac{1}{\rho},\rho)\|_{L^\infty_T(L^\infty)}\leq C(T),\\
&\|v\|_{L^\infty_T(L^\infty)}\leq C(T),\\
&{\cal E}(\rho,u)(T)\leq C_1,\,{\cal E}_1(\rho,v)(T)\leq C_1
\end{aligned}
\end{cases}
\label{estimcrucial1}
\end{equation}
with $p$, $\e>0$ defined as in the Theorem \ref{theo4} and $C$ locally bounded.
From (\cite{Fourier}), we deduce that the solution $(\rho,u,v)$ is unique. It implies in particular that $(\rho,u,v)$ is a global strong solution for the system (\ref{0.1}) with initial data $(\rho_0,u_0,v_0)$.
\section{Proof of the Theorem \ref{theo2}}
\label{section4}
\subsubsection*{Existence of approximate solutions and conservation of the geometrical structure}
Let $(\rho_0,u_0,v_0)$ verifying the assumptions of the Theorem \ref{theo2} with $\frac{1}{\rho_0}\in L^\infty(\R^N)$ (we will treat the case where $v_0$ is in $L^1(\R^N)$ but $\frac{1}{\rho_0}$ is not in $L^\infty(\R^N)$ later, see the Remark \ref{reminio}), we define
now the regularizing initial data:
\begin{equation}
\begin{cases}
\begin{aligned}
&\rho_0^n=\varphi_n j_n*\rho_0+\frac{1}{n}, m_{0,1}^n=\varphi_n j_n*(\rho_0 v_0)\\
& m_0^n=m_{0,1}^n
-2\mu\n\rho_0^n\\
&u_0^n=\frac{m_0^n}{\rho_0^n}, v_0^n=\frac{m_{0,1}^n}{\rho_0^n}
\end{aligned}
\end{cases}
\label{definitu}
\end{equation}
with $j_n$ a regularizing kernel and $\varphi_n$ a truncature. More precisely $j_n=n^Nj(n\cdot)$ and $j$ belongs to $C^{\infty}_0(\R^N)$, is radial (in particular $j_n(\cdot)=j^1_n(|\cdot|)$) and verifies:
$$0\leq j\leq 1, \,\mbox{supp}j\in B(0,1)\;\;\mbox{and}\;\;\int_{\R^N}j(x) dx=1.$$
In addition we take $\va_n(x)=\va(\frac{x}{n})$ with $\va\in C^{\infty}_0(\R^N)$ and radial verifying:
$$\va=1\;\;\;\mbox{on}\;B(0,1),\;\;0\leq\va\leq 1\;\;\mbox{and}\;\;\mbox{supp}\va\in B(0,2).$$
We observe in particular that we have:
\begin{equation}
\begin{cases}
\begin{aligned}
&\rho_0^n\rightarrow _{n\rightarrow+\infty}\rho_0\;\;\mbox{weakly * in}\;L^\infty(\R^N),\\
&m_{0,1}^n\rightarrow _{n\rightarrow+\infty}\rho_0 v_0\;\;\mbox{weakly * in}\;L^\infty(\R^N).
\end{aligned}
\end{cases}
\label{confaible}
\end{equation}
First we observe easily that $(\rho_0^n)_{n\in\mathbb{N}}$ and $(m_{0,1}^n)_{n\in\mathbb{N}}$ are uniformly bounded in $L^\infty(\R^N)$.
Now let $\psi\in L^1(\R^N)$, we have since $j_n(y-z)=j_n(z-y)$ (indeed $j_n$ is radial):
$$
\begin{aligned}
\int_{\R^N}\rho_0^n(y)\psi(y) dy&=\int_{\R^N}\int_{\R^N}\va_n(y)j_n(y-z)\rho_0(z)\psi(y) dy dz+\frac{1}{n}\int_{\R^N}\psi(y) dy\\
&=\int_{\R^N}\int_{\R^N}\rho_0(y) (j_n*(\va_n\,\psi)) (y) dy +\frac{1}{n}\int_{\R^N}\psi(y) dy
\end{aligned}
$$
Next we observe that $(j_n*(\va_n\,\psi))_{n\in\mathbb{N}}$ converges strongly to $\psi$ in $L^1(\R^N)$. It implies then that:
$$
\begin{aligned}
\int_{\R^N}\rho_0^n(y)\psi(y) dy\rightarrow_{n\rightarrow+\infty} \int_{\R^N}\int_{\R^N}\rho_0(y) \psi(y) dy.
\end{aligned}
$$
We proceed similarly for $(m_{0,1}^n)_{n\in\mathbb{N}}$. Now since $j_n$, $\va_n$ and $\rho_0$ are radial , we deduce that $\rho_0^n$ is radial as convolution of radial functions and product of radial functions. We recall now that $v_0$ can be written as follows:
$$\rho_0 v_0(x)=\frac{x}{|x|}\rho_0(|x|)v^1_0(|x|)=\n_x m^2_0(|x|),$$ 
with $v^1_0\in L^\infty(\R)$ and $(m^2_0)'=\rho_0 v^1_0$. We have now:
$$
\begin{aligned}
&\rho_0^n v_0^n(x)=\va_n(x)\int_{\R^N}j^1_n(|x-y|)\n_y[ m^2_0(|y|)] dy\\
&=-\va_n(x)\int_{\R^N}\n_y \,j^1_n(|x-y|)m^2_0(|y|) dy=\va_n(x)\int_{\R^N}\n_x \,j^1_n(|x-y|)m^2_0(|y|) dy\\
&=\va_n(x)\n_x(j^1_n*m_0^2(|\cdot|)).
\end{aligned}
$$  
Since $j_n$ and $m_0^2(|\cdot |)$ are radial, we deduce that $j_n*m_0^2(|\cdot|))$ is radial (with $j_n*m_0^2(|\cdot|)(x)=j_n^a(|x|)$)
 and belongs to $C^\infty(\R^N)$. We deduce then that we can write $\rho_0^n v_0^n$ as follows:
$$\rho_0^n v_0^n(x)=\frac{x}{|x|} \va_n(x)(j_n^a)'(|x|).$$ 
It implies in particular that $m_0^n$ and $u_0^n$ are also radially symmetric. Now we observe since $\rho_0$ is positive and belongs to $L^\infty(\R^N)$ that $\rho_0^n$ verifies:
\begin{equation}
\frac{1}{n}\leq\rho_0^n\leq\|\rho_0\|_{L^\infty}+\frac{1}{n}.
\end{equation}
In addition we can verify that $\rho_0^n-\frac{1}{n}$ belongs to $L^2(\R^N)\cap L^\infty(\R^N)$
and by interpolation we deduce that $(\rho_0^n-\frac{1}{n}) $ is in $L^\gamma_2(\R^N)$. In a similar way, we observe that $\n\sqrt{\rho_0^n}=\frac{\n\rho_0^n}{2\sqrt{\rho_0^n}}$ belongs to $L^2(\R^N)$ and $v_0^n$ is also in $L^2(\R^N)$.
 It implies in particular that $(\rho_0^n,v_0^n,u_0^n)$ are of finite energy for the entropy ${\cal E}$ and ${\cal E}_1$. \\
We verify in fact using interpolation and standard estimate for the convolution  that $(\rho_0^n-\frac{1}{n},u_0^n)$ verify the conditions of the Theorem \ref{theo4}. Indeed for example, we observe that $(\rho_0^n-\frac{1}{n})$ belongs to $W^{k,p}(\R^N)$ with any $k\in\mathbb{N}$ and $p>N$, it implies by interpolation in Besov space that $(\rho_0^n-\frac{1}{n})\in \dot{B}^{\NN}_{p,1}(\R^N)\cap \dot{B}^{\NN+\e}_{p,1}(\R^N)$
with $\e>0$ arbitrary small. Similarly $v_0^n$ is in $L^1(\R^N)\cap L^p(\R^N)$ which implies by interpolation that $v_0^n$ is in $\dot{B}^{\NN-1}_{p,1}(\R^N)\cap \dot{B}^{\NN-1+\e}_{p,1}(\R^N)$ with $\e>0$ such that $\NN-1+\e<0$. 
\\
 From the Theorem \ref{theo4}, we know that there exists a global strong solution $(\rho_n,u_n,v_n)$ for the system (\ref{0.1}) with initial data $(\rho_0^n,u_0^n,v_0^n)$ and with $P(\rho)=a\rho^\gamma$ ($\gamma\geq 2$) if $N=2$ and with the following approximate regular pressure $P_n(\rho)$ if $N=3$ defined as follows:
\begin{equation}
\begin{cases}
\begin{aligned}
&P_n(\rho)=P(\rho)\;\;\mbox{if}\;\;\rho\leq n\\
&P_n(\rho)=a 2^{\gamma-1}\rho n^{\gamma-1}\;\;\mbox{if}\;\;\rho\geq 2n
\end{aligned}
\end{cases}
\label{defPn}
\end{equation}
with for $C>0$ independent on $n$:
$$P'_n(\rho)\leq Ca(\gamma-1)n^{\gamma-1}\;\;\;\mbox{on}\;\;[n,2n].$$
In addition we choose $P_n$ such that $P_n$ is increasing. It implies in particular that for $C>0$ large enough and independent on $n$, we have:
\begin{equation}
\begin{cases}
\begin{aligned}
&\|\frac{P'_n(\rho)}{\rho}\|_{L^\infty}\leq C\|\rho^{\gamma-2}\|_{L^\infty}\\
&\|P'_n(\rho)\|_{L^\infty}\leq C\\
&\|P_n(\rho)\|_{L^\infty}\leq C\|\rho^\gamma\|_{L^\infty}.
\end{aligned}
\end{cases}
\label{Pressiontech}
\end{equation}
Indeed we have observed (see the Remark \ref{apression}) that we have a global strong solution provided that  the pressure $P_n$ satisfies $P_n'\in L^\infty(\R)$.\\
We know more precisely from the Theorem \ref{theo4} that $(\rho_n-\frac{1}{n},u_n)$ verifies for any $n\in\mathbb{N}$, $p>N$ and $\e>0$ arbitrary small:
\begin{equation}
\begin{cases}
\begin{aligned}
&(\rho_n-\frac{1}{n},u_n)\in \widetilde{L}_{loc}^{\infty}(\dot{B}^{\NN}_{p,1}\cap \dot{B}^{\NN+\e}_{p,1})\times( \widetilde{L}^{\infty}_{loc}(\dot{B}^{\NN-1}_{p,1}\cap B^{\NN-1+\e}_{p,1})\cap  \widetilde{L}^{1}_{loc}(\dot{B}^{\NN+1}_{p,1}\cap \dot{B}^{\NN+1+\e}_{p,1}))^N,\\
&(\rho_n,\frac{1}{\rho_n},v_n)\in L^\infty_{loc}(L^\infty(\R^N))^{N+2}
\end{aligned}
\end{cases}
\label{estimo1}
\end{equation}
The previous estimates allow in particular to conserve any regularity (as it is well known for compressible Navier-Stokes equation), we have in particular for any $n\in\mathbb{N}$, $p>N$, $s>\NN$:
\begin{equation}
\begin{aligned}
&(\rho_n-\frac{1}{n},u_n)\in \widetilde{L}_{loc}^{\infty}(B^{s}_{p,1})\times( \widetilde{L}^{\infty}_{loc}(B^{s-1}_{p,1})\cap  \widetilde{L}^{1}_{loc}(B^{s+1}_{p,1}))^N.
\end{aligned}
\label{estimo2}
\end{equation}
It implies in particular that we have for any $n\in\mathbb{N}$ and $m_n=\rho_n u_n$:
\begin{equation}
\begin{aligned}
&(m_n,v_n)\in L^\infty_{loc}(\R^+,L^\infty(\R^N)).
\end{aligned}
\label{estimo3}
\end{equation}
It is important to point out that all these estimates are not uniform in $n$.\\
We are going now to prove uniform estimates in $n$ on the solutions $(\rho_n,m_n,v_n)_{n\in\mathbb{N}}$ which is in $C^\infty(\R^+\times\R^N)$ since the Besov norm are propagated. When we will have sufficiently regular uniform estimates, we will verify that the sequence $(\rho_n,m_n,v_n)_{n\in\mathbb{N}}$ converges up to a subsequence in a suitable functional space to a weak solution $(\rho,m,v)$ of (\ref{0.1m}) in finite time.\\
\subsection*{Uniform estimates in $n$ on $(\rho_n,v_n,m_n)_{n\in\mathbb{N}}$}
We want now to prove uniform estimates in $n$ on the sequel $(\rho_n,v_n,m_n)_{n\in\mathbb{N}}$ in $Y_T$ with $T>0$ sufficiently small, for the simplicity in the sequel we forget the subscript $n$.
\subsubsection*{$bmo_T^{-1}$ estimates on the momentum $m$}
We are going to estimate the momentum $m$ in the following space ${\cal E}_T$ for $0<T<1$  sufficiently small (we will define it later):
\begin{equation}
\|m\|_{{\cal E}_T}=\sup_{0<t<T}\sqrt{t}\|m(t)\|_{L^\infty(\R^N)}+\big(\sup_{x\in\R^N,\,0<t<T}t^{-\N}\int^t_0 \int_{y\in B(x,\sqrt{t})}|m(t,y)|^2 dt \,dy\big)^{\frac{1}{2}}
\label{3}
\end{equation}
From (\ref{estimo3}), we recall that $m$ belongs to $L^\infty_{loc}(L^\infty(\R^N))$, in particular it implies that $\|m\|_{{\cal E}_T}$ is well defined. In the sequel we will use also the quantities $\|v\|_{L^\infty([0,T],L^\infty(\R^N))}$ and $\|\rho\|_{L^\infty([0,T],L^\infty(\R^N))}$ which are also well defined from (\ref{estimo1}).\\
We recall that we have from (\ref{0.3}) (we take $2\mu=1$ in order to simplify the notation):
\begin{equation}
m(t,\cdot)=e^{ t\D}m_0-\int^t_0e^{(t-s)\D}({\rm div}(v\otimes m)-\n P(\rho))(s,\cdot) ds.
\label{Duha}
\end{equation}
Since $m_0=\rho_0 v_0+2\mu\n\rho_0$ and $(\rho_0,v_0)$ belongs to $L^\infty(\R^N)$ which is embedded in $bmo(\R^N)$, it implies then that $m_0$ belongs to $bmo^{-1}(\R^N)$.
  We have then by definition of $bmo^{-1}(\R^N)$ (see \cite{Lemarie}) that $e^{t\D}m_0$ belongs to ${\cal E}_T$ when $T\leq 1$. More precisely we get for $C>0$:
\begin{equation}
\|e^{ t\D}m_0\|_{{\cal E}_T}\leq C \|m_0\|_{bmo^{-1}},
\label{Duha1}
\end{equation}
with $\|m_0\|_{bmo^{-1}}=(\sup_{0<t<1, x\in\R^N}\frac{1}{t^\N}\int^t_0\int_{B(x,\sqrt{t})}|e^s\D m_0(x)|^2 dx ds)$.
We denote now by $B(m,v)=\int^t_0e^{(t-s)\D}{\rm div}(v\otimes m)ds.$
We are interesting in proving now that there exists $C>0$ such that:
\begin{equation}
\begin{aligned}
&\|B(m,v)\|_{{\cal E}_T}\leq C \sqrt{T} \|m\|_{{\cal E}_T}\|v\|_{L^\infty_T(L^\infty(\R^N))}.
\end{aligned}
\label{cru}
\end{equation}
To do this we are going to follow the approach developed by Lemari\'e in \cite{Lemarie}.
Next let us start with the $L^\infty$ norm of $B(m,v)$, we split the integral in two domains $0<s<\frac{t}{2}$ and $\frac{t}{2}\leq s<t$. 
We have then by integration by parts and making the change of unknown $u=\frac{y-x}{2\sqrt{t-s}}$  and for $C,C_1,C_2,C_3,C_4>0$ large enough:
\begin{equation}
\begin{aligned}
&|\big(\int^t_{\frac{t}{2}} e^{(t-s)\D}{\rm div}(v\otimes m)ds\big)(x)|=|\int^t_{\frac{t}{2}} \int_{\R^N}\frac{1}{(4\pi (t-s))^{\N}}e^{-\frac{|x-y|^2}{4(t-s)}}{\rm div}(v\otimes m)(s,y) dy ds|\\
&\leq C\int^t_{\frac{t}{2}} \int_{\R^N}\frac{1}{(4\pi (t-s))^{\frac{N}{2}+1}}|x-y|e^{-\frac{|x-y|^2}{4(t-s)}}|v(s,y)|\,|m(s,y)| dy ds\\
&\leq C_1 \int^t_{\frac{t}{2}} \int_{\R^N}\frac{1}{(4\pi (t-s))^{\frac{N+1}{2}}}\frac{1}{(1+\frac{|x-y|}{2\sqrt{t-s}})^{N+1}}(v(s,y)\otimes m(s,y)) dy ds\\
&\leq C_2 |\int^t_{\frac{t}{2}}\int_{\R^N} \frac{1}{\sqrt{t-s}}\frac{1}{(1+|u|)^{N+1}}|v(s,x+2\sqrt{t-s}\,u)|\,| m(s,x+2\sqrt{t-s}\,u))| du ds\\
&\leq C_3   \|m\|_{{\cal E}_T}\|v\|_{L^\infty_T(L^\infty(\R^N)}   \int^t_{\frac{t}{2}} \frac{1}{\sqrt{t-s}}\frac{1}{\sqrt{s}} ds\\
&\leq C_4  \|m\|_{{\cal E}_T}\|v\|_{L^\infty_T(L^\infty(\R^N)}.
\end{aligned}
\label{kestim1}
\end{equation}
It implies that we have for any $0<t< T$, there exists $C>0$ such that:
\begin{equation}
|\sqrt{t}\int^t_{\frac{t}{2}} e^{(t-s)\D}{\rm div}(v\otimes m)ds|\leq C \sqrt{T} \|m\|_{{\cal E}_T}\|v\|_{L^\infty_T(L^\infty(\R^N)} .
\label{estim1}
\end{equation}
Next in a similar way using the third line of (\ref{kestim1}), we have for $C,C_1
C_2>0$ large enough:
\begin{equation}
\begin{aligned}
&|\big(\int^{\frac{t}{2}}_{0} e^{(t-s)\D}{\rm div}(v\otimes m)ds\big)(x)|
\leq C\int^{\frac{t}{2}}_{0} \int_{\R^N}\frac{1}{(2\sqrt{t-s}+|x-y|)^{N+1}}|v(s,y)|\,| m(s,y)| dy ds\\
&\leq C\int^{\frac{t}{2}}_{0} \int_{\R^N}\frac{1}{(\sqrt{2 t}+|x-y|)^{N+1}}|v(s,y)|\,| m(s,y)| dy ds\\
&\leq C_1\int^{\frac{t}{2}}_{0} \int_{\R^N}\frac{1}{t^{\frac{N+1}{2}}}\frac{1}{(1+\frac{|x-y|}{\sqrt{2t}})^{N+1}}|v(s,y)|\,| m(s,y)| dy ds\\
&\leq C_1\sum_{k\in\mathbb{Z}^N}\int^{\frac{t}{2}}_{0} \int_{x-y\in\sqrt{2t}(k+[0,1]^N)}\frac{1}{t^{\frac{N+1}{2}}}\frac{1}{[\max(1+|k|-\sqrt{N},1)]^{N+1}}|v(s,y)|\, |m(s,y)| dy ds\\
&\leq C_1\sum_{k\in\mathbb{Z}^N}\int^{\frac{t}{2}}_{0} \int_{|y-x+\sqrt{2t}k|\leq \sqrt{2N}\sqrt{t}}\frac{1}{t^{\frac{N+1}{2}}}\frac{1}{[\max(1+|k|-\sqrt{N},1)]^{N+1}}|v(s,y)|\,|m(s,y)| dy ds\\
&\leq C_1\frac{1}{t^{\frac{N+1}{2}}}\sum_{k\in\mathbb{Z}^N}\frac{\|v\|_{L^\infty((0,t)\times\R^N)}}{[\max(1+|k|-\sqrt{N},1)]^{N+1}} \big(\int^{\frac{t}{2}}_{0} \int_{|y-x+\sqrt{2t}k|\leq \sqrt{2N}\sqrt{t}}|m(s,y)|^2 dy ds\big)^{\frac{1}{2}} t^{\frac{N}{4}+\frac{1}{2}}\\
&\leq C_2 \big(\sup_{x\in\R^N,\,0<t<T}t^{-\N}\int^t_0\int_{y\in B(x,\sqrt{t})}|m(s,y)|^2 ds\,dy\big)^{\frac{1}{2}}\|v\|_{L^\infty((0,t)\times\R^N)}.
\end{aligned}
\end{equation}
It implies that we have for any $0<t< T$, there exists $C>0$ such that:
\begin{equation}
|\sqrt{t}\int_0^{\frac{t}{2}} e^{(t-s)\D}{\rm div}(v\otimes m)ds|\leq C\sqrt{T}\|m\|_{{\cal E}_T}\|v\|_{L^\infty_T(L^\infty(\R^N)} .
\label{estim2}
\end{equation}
We now turn to the estimate on the local $L^2$ norm. We are going to deal with the following term:
$$
I=\int^t_0\int_{|x-x_0|\leq\sqrt{t}}|\big(\int^s_0 e^{(s-\tau)\D}{\rm div}(m\otimes v) d\tau\big)(x) |^2 ds dx
$$
We denote by $M$ the following function for $s\in(0,t)$:
$$J(s,x)=\big(\int^s_0 e^{(s-\tau)\D}{\rm div}(m\otimes v) d\tau\big)(x).$$
We now split the function $V$ as follows $J=M_1-M_2-M_3$ with:
$$
\begin{cases}
\begin{aligned}
&M_1=\int^s_0 e^{(s-\tau)\D}{\rm div}((1-\chi_{t,x_0})m\otimes v ) d\tau \\
&M_2=(-\D)^{-\frac{1}{2}}{\rm div}\int^s_0 e^{(s-\tau)\D}\D \big((-\D)^{-\frac{1}{2}}(Id-e^{\sigma\D})(\chi_{t,x_0} m\otimes v ) \big) d\tau\\
&M_3=(-\D)^{-\frac{1}{2}}{\rm div} (-\D)^{\frac{1}{2}}e^{s\D} (\int^s_0 \chi_{t,x_0} m\otimes v \, d\tau).
\end{aligned}
\end{cases}
$$
Here we have use the function $\chi_{t,x_0}(y)=1_{B(x_0,10\sqrt{t})}(y)$.
Next since $s<t$, $|x-x_0|<\sqrt{t}$ and $|y-x_0|\geq 10\sqrt{t}$,  
we have for $C,C_1>0$ large enough ( we have used the fact that $|x_0-y|\leq 2|x-y|$):
$$
\begin{aligned}
&|M_1|\leq C\int_0^s\int_{|y-x_0|\geq 10\sqrt{t}}\frac{1}{(\sqrt{s-\tau}+|x-y|)^{N+1}}|m(\tau,y)||v(\tau,y)| d\tau dy \\
&\leq C_1\int_0^t\int_{|y-x_0|\geq 10\sqrt{t}}\frac{1}{|x_0-y|^{N+1}}|m(\tau,y)||v(\tau,y)| d\tau dy \\
&\leq C\sum_{k\in\mathbb{Z}^N, |k|\geq 4}\int^t_0\int_{y=x_0+\sqrt{t}k+\sqrt{t}[0,1]^N}\frac{1}{(\sqrt{t} \max(10, |k|-\sqrt{N}))^{N+1}}|m(\tau,y)|\,|v(\tau,y)| d\tau dy \\
&\leq \frac{C}{t^{\frac{N+1}{2}}}\|v\|_{L^\infty_s(L^\infty(\R^N)}\sum_{k\in\mathbb{Z}^N, |k|\geq 4} \frac{1}{(\max(10, |k|-\sqrt{N}))^{N+1}} \big(\int^{t}_{0} \int_{|y-x_0-\sqrt{t}k|\leq \sqrt{t}\sqrt{N}}|m(s,y)|^2 dy ds\big)^{\frac{1}{2}} t^{\frac{N}{4}+\frac{1}{2}}\\
&\leq C \|v\|_{L^\infty_T (L^\infty(\R^N)}\|m\|_{{\cal E}_T}.
\end{aligned}
$$
It implies that:
$$
\begin{aligned}
\int^t_0 \int_{|x-x_0|\leq\sqrt{t}} |M_1(s,x)|^2 ds \, dx\leq C t^{\N+1}\|v\|^2_{L^\infty_T (L^\infty(\R^N)}\|m\|^2_{{\cal E}_T}.
\end{aligned}
$$
We deduce that for any $0<t<T$:
\begin{equation}
\begin{aligned}
&\big( \frac{1}{t^{\N}} \int^t_0 \int_{|x-x_0|\leq\sqrt{t}} |M_1(s,x)|^2 ds \, dx\big)^{\frac{1}{2}}\leq C\sqrt{T}\|v\|_{L^\infty_T (L^\infty(\R^N)}\|m\|_{{\cal E}_T}.
\end{aligned}
\label{estim3}
\end{equation}
To estimate $M_2$ and $M_3$ as in \cite{Lemarie}, we note $M(\tau,y)=\chi_{t,x_0}(y)m(\tau,y)\otimes v(\tau,y)$. Using the boundness of the Riez transform on $L^2(\R^N)$, we obtain since $M_2$ is the sum of terms of the form $R_i \int^s_0 e^{(s-\tau)\D}\D ((-\D)^{-\frac{1}{2}}(Id-e^{\sigma\D})(\chi_{t,x_0} m_k v_j ) \big) d\tau$ with $R_i$ the Riesz transform:
$$
\int^t_0\int_{\R^N}|M_2|^2 ds dx\leq C \int^t_0\int_{\R^N}|\int^s_0 e^{(s-\tau)\D}\D(-\D)^{-\frac{1}{2}}(Id-e^{\tau\D})M\,d\tau|^2 ds \, dx
$$
Now using the maximal $L^2_t(L^2_x)$ regularity of the heat kernel, we get:
$$\int^t_0\int_{\R^N}|M_2|^2 ds dx\leq C  \int^t_0\int_{\R^N}|(-\D)^{-\frac{1}{2}}(Id-e^{\tau\D})M|^2 d\tau \, dx
$$
Since the function $z\rightarrow \frac{1-e^{-\tau z^2}}{\sqrt{\tau}z}$ is uniformly bounded on $\R^+$ in $\tau$, we deduce using Plancherel theorem that for $C>0$ large enough we have:
\begin{equation}
\begin{aligned}
\int^t_0\int_{\R^N}|M_2|^2 ds dx&\leq C  \int^t_0\int_{\R^N}\tau |M|^2 d\tau \, dx\\
&\leq C\|\tau M\|_{L^\infty((0,t)\times\R^N)} \|M\|_{L^1((0,t)\times\R^N)}.
\end{aligned}
\label{1estim}
\end{equation}
Next we have using the definition of $M$ and for $0<t<T$:
\begin{equation}
\begin{cases}
\begin{aligned}
&\|\tau M\|_{L^\infty((0,t)\times\R^N)}\leq \sqrt{t}\|m\|_{{\cal E}_T}\|v\|_{L^\infty((0,T)\times\R^N)}\\
&\|M\|_{L^1((0,t)\times\R^N)}\leq (\int^t_0\int_{|x-x_0|\leq10\sqrt{t}}|m(s,x)|^2 dx\, ds)^{\frac{1}{2}}t^{\frac{N}{4}+\frac{1}{2}}\,\|v\|_{L^\infty((0,T)\times\R^N)} \\
\end{aligned}
\end{cases}
\label{2estim}
\end{equation}
From (\ref{1estim}) and (\ref{2estim}) we deduce that there exists $C>0$ such that:
\begin{equation}
\begin{aligned}
t^{-\N}\int^t_0\int_{\R^N}|M_2|^2 ds dx&\leq C 
t \|m\|_{{\cal E}_T}\|v\|^2_{L^\infty((0,T)\times\R^N)}(\int^t_0\int_{|x-x_0|\leq10\sqrt{t}}|m(s,x)|^2 dx\, ds)^{\frac{1}{2}}t^{-\frac{N}{4}}\\
&\leq C 
t \|m\|^2_{{\cal E}_T}\|v\|^2_{L^\infty((0,T)\times\R^N)}.
\end{aligned}
\label{3estim}
\end{equation}
It implies that we have for $C>0$ since $t\in(0,T)$:
\begin{equation}
\begin{aligned}
(t^{-\N}\int^t_0\int_{\R^N}|M_2|^2 ds dx)^{\frac{1}{2}}&\leq C 
\sqrt{T} \|m\|_{{\cal E}_T}\|v\|_{L^\infty((0,T)\times\R^N)}.
\end{aligned}
\label{4estim}
\end{equation}
Similarly we have now for $C>0$ large enough:
$$
\int^t_0\int_{\R^N}|M_3|^2 ds dx\leq C\int^t_0\int_{\R^N}|(-\D)^{\frac{1}{2}}e^{s\D}(\int^s_0 M d\tau)|^2 ds\, dx.
$$
We make the change of variables $s=t s_1$, $\tau=t\tau_1$, $x=\sqrt{t} x_1$, $y=\sqrt{t}y_1$ and we have for $C>0$ large enough:
\begin{equation}
\int^t_0\int_{\R^N}|M_3|^2 ds dx\leq C t^{2+\N} \int^1_0\int_{\R^N}|(-\D)^{\frac{1}{2}}e^{s_1\D}(\int^{s_1}_0 N d\tau_1)|^2 ds_1\, dx_1,
\label{5estim}
\end{equation}
with $N(\tau_1,y_1)=M(t\tau_1,\sqrt{t} y_1)$. We have in particular used the fact that the operator $\sqrt{-\D}$ has a kernel wich is homogeneous of degree $-N-1$. We recall now the following Lemma (see \cite{Lemarie} p 163).
\begin{lem}
For $\alpha$ a function defined on $(0,1)\times\R^N$, let $A(\alpha)$ defined as follows:
$$A(\alpha)=\sup_{x_0\in\R^N, 0<t<1} t^{-\N}\int^t_0\int_{|x-x_0|<\sqrt{t}}|\alpha(s,x)| ds dx $$
Then for $\beta(t,x)=(-\D)^{\frac{1}{2}}e^{t\D}\int^t_0\alpha ds$, the following inequality holds:
$$\int^1_0\int_{\R^N} |\beta(s,x)|^2 ds\, dx\leq C A(\alpha) \int^1_0\int_{\R^N}|\alpha(s,x)| ds\, dx.$$
\label{lemlema}
\end{lem}
Using the Lemma \ref{lemlema}, we get then:
\begin{equation}
\int^t_0\int_{\R^N}|M_3|^2 dx\, ds\leq C t^{2+\N} A(N)\|N\|_{L^1((0,1)\times\R^N)}.
\label{M3a}
\end{equation}
Next we have for any $0<t<T$:
\begin{equation}
\begin{aligned}
&\|N\|_{L^1((0,1)\times\R^N)}=\int^1_0\int_{\R^N}\chi_{t,x_0}(\sqrt{t} y_1) |m(t \tau_1,\sqrt{t} y_1) v(t\tau_1,\sqrt{t} y_1) | d\tau_1 dy_1\\
&\leq\frac{1}{t^{\frac{N}{2}+1}} \int^t_0\int_{|y-x_0|\leq10\sqrt{t}}|m(u,y)v(u,y)| du dy\\
&\leq \frac{1}{t^{\N+1}} \big(\int^t_0\int_{|y-x_0|\leq10\sqrt{t}}|m(u,y)|^2 du dy\big)^{\frac{1}{2}}\|v\|_{L^\infty((0,t)\times\R^N)}t^{\frac{N}{4}+\frac{1}{2}}\\
&\leq \frac{1}{\sqrt{t}}\|m\|_{{\cal E}_T}\|v\|_{L^\infty((0,T),L^\infty(\R^N))}.
\end{aligned}
\label{N1}
\end{equation}
Furthermore by definition, we have:
\begin{equation}
\begin{aligned}
&\|A(N)\|=\sup_{x_0\in\R^N, 0<s<1}s^{-\N}\int^s_0\int_{|y_1-x_0|<\sqrt{s}}\chi_{t,x_0}(\sqrt{t} y_1) |m(t \tau_1,\sqrt{t} y_1) v(t\tau_1,\sqrt{t} y_1) | d\tau_1 dy_1\\
&\leq \sup_{x_0\in\R^N, 0<s<1}t^{-1}(ts)^{-\N}\int^{ts}_0\int_{|y-x_0\sqrt{t}|<\sqrt{s t}}\chi_{t,x_0}(y) |m(u,y) v(u,y) | du dy\\
&\leq \sup_{x_0\in\R^N, 0<s<1}t^{-1}(ts)^{-\frac{N}{4}}\big(\int^{ts}_0\int_{|y-x_0\sqrt{t}|<\sqrt{s t}} |m(u,y) |^2 du dy\big)^{\frac{1}{2}}(t s)^{\frac{1}{2}}\|v\|_{L^\infty((0,t)\times\R^N)}\\
&\leq \frac{1}{\sqrt{t}}\|m\|_{{\cal E}_T}\|v\|_{L^\infty((0,T),L^\infty(\R^N))}.
\end{aligned}
\label{N2}
\end{equation}
Combining (\ref{M3a}), (\ref{N1}) and (\ref{N2}) we deduce that there exists $C>0$ such that for any $0<t<T$:
\begin{equation}
\big(\frac{1}{t^\N}\int^t_0\int_{\R^N}|M_3|^2 dx\, ds\big)^{\frac{1}{2}}\leq C\sqrt{T} \|m\|_{{\cal E}_T}\|v\|_{L^\infty((0,T),L^\infty(\R^N))}.
\label{M3ab}
\end{equation}
Combining (\ref{estim1}), (\ref{estim2}), (\ref{estim3}), (\ref{4estim}) and (\ref{M3ab}) we obtain (\ref{cru}). We are now interested in estimating the term:
$$B_1(\rho)=\int^t_0e^{(t-s)\D}\n P(\rho) ds.$$
In a similar way as in \cite{Koch}, we know that there exists $C>0$ such that for any $T>0$ we have:
\begin{equation}
\begin{aligned}
&\|B_1(\rho)\|_{{\cal E}_T}\leq \|P(\rho)\|_{Y_T}
\end{aligned}
\label{cru1}
\end{equation}
with:
$$\|P(\rho)\|_{Y}=\sup_{0<t<T}t\|P(\rho)(t,\cdot)\|_{L^\infty(\R^N)}+\sup_{x\in\R^N,\,0<t<T}t^{-\N}\int^t_0\int_{y\in B(x,\sqrt{t})}|P(\rho)(t,y)| dt \,dy.$$
It implies then that:
\begin{equation}
\begin{aligned}
&\|B_1(\rho)\|_{{\cal E}_T}\leq T \|P(\rho)\|_{L^\infty((0,T),L^\infty(\R^N))}.
\end{aligned}
\label{cru2}
\end{equation}
Combining (\ref{Duha}), (\ref{cru}) and (\ref{cru2}), we deduce that there exists $C>0$ large enough such that:
\begin{equation}
\|m\|_{{\cal E}_T}\leq C\|m_0\|_{bmo^{-1}}+C\sqrt{T}\big(\|m\|_{{\cal E}_T}\|v\|_{L^\infty((0,T),L^\infty(\R^N))}+\sqrt{T}\|P(\rho)\|_{L^\infty((0,T),L^\infty(\R^N))}\big).
\label{controle1}
\end{equation}
\subsubsection*{$L^\infty$ estimates on the effective velocity $v$}
We recall that the effective velocity $v$ satisfies (since ${\rm curl}v=0$ on $\R^+\times\R^N$):
$$\rho \p_{t}v+\rho u\cdot\n v+\frac{P'(\rho)}{2\mu}\rho v=\frac{P'(\rho)}{2\mu}\rho u.$$
Dividing by $\rho$ the previous equation (it has a sense since in fact we consider the solution $(\rho_n,u_n,v_n)$ which is regular and which satisfies $\frac{1}{\rho_n}\in L^\infty(\R^+,L^\infty(\R^N))$) , we have:
$$\p_{t}v+ u\cdot\n v+\frac{P'(\rho)}{2\mu} v=\frac{P'(\rho)}{2\mu\rho} m.$$
Since we have a damped transport equation, we deduce that for any $t\in\R^+$ we have for $C>0$ large enough:
$$\|v\|_{L^\infty((0,t),L^\infty(\R^N)}\leq\|v_0\|_{L^\infty(\R^N)}+ C\int^t_0 \|\frac{P'(\rho)}{\rho}(s,\cdot)\|_{L^\infty(\R^N)}\|m(s,\cdot)\|_{L^\infty(\R^N))} ds.$$
Using (\ref{Pressiontech}) and since $\gamma\geq 2$, we obtain then that for any $T>0$ :
\begin{equation}
\|v\|_{L^\infty((0,T),L^\infty(\R^N)}\leq \|v_0\|_{L^\infty(\R^N)}+C\sqrt{T} \|\rho\|_{L^\infty((0,T),L^\infty(\R^N))}^{\gamma-2}\|m\|_{{\cal E}_T}.
\label{controle2}
\end{equation}
\subsubsection*{Estimates in Besov spaces on $\rho u$}
From the equation on $\rho u$ in (\ref{0.3}) we deduce using proposition \ref{chaleur} that there exists $C>0$ such that for any $0<T\leq 1$ we have:
\begin{equation}
\begin{aligned}
&\|\rho u\|_{\widetilde{L}^\infty_T(B^{-1}_{\infty,\infty})}+\|\rho u\|_{L^1_T(\dot{B}^{1}_{\infty,\infty})}\leq C(\|m_0\|_{B^{-1}_{\infty,\infty}}+\|\rho u\otimes v\|_{\widetilde{L}^1_T(B^{0}_{\infty,\infty})}+T\|P(\rho)\|_{L^\infty_T(L^\infty)}).
\end{aligned}
\label{mom}
\end{equation}
It is important to mention that the left hand side of (\ref{mom}) is well defined since $\rho u$ belongs to $\widetilde{L}^\infty_T(\dot{B}^{s}_{\infty,\infty})$ for any $s\geq 0$ and $\rho u$ is in $L^\infty_T(L^2)$ then by Besov embedding in $\widetilde{L}^\infty(\dot{B}^{-\N}_{\infty,\infty})$. In addition $m_0$ belongs to $B^{-1}_{\infty,\infty}(\R^N)$ since $m_0=\rho_0 v_0+2\mu\n\rho_0$ with $\rho_0\in L^\infty(\R^N)$ and $\rho_0v_0\in L^\infty(\R^N)$.\\
We are going now to estimate the term $ \|\rho u\otimes v\|_{\widetilde{L}^1_T(B^{0}_{\infty,\infty})}$.
It comes using paraproduct laws:
\begin{equation}
\begin{aligned}
&\|T_v (\rho u)\|_{\widetilde{L}^2_T(B^{0}_{\infty,\infty})}\leq \|v\|_{L^\infty_T(L^\infty)}\|\rho u\|_{\widetilde{L}^2_T(B^{0}_{\infty,\infty})}.
\end{aligned}
\label{para1}
\end{equation}
Using again paraproduct law, we deduce from Besov embedding that we have for $C,C_1,C_2>0$ large enough:
\begin{equation}
\begin{aligned}
\|T_{\rho u} v\|_{\widetilde{L}^{2-\alpha(\e)}_T(B^{0}_{\infty,\infty})}&\leq C\|v\|_{\widetilde{L}_T^\infty(B^{0}_{\infty,\infty})}\|\rho u\|_{L^{2-\alpha(\e)}_T(L^\infty)}, \\
&\leq C_1 \|v\|_{L^\infty_T(L^\infty)}\|\rho u\|_{L^{2-\alpha(\e)}_T(B^{0}_{\infty,1})},\\
&\leq C_2\|v\|_{L^\infty_T(L^\infty)}\|\rho u\|_{\widetilde{L}^{2-\alpha(\e)}_T(B^{\e}_{\infty,\infty})},
\end{aligned}
\label{para4}
\end{equation}
with $\alpha(\e)=\frac{2\e}{1+\e}$
Finally we have (see \cite{Danchin}):
\begin{equation}
\begin{aligned}
&\|R(\rho u,v) \|_{\widetilde{L}^{2-\alpha(\e)}_T(B^{0}_{\infty,\infty})}
\leq \|\rho u \|_{\widetilde{L}^{2-\alpha(\e)}_T(B^{0}_{\infty,1})
} \|v\|_{\widetilde{L}^\infty_T(B^0_{\infty,\infty})}\\
&\leq \|v\|_{L^\infty_T(L^\infty)}\|\rho u\|_{\widetilde{L}^{2-\alpha(\e)}_T(B^{\e}_{\infty,\infty})}
\end{aligned}
\label{para5}
\end{equation}
Combining (\ref{mom}), (\ref{para1}), (\ref{para4}) and (\ref{para5}) we have using interpolation and Young inequality for $C>0$ large enough and $\beta>0$, then we get for any $0<T<1$:
\begin{equation}
\begin{aligned}
&\|\rho u\|_{\widetilde{L}^\infty_T(B^{-1}_{\infty,\infty})}+\|\rho u\|_{\widetilde{L}^1_T(B^{1}_{\infty,\infty})}\leq C\big(\|m_0\|_{B^{-1}_{\infty,\infty}}+T\|P(\rho)\|_{L^\infty_T(L^\infty)}\\
&+T^\beta \|v\|_{L^\infty_T(L^\infty(\R^N))}(\|\rho u\|_{\widetilde{L}^\infty_T(B^{-1}_{\infty,\infty})}+\|\rho u\|_{\widetilde{L}^1_T(B^{1}_{\infty,\infty})})\big).
\end{aligned}
\label{momfin}
\end{equation}
\subsubsection*{$L^\infty$ estimates on the density $\rho$}
We recall that the density $\rho$ verifies on $\R^+\times\R^N$:
$$\p_t\rho-2\mu\D\rho+{\rm div}(\rho v)=0.$$
Using the Duhamel formula, we have then:
$$\rho(t,x)=e^{2\mu t\D}\rho_0+\int^t_0e^{2\mu(t-s)\D}{\rm div}(\rho v) ds.$$
Using the maximum principle, it yields for $C>0$ and $C_1>0$ large enough:
$$
\begin{aligned}
&|\rho(t,x)|=|e^{2\mu t\D}\rho_0(x)+C\int^t_0\int_{\R^N}\sum_{i+1}^N\frac{x_i-y_i}{\sqrt{4(t-s)}}e^{\frac{-|x-y|^2}{4(t-s)}}\frac{1}{(t-s)^{\frac{N+1}{2}}}\rho(s,y)v_i(s,y) ds dy|\\
&\leq \|\rho_0\|_{L^\infty}+C_1\int^t_0\frac{1}{\sqrt{t-s}}\|\rho(s,\cdot)\|_{L^\infty}\|v(s,\cdot)\|_{L^\infty}ds
\end{aligned}
$$
We deduce that for any $(t,x)\in\R^+\times\R^N$ we have:
$$
\begin{aligned}
&\rho(t,x)\leq \|\rho_0\|_{L^\infty}+C\int^t_0\frac{1}{\sqrt{t-s}}\|\rho(s,\cdot)\|_{L^\infty}\|v(s,\cdot)\|_{L^\infty}ds
\end{aligned}
$$
It implies using Gronwall Lemma that for any $T>0$ we have for $C>0$ large enough:
\begin{equation}
\|\rho\|_{L^\infty_T(L^\infty(\R^N))}\leq C\|\rho_0\|_{L^\infty(\R^N)}e^{C\sqrt{T}\|v\|_{L^\infty_T(L^\infty(\R^N))}}.
\label{rhoLinf}
\end{equation}
\subsection*{Existence in finite time of a solution for the system (\ref{0.1m})}
In the previous section we have defined a sequence of solution $(\rho_n,u_n,v_n)_{n\in\mathbb{N}}$ which is a unique global solution of the system (\ref{0.1m}) for the initial data $(\rho_0^n,u_0^n,v_0^n)$. We define now the following norm $\|\cdot\|_{Z_T}$ for $T>0$ as follows (in fact $Z_T=Y_T$):
\begin{equation}
\begin{aligned}
\|(\rho,v,m)\|_{Z_T}=&\|\rho\|_{L^\infty_T(L^\infty(\R^N))}+\|v\|_{L^\infty_T(L^\infty(\R^N))}+\|m\|_{{\cal E}_T}\\
&+\|m\|_{\widetilde{L}^\infty_T(B^{-1}_{\infty,\infty})}+\|m\|_{\widetilde{L}^1_T(B^{1}_{\infty,\infty})}.\\[2mm]
\|(\rho,v,m)\|_{Z_{\infty}}=&\|\rho\|_{L^\infty(\R^+,L^\infty(\R^N))}+\|v\|_{L^\infty(\R^+,L^\infty(\R^N))}+\|m\|_{{\cal E}_{\infty}}\\
&+\|m\|_{\widetilde{L}^\infty(\R^+,B^{-1}_{\infty,\infty})}+\|m\|_{\widetilde{L}^1(\R^+,B^{1}_{\infty,\infty})}.
\end{aligned}
\end{equation}
We recall that since  $(\rho_n,u_n,v_n)_{n\in\mathbb{N}}$ is a regular solution, $\|(\rho_n,v_n,m_n)\|_{Z_T}$ is well defined for any $0<T\leq 1$ and any $n\in\mathbb{N}$. We deduce from (\ref{controle2}), (\ref{controle1}) and (\ref{momfin}) that there exists $\beta>0$, $C>0$ such that for any $0<T<1$ we have:
\begin{equation}
\begin{aligned}
&\|(\rho_n,v_n,m_n)\|_{Z_T} \leq  \|\rho_n\|_{L^\infty_T(L^\infty(\R^N))} +C\biggl(\|v^n_0\|_{L^\infty(\R^N)}+C\sqrt{T} \|\rho_n\|_{L^\infty((0,T),L^\infty(\R^N))}^{\gamma-2}\|m_n\|_{{\cal E}_T}\\
&+ \|m^n_0\|_{bmo^{-1}}+C\sqrt{T}\big(\|m_n\|_{{\cal E}_T}\|v_n\|_{L^\infty((0,T),L^\infty(\R^N))}+\sqrt{T}\|P_n(\rho_n)\|_{L^\infty((0,T),L^\infty(\R^N))}\big)\\
&+\|\rho^n_0 u^n_0\|_{B^{-1}_{\infty,\infty}}+T\|P_n(\rho_n)\|_{L^\infty_T(L^\infty)}\\
&\hspace{3cm}+T^\beta \|v_n\|_{L^\infty_T(L^\infty(\R^N))}(\|\rho_n u_n\|_{\widetilde{L}^\infty_T(B^{-1}_{\infty,\infty})}+\|\rho_n u_n\|_{\widetilde{L}^1_T(B^{1}_{\infty,\infty})})\biggl).
\end{aligned}
\end{equation}
It implies that we have using again (\ref{rhoLinf}) and (\ref{Pressiontech}) for $C>0$ large enough:
\begin{equation}
\begin{aligned}
&\|(\rho_n,v_n,m_n)\|_{Z_T} \leq C\|\rho^n_0\|_{L^\infty(\R^N)}e^{C\sqrt{T}\|(\rho_n,v_n,m_n)\|_{Z_T}}+C\big(\|v^n_0\|_{L^\infty(\R^N)}\\
&+C\sqrt{T}\|\rho^n_0\|_{L^\infty(\R^N)}^{\gamma-2}e^{C(\gamma-2)\sqrt{T}\|(\rho_n,v_n,m_n)\|_{Z_T}}\|(\rho_n,v_n,m_n)\|_{Z_T}+ \|m^n_0\|_{bmo^{-1}}\\
&+C\sqrt{T}\big(\|(\rho_n,v_n,m_n)\|^2_{Z_T} +C\sqrt{T}\|\rho^n_0\|_{L^\infty(\R^N)}^{\gamma}e^{C\gamma\sqrt{T}\|(\rho_n,v_n,m_n)\|_{Z_T}}\big)\\
&+\|\rho_0^n u_0^n\|_{B^{-1}_{\infty,\infty}}+CT \|\rho_0^n\|_{L^\infty(\R^N)}^{\gamma} e^{C\gamma\sqrt{T}\|(\rho_n,v_n,m_n)\|_{Z_T}}+T^\beta \|(\rho_n,v_n,m_n)\|^2_{Z_T}\big).
\end{aligned}
\label{supcrua1}
\end{equation}
We denote by $T_n^*$ the following time such that:
\begin{equation}
\begin{aligned}
&T_n^*=\sup\{t>0, \forall s\in[0,t]\;\|(\rho_n,v_n,m_n)\|_{Z_s}< C_1 C(\|\rho_0\|_{L^\infty(\R^N)}+
\|\rho_0 v_0\|_{L^\infty(\R^N)}\\
&+\|\frac{1}{\rho_0}\|_{L^\infty(\R^N)}\|\rho_0 v_0\|_{L^\infty(\R^N)}(1+\|\rho_0\|_{L^\infty(\R^N)}))
\end{aligned}
\label{Tmax}
\end{equation}
with $C_1>3$ large enough.
In the sequel we will note $M_1=C(\|\rho_0\|_{L^\infty(\R^N)}+
\|\rho_0 v_0\|_{L^\infty(\R^N)}
+\|\frac{1}{\rho_0}\|_{L^\infty(\R^N)}\|\rho_0 v_0\|_{L^\infty(\R^N)}(1+\|\rho_0\|_{L^\infty(\R^N)}))$. Let us prove now that for any $n\in\mathbb{N}$, we have $T_n^*>0$. First it is clear by definition of $(\rho_0^n,m_0^n,m_{0,1}^n,v_0^n)$ (see (\ref{definitu})) that we have:
\begin{equation}
\begin{cases}
\begin{aligned}
&\|\rho_0^n\|_{L^\infty}\leq\|\rho_0\|_{L^\infty},\;\|m_{0,1}^n\|_{L^\infty}\leq \|\rho_0v_0\|_{L^\infty}\\
&\|m^n_0\|_{bmo^{-1}}\leq 2\mu\|\rho_0\|_{L^\infty}+\|\rho_0 v_0\|_{L^\infty}\\
&\|m^n_0\|_{B^{-1}_{\infty,\infty}}\leq 2\mu\|\rho_0\|_{L^\infty}+\|\rho_0 v_0\|_{L^\infty}\\
&\|v^n_0\|_{L^\infty}\leq \|\frac{1}{\rho_0}\|_{L^\infty}\|\rho_0v_0\|_{L^\infty}.
\end{aligned}
\end{cases}
\label{estimop}
\end{equation}
Now it is easy to observe that for any $t>0$ we have for $C_2>0$ large enough:
$$\|\rho_n(t,\cdot)\|_{L^\infty}\leq \|\rho^n_0\|_{L^\infty}e^{C_2\int^t_0\|\n u_n(s,\cdot)\|_{L^\infty} ds}. $$
Since $(\rho_n,m_n,v_n)$ is a global strong solution, we know that $u_n\in \widetilde{L}_{loc}^{\frac{2}{2-\e}}(\R^+,\dot{B}^{\NN+1}_{p,1})$. By Besov embedding, we get for a constant $C(n,t)$ large enough and depending on $n$ and $t$ that:
\begin{equation}
\|\rho_n(t,\cdot)\|_{L^\infty}\leq \|\rho_0\|_{L^\infty}\exp(t^{\frac{\e}{2}} \|\n u_n\|_{L_t^{\frac{2}{2-\e}}(\dot{B}^\NN_{p,1})} ). 
\label{gutech}
\end{equation}
It implies in particular that for $t^1_n>0$ sufficiently small and depending on $n$, we have from (\ref{estimop}) and (\ref{gutech}) for $t\in(0,t^1_n]$:
\begin{equation}
\|\rho_n(t,\cdot)\|\leq 2C\|\rho_0\|_{L^\infty}\leq 2C \|\frac{1}{\rho_0}\|_{L^\infty}\|\rho_0v_0\|_{L^\infty}.
\label{t1}
\end{equation}
We can observe now that:
$$\|m_n\|_{{\cal E}_t}\rightarrow_{t\rightarrow 0}0.$$
Indeed since $(m_n)_{n\in\mathbb{N}}$ satiesfies (\ref{estimo3}), it implies in particular that for $t>0$ and $C_2>0$ large enough:
\begin{equation}
\begin{aligned}
&\|m_n\|_{{\cal E}_t}\leq C_2 \sqrt{t}\|m_n\|_{L^\infty([0,t],L^\infty(\R^N))}.
\end{aligned}
\label{t3}
\end{equation}
Combining (\ref{estimop}), (\ref{t1}), (\ref{t3}) and (\ref{controle2}), we can show that for $t^2_n>0$ sufficiently small, we have for any $t\in(0,t^2_n)$:
\begin{equation}
\|v_n(t,\cdot)\|\leq 2C\|v_0^n\|_{L^\infty}.
\label{t2}
\end{equation}
Finally we have for $t\in(0,\min(t^1_n,t^2_n))$:
\begin{equation}
\|m_n(t,\cdot)\|_{B^{-1}_{\infty,\infty}}\leq \|\rho_n v_n(t,\cdot)\|_{L^\infty}+2\mu \|\rho_n(t,\cdot)\|_{L^\infty}\leq C_2M_1,
\label{t2a}
\end{equation}
for $C_2>0$ large enough.
As previously, we can show easily that:
\begin{equation}
\|m_n\|_{L^1_t(B^{1}_{\infty,\infty})}\rightarrow_{t\rightarrow0}0.
\label{t2ab}
\end{equation}
From (\ref{t1}), (\ref{t3}), (\ref{t2}), (\ref{t2a}) and (\ref{t2ab}), we deduce that  for any $n\in\mathbb{N}$ we have $T_n^*>0$ provided that $C_1>0$ is large enough. We define now in the sequel $T$ as follows:
\begin{equation}
\begin{cases}
\begin{aligned}
&0<T\leq 1\\
&e^{C_1C\sqrt{T}M_1}\leq\frac{6}{5}\\
&CC_1 M_1\sqrt{T} \|\rho_0\|_{L^\infty(\R^N)}^{\gamma-2}e^{C_1C(\gamma-2)\sqrt{T}M_1}\leq \frac{1}{20 C}M_1\\
& C\sqrt{T} C_1^2 M_1^2\leq\frac{1}{20C}M_1\\
&C^2 T \|\rho_0\|_{L^\infty(\R^N)}^{\gamma}e^{C_1C\gamma\sqrt{T}M_1}\leq\frac{1}{20C}M_1\\
&CT  \|\rho_0\|_{L^\infty(\R^N)}^{\gamma} e^{C_1 C \gamma\sqrt{T}M_1}\leq\frac{1}{20C}M_1\\
&T^\beta C_1^2 M_1^2\leq\frac{1}{20C}M_1.
\end{aligned}
\end{cases}
\label{defM1}
\end{equation}
Assume by absurd that $T_n^*<T$, it implies in particular from (\ref{supcrua1}), (\ref{estimop}) and (\ref{defM1}) that:
\begin{equation}
\begin{aligned}
&\|(\rho_n,v_n,m_n)\|_{Z_{T_n^*}} \leq \frac{6}{5}C\|\rho^n_0\|_{L^\infty(\R^N)}+C\big(\|v^n_0\|_{L^\infty(\R^N)}\\
&+\frac{1}{20 C }M_1+ \|m_0^n\|_{bmo^{-1}}+\frac{1}{20C}M_1+\frac{1}{20 C }M_1+\|m_0^n\|_{B^{-1}_{\infty,\infty}}\\
&+\frac{1}{20 C }M_1+\frac{1}{20 C }M_1\big).
\end{aligned}
\label{supcrua2}
\end{equation}
From (\ref{estimop}) and (\ref{supcrua2}) we have:
\begin{equation}
\begin{aligned}
\|(\rho_n,v_n,m_n)\|_{Z_{T_n^*}} &\leq \frac{6}{5}M_1+\frac{1}{4}M_1+C(\|\frac{1}{\rho_0}\|_{L^\infty}\|\rho_0\|_{L^\infty}\|v_0\|_{L^\infty}+4\mu\|\rho_0\|_{L^\infty}+\|\rho_0\|_{L^\infty}\|v_0\|_{L^\infty})\\
&<\frac{C_1}{2}M_1,
\end{aligned}
\label{supcrua3}
\end{equation}
provided that $C_1>3$ is large enough. Following the same argument than previously we can show that:
\begin{equation}
\|(\rho_n,v_n,m_n)\|_{Z_t}=C^1_n(t),
\label{impbootad}
\end{equation}
with $C^1_n$ depending on $n$ and continuous. It is a consequence of the fact that $(m_n,\rho_n-\frac{1}{n})$ belong to $C([0,T],\dot{B}^{\NN}_{p,1})$ for any $T>0$ (because we have a sequence of regular solutions). It implies then that the function $t\rightarrow \sqrt{t}\|m_n(t,\cdot)\|_{L^\infty}$ is continuous. Setting $A_n(s)=\sup_{0<t\leq s} \sqrt{t}\|m_n(t,\cdot)\|_{L^\infty}$, we can observe that $A_n$ is continuous on $(0,+\infty)$. Indeed let $T>0$ and $\e>0$ small enough, we have $A_n(T+\e)=A_n(T)$ or $A_n(T+\e)>A_n(T)$. The last case implies that:
$$A_n(T+\e)=\sup_{t\in(T,T+\e]}\sqrt{t}\|m_n(t,\cdot)\|_{L^\infty}\rightarrow_{\e\rightarrow0}\sqrt{T}\|m_n(T,\cdot)\|_{L^\infty}\leq A_n(T).$$
Similarly the function $B_n(T)=\sup_{x_0\in\R^N,0<s\leq T}\frac{1}{s^\N}\int^s_0\int_{B(x_0,\sqrt{s})}|m_n(s',x)|^2 ds' dx$ is continuous. Again for $\e>0$ we have for $T>0$, $B_n(T+\e)=B_n(T)$ or $B_n(T+\e)>B_n(T)$. In this last case we have:
$$B_n(T+\e)=\sup_{x_0\in\R^N,T<s\leq T+\e}\frac{1}{s^\N}\int^s_0\int_{B(x_0,\sqrt{s})}|m_n(s',x)|^2 ds' dx$$
For $s\in(T,T+\e]$, we have:
$$
\begin{aligned}
&|\frac{1}{s^\N}\int^s_0\int_{B(x_0,\sqrt{s})}|m_n(s',x)|^2 ds' dx-\frac{1}{T^\N}\int^T_0\int_{B(x_0,\sqrt{T})}|m_n(s',x)|^2 ds' dx|\\
&\leq|(\frac{1}{s^\N}-\frac{1}{T^\N})\int^T_0\int_{B(x_0,\sqrt{T})}|m_n(s',x)|^2 dx ds'|+
|\frac{1}{s^\N}\int^T_0\int_{B(x_0,\sqrt{s})\backslash B(x_0,\sqrt{T})}|m_n(s,x)|^2 dxds|\\
&+|\frac{1}{s^\N}\int^s_T\int_{B(x_0,\sqrt{s})}|m_n(s,x)|^2 dxds|
\end{aligned}
$$
Using the fact that $m_n$ is bounded in $L^\infty([0,T],L^\infty\cap L^2)$ for any $T>0$, we deduce then that:
$$
\begin{aligned}
&|\frac{1}{s^\N}\int^s_0\int_{B(x_0,\sqrt{s})}|m_n(s',x)|^2 ds' dx-\frac{1}{T^\N}\int^T_0\int_{B(x_0,\sqrt{T})}|m_n(s',x)|^2 ds' dx|\leq C(\e,n).
\end{aligned}
$$
with $C(\e,n)$ depending on $n$ and $\e$ and with $C(\e,n) \rightarrow_{\e\rightarrow0}0$. It proves that $t\rightarrow \|m_n\|_{{\cal E}_t}$ is continuous. Similarly $t\rightarrow \|\rho_n\|_{L^\infty_t(L^\infty)}$, $t\rightarrow \|v_n\|_{L^\infty_t(L^\infty)}$, $t\rightarrow \|m_n\|_{L^\infty_t(L^\infty)}$ and $t\rightarrow \|\n m_n\|_{L^\infty_t(L^\infty)}$ is continuous since $(\rho_n,m_n,v_n)$ is a regular solution on $(0,+\infty)$. It implies then that $t\rightarrow \|(\rho_n,v_n,m_n)\|_{Z_t}$ is continuous on $(0,+\infty)$.\\
From (\ref{supcrua3}), (\ref{impbootad}) and the definition of $T_n^*$, we deduce that we have an absurdity. In other words we have proved that for any $n\in\mathbb{N}$, we have:
\begin{equation}
\begin{cases}
\begin{aligned}
&T_n^*\geq T\\
&\|(\rho_n,v_n,m_n)\|_{Z_{T_n^*}} \leq C_1 M_1.
\end{aligned}
\end{cases}
\label{infofin}
\end{equation}
\begin{remarka}
Let us assume now that $\frac{1}{\rho_0}$ is not in $L^\infty(\R^N)$. We only assume that $v_0$ is in $L^1(\R^N)$. It suffices then to apply the same ideas as previously, we have only to change our regularizing process on the initial data. More precisely, we set:
\begin{equation}
\begin{cases}
\begin{aligned}
&\rho_0^n=\varphi_n j_n*\rho_0+\frac{1}{n}, v_0^n=\varphi_n j_n*v_0\\
&m_{1,0}^n=\rho_0^n v_0^n,\;m_0^n=m_{0,1}^n-2\mu\n\rho_0^n\\
&u_0^n=\frac{m_0^n}{\rho_0^n}.
\end{aligned}
\end{cases}
\label{definitua1}
\end{equation}
with $j_n$ and $\varphi_n$ defined as previously. We can check that :
\begin{equation}
\begin{cases}
\begin{aligned}
&\rho_0^n\rightarrow _{n\rightarrow+\infty}\rho_0\;\;\mbox{weakly * in}\;L^\infty(\R^N),\\
&v_0^n\rightarrow _{n\rightarrow+\infty}v_0\;\;\mbox{strongly in}\;L^1(\R^N).
\end{aligned}
\end{cases}
\label{confaible6}
\end{equation}
Similarly for any $\psi\in C_0^\infty(\R^N)$ we have:
\begin{equation}
\int_{\R^N}\rho_0^n v_0^n\psi dx\rightarrow_{n\rightarrow+\infty}  \int_{\R^N}\rho_0 v_0\psi dx.
\label{distribuy}
\end{equation}
\label{reminio}
In addition we have for any $n\in\mathbb{N}$:
\begin{equation}
\begin{aligned}
&\|v_0^n\|_{L^\infty}\leq \|v_0\|_{L^\infty}, \;\|\rho_0^n\|_{L^\infty}\leq\|\rho_0\|_{L^\infty}\\
&\|m_0^n\|_{bmo^{-1}\cap B^{-1}_{\infty,\infty}}\leq 4\mu\|\rho_0\|_{L^\infty}+2\|\rho_0\|_{L^\infty}\|v_0\|_{L^\infty}.
\end{aligned}
\end{equation}
\end{remarka}
\subsection*{Existence of weak solution in finite time for (\ref{0.1m}) when $\frac{1}{\rho_0}$ is not in $L^\infty(\R^N)$}
It remains now to prove that the sequence $(\rho_n,m_n,v_n)_{n\in\mathbb{N}}$ converges on $[0,T]$ up to a subsequence to a weak solution $(\rho,m,v)$ of the system (\ref{0.1m}) on $[0,T]$. From (\ref{infofin}) we have for any $n\in\mathbb{N}$ and $T$ defined as in (\ref{defM1}) that:
\begin{equation}
\|\rho_n u_n\|_{\widetilde{L}^1_T(B^{1}_{\infty,\infty})}+\|\rho_n u_n\|_{\widetilde{L}^\infty_T(B^{-1}_{\infty,\infty})}\leq C_1 M_1,
\label{hyperimp}
\end{equation}
with $C_1>0$ large enough.
From (\ref{hyperimp}) and interpolation arguments, we deduce that
the sequence $(\rho_n u_n)_{n\in\mathbb{N}}$ is uniformly bounded in $\widetilde{L}_T^{2-\alpha(\e)}(B^{\e}_{\infty,\infty})$ for $\e>0$ small enough and with $\alpha(\e)=\frac{2\e}{1+\e}$. By interpolation, we obtain that he sequence $(\rho_n u_n)_{n\in\mathbb{N}}$ is uniformly bounded in $L_T^{2-\alpha(\e)}(B^{\frac{\e}{2}}_{\infty,1})$.
We recall in addition that $B^{\frac{\e}{2}}_{\infty,1}$ is a Banach space. In addition we have:
$$\p_t(\rho_n u_n)=-{\rm div}(v_n\otimes (\rho_n u_n))+2\mu\D (\rho_n u_n)-\frac{a\gamma}{2\mu}\rho_n^\gamma v_n+\frac{a\gamma}{2\mu}\rho_n^\gamma u_n$$
From (\ref{hyperimp}) and (\ref{infofin}), we deduce that $\p_t(\rho_n u_n)$ is uniformly bounded in $L^{2-\alpha(\e)}_T(B^{-2}_{\infty,\infty}+B^{-1}_{\infty,\infty}+B^{0}_{\infty,\infty})$. We have used the fact that  $(\rho_n u_n)_{n\in\mathbb{N}}$ is uniformly bounded in $L_T^{2-\alpha(\e)}(B^{\frac{\e}{2}}_{\infty,1})$, it implies by interpolation that $(\rho_n u_n)_{n\in\mathbb{N}}$ is uniformly bounded in $L_T^{2-\alpha(\e)}(L^\infty)$. \\
We recall now the classical lemma of Aubin-Lions.
\begin{lemme}
\label{Aubin}
Let $X_0, X, X_1$ Banach spaces. Assume that $X_0$ is compactly embedded in $X$ and $X$ is continuously embedded in $X_1$. Let $1\leq p,q\leq +\infty$. We set for $T>0$:
$$W_T=\{u\in L^p([0,T],X_0),\,\frac{d}{dt}u\in L^q([0,T],X_1)\}.$$
Then if:
\begin{itemize}
\item $p<+\infty$ then the embedding of $W_T$ into $L^p([0,T],X)$ is compact.
\item $p=+\infty$ and $q>1$, then the embedding of $W_T$ into $C([0,T],X)$ is compact.
\end{itemize}
\label{Aubin}
\end{lemme}
We have seen that $(\rho_n u_n)_{n\in\mathbb{N}}$ is uniformly bounded in $L_T^{2-\alpha(\e)}(B^{\frac{\e}{2}}_{\infty,1})$.  In addition $\p_t(\rho_n u_n)$ is uniformly bounded in $L^{2-\alpha(\e)}_T(B^{-2}_{\infty,\infty}+B^{-1}_{\infty,\infty}+B^{0}_{\infty,\infty})$. Since $\varphi B^{\frac{e}{2}}_{\infty,1}$ is compactly embedded in $B^0_{\infty,1}\h L^\infty$ for any $\varphi\in C^\infty_c(\R^N)$, we deduce using the Aubin-Lions lemma that $(\rho_n u_n)_{n\in\mathbb{N}}$ converges up to a subsequence to $m$ in $L_T^{2-\alpha(\e)}(L_{loc}^\infty(\R^N))$. We note in the sequel $\alpha_3(\e)=2-\alpha(\e)$.\\
From (\ref{infofin}), we deduce that $(v_n)_{n\in\mathbb{N}}$ converges up to a subsequence weakly * to $v$ in $L^{\alpha_3(\e)'}([0,T],L^\infty(\R^N))$ with $\frac{1}{\alpha_3(\e)'}+\frac{1}{\alpha_3(\e)}=1$. Furthermore since $(\rho_n u_n)_{n\in\mathbb{N}}$ converges up to a subsequence to $m$ in $L_T^{\alpha_3(\e)}(L_{loc}^\infty)$, we deduce that up to a subsequence for any $\varphi\in C^\infty_c([0,T]\times\R^N)$ we have for any $i,j\in\{1\cdots N\}$:
\begin{equation}
\int^T_0\int_{\R^N}\rho_n u_n^{i}v_n^j \varphi dx dt\rightarrow_{n\rightarrow+\infty}\int^T_0\int_{\R^N} m_i v_j \varphi dx dt.
\label{compac1}
\end{equation}
Let us deal now with the pressure term $P_n(\rho_n)$. By definition of $v_n$ we have $\n\rho_n=\frac{1}{2\mu}(\rho_n v_n-\rho_n u_n),$
we deduce from (\ref{infofin}) and since $(\rho_n u_n)_{n\in\mathbb{N}}$ is uniformly bounded in $L_T^{2-\alpha(\e)}(L^\infty)$ that $(\n\rho_n)_{n\in\mathbb{N}}$ is uniformly bounded in $L^{2-\alpha(\e)}_{T}(L^\infty(\R^N))$ for $\e>0$ sufficiently small. It implies in particular that $(\rho_n)_{n\in\mathbb{N}}$ is uniformly bounded in $L^{2-\alpha(\e)}_{T}(W^{1,\infty}(\R^N))$ for $\e>0$ sufficiently small. Furthermore we deduce that $(\p_t\rho_n)_{n\in\mathbb{N}}$ is uniformly bounded in $L^{2-\alpha(\e)}_{T}(W^{-1,\infty}(\R^N))$. Using again the Aubin-Lions lemma we deduce that $(\rho_n)_{n\in\mathbb{N}}$ converges up to a subsequence to $\rho$ in $L^{2-\alpha(\e)}([0,T];L^\infty_{loc}(\R^N))$. Up to a subsequence, it implies that $(\rho_n)_{n\in\mathbb{N}}$ converges almost everywhere to $\rho$. Using (\ref{infofin}) and the fact that $(\rho_n)_{n\in\mathbb{N}}$ converges almost everywhere up to a subsequence to $\rho$, we deduce using the dominated convergence theorem that for any $\varphi\in C^\infty_c([0,T]\times\R^N)$ we have:
\begin{equation}
\int^T_0\int_{\R^N}P_n(\rho_n) \varphi dx dt\rightarrow_{n\rightarrow+\infty}\int^T_0\int_{\R^N} P(\rho)\varphi dx dt.
\label{compac2}
\end{equation}
From (\ref{infofin}) we have used in particular the fact that $P_n(\rho_n)=P(\rho_n)$ on $[0,T]$ for $n$ sufficiently large.
Multiplying now the third equation of (\ref{0.3}) by $\varphi\in C^\infty_c([0,T]\times\R^N)$ and integrating by parts, we have:
\begin{equation}
\begin{aligned}
&\int_{\R^N}\rho_nu_n(T,x)\cdot\varphi(T,x) dx-\int_{\R^N}\rho_nu_n(0,x)\cdot\varphi(0,x) dx-\int^T_0\int_{\R^N}\rho_n u_n(s,x)\cdot \p_s \varphi(s,x) ds dx\\
&-\sum_{i,j=1}^N \int^T_0\int_{\R^N}\rho_n v_n^{i} u_n ^j(s,x)\p_i \varphi_j(s,x) ds dx-2\mu\sum_{i,j=1}^N \int^T_0\int_{\R^N}\rho_n u_n^j(s,x)\p_{ii}\varphi_j(s,x) ds dx\\
&-\int^T_0 \int_{\R^N} P_n(\rho_n)(s,x){\rm div} \varphi(s,x) ds dx=0.
\end{aligned}
\end{equation}
From (\ref{compac1}), (\ref{compac2}) and since  $(\rho_n u_n)_{n\in\mathbb{N}}$ converges up to a subsequence to $m$ in $L_T^{2-\alpha(\e)}(L_{loc}^\infty)$, we have:
\begin{equation}
\begin{aligned}
&\int_{\R^N}m(T,x)\cdot\varphi(T,x) dx-\int_{\R^N}\rho_0 v_0(x)\cdot\varphi(0,x) dx-2\mu\int_{\R^N}\rho_0(x){\rm div}\va(0,x)dx\\
&-\int^T_0\int_{\R^N}m(s,x)\cdot \p_s \varphi(s,x) ds dx-\sum_{i,j=1}^N \int^T_0\int_{\R^N}m_j v_i (s,x)\p_i \varphi_j(s,x) ds dx\\
&-2\mu\sum_{i,j=1}^N \int^T_0\int_{\R^N}m_j (s,x)\p_{ii}\varphi_j(s,x) ds dx-\int^T_0 \int_{\R^N} P(\rho)(s,x){\rm div} \varphi(s,x) ds dx=0.
\end{aligned}
\label{solfaible1}
\end{equation}
We have use here the fact that for $\e>0$ small enough and any $\psi\in C^\infty_c(\R^N)$:
\begin{equation}
\begin{cases}
\begin{aligned}
&\rho_n^0u_n^0\rightarrow_{n\rightarrow+\infty}\rho_0u_0\;\;\mbox{in}\;{\cal D}'(\R^N),\\
&\psi \rho_n u_n(T,\cdot)\rightarrow_{n\rightarrow+\infty}\psi m(T,\cdot)\;\;\mbox{in}\;B^{-1-\e}_{\infty,\infty}.
\end{aligned}
\end{cases}
\label{impo5}
\end{equation}
Indeed we recall that $ \rho_n^0u_n^0=m_{0,1}^n-2\mu\n\rho_0^n$ and from
(\ref{distribuy}), we know that $(m_{0,1}^n)_{n\in\mathbb{N}}$ and $(\rho_0^n)_{n\in\mathbb{N}}$
converge respectively to $\rho_0v_0$ in the sense of distributions and $\rho_0$ weakly * in $L^\infty(\R^N)$. It implies then that 
$(\rho_n^0u_n^0)_{n\in\mathbb{N}}$ converges to $\rho_0 u_0$ in ${\cal D}'(\R^N)$.
\\
We recall now that $(\rho_n u_n)_{n\in\mathbb{N}}$ is uniformly bounded in $L^\infty([0,T],B^{-1}_{\infty,\infty})$ and that $(\p_t(\rho_n u_n))_{n\in\mathbb{N}}$ is uniformly bounded in $L^{2-\alpha(\e)}_T(B^{-2}_{\infty,\infty}+B^{-1}_{\infty,\infty}+B^{0}_{\infty,\infty})$. From the Aubin-Lions Lemma we deduce that for any $\psi\in C^\infty_c(\R^N)$, $(\psi\rho_n u_n)_{n\in\mathbb{N}}$ converges to $\psi \rho u$ up to a subsequence in $C([0,T],B^{-1-\e}_{\infty,\infty})$ for $\e>0$ small enough. It implies in particular
(\ref{impo5}). We can now pass to the limit in the term $\int_{\R^N}\rho_nu_n(T,x)\cdot\varphi(T,x) dx$ using the fact that the dual space of $B^{-1-\e}_{\infty,\infty}(\R^N)$ is $B^{1+\e}_{1,1}(\R^N)$ and writing $\varphi(T,\cdot)$ as follows $\varphi(T,\cdot)=\varphi(T,\cdot)\psi$ with $\psi\in C^\infty_c(\R^N)$ and $\psi$ equal to $1$ in a neighborhood of $\mbox{supp} \varphi(T,\cdot)$.
\\
\\
We have seen that  $(\rho_n)_{n\in\mathbb{N}}$ converges up to a subsequence to $\rho$ in $L^{2-\alpha(\e)}([0,T];L^\infty_{loc}(\R^N))$ and $(v_n)_{n\in\mathbb{N}}$ converges up to a subsequence weakly * to $v$ in $L^\infty([0,T],L^\infty(\R^N))$. It implies that  for any $\varphi\in C^\infty_c([0,T]\times\R^N)^N$ we have:
\begin{equation}
\int^T_0\int_{\R^N}\rho_n v_n \cdot\varphi dx dt\rightarrow_{n\rightarrow+\infty}\int^T_0\int_{\R^N} \rho v\cdot\varphi dx dt.
\label{compac3}
\end{equation}
Proceeding as previously we deduce that the first and third equations of (\ref{0.1m}) are verified by $(\rho,m,v)$ in the sense of distribution.
\subsection*{$(\rho,v,m)$ is in $Y_T$}
We know from (\ref{infofin}) that it exists $C>0$ such that for all $n\in\mathbb{N}$ we have:
\begin{equation}
\sup_{0<t<T,x\in\R^N}\frac{1}{t^\N}\int^t_0\int_{B (x,\sqrt{t})}|\rho_nu_n(s,x)|^2 ds dx\leq C
\label{imput}
\end{equation}
We have seen that up to a subsequence $(\rho_n u_n)_{n\in\mathbb{N}}$ converges almost everywhere on $(0,T)\times\R^N$ to $m$, we deduce then from (\ref{imput}) and the Fatou lemma that:
\begin{equation}
\sup_{0<t<T,x\in\R^N}\frac{1}{t^\N}\int^t_0\int_{B (x,\sqrt{t})}|m(s,x)|^2 ds dx\leq C
\label{imput1}
\end{equation}
From (\ref{imput1}) we deduce in particular that $m$ is in $L^2([0,T],L^2_{loc}(\R^N))$.\\
Similarly since $(\rho_n u_n)_{n\in\mathbb{N}}$ converges almost everywhere on $(0,T)\times\R^N$ to $m$, we deduce from (\ref{infofin}) that it exists $C>0$ large enough such that:
\begin{equation}
\sup_{t\in(0,T)} \sqrt{t}\|m(t,\cdot)\|_{L^\infty(\R^N)}\leq C
\end{equation}
In a similar way, we prove that $\rho$ is in $L^\infty([0,T],L^\infty(\R^N))$. Since the sequence $(v_n)_{n\in\mathbb{N}}$ converges weakly * to $v$ in $L^\infty([0,T]\times\R^N)$, we deduce that $v$ belongs to $L^\infty([0,T]\times\R^N)$.\\
Similarly we can prove that $m$ belongs to $\widetilde{L}^\infty([0,T],B^{-1}_{\infty,\infty})\cap \widetilde{L}^1((0,T),B^{1}_{\infty,\infty}(\R^N))$.
\subsection*{Existence of weak solution in finite time for (\ref{0.1}) when $\frac{1}{\rho_0}\in L^\infty(\R^N)$}
In this case, the main difference is that we are interesting in proving that the sequence $(\rho_nu_n\otimes u_n)_{n\in\mathbb{N}}$ converges in the sense of the distribution to a limit $\rho u\otimes u$. The first thing is to give a sense to $u$ and not only to $m$.\\
As previously we observe that $(\rho_n u_n)_{n\in\mathbb{N}}$ converges up to a subsequence to $m$ in $L_T^{\alpha_3(\e)}(L_{loc}^\infty)$ and that up to a subsequence  $(\rho_n u_n)_{n\in\mathbb{N}}$ converges almost everywhere to $m$ on $(0,T)\times\R^N$.
We are going now to prove that $(\frac{1}{\rho_n})_{n\in\mathbb{N}}$ remains bounded in $L^\infty([0,T],L^\infty(\R^N))$. We have seen that $(v_n)_{n\in\mathbb{N}}$ is uniformly bounded in $L^\infty([0,T],L^\infty(\R^N))$, furthermore $(\rho_n)_{n\in\mathbb{N}}$ verifies on $[0,T]$:
$$\p_t \rho_n-2\mu\D\rho_n+{\rm div}(\rho_n v_n)=0.$$
Using the maximum principle (see \cite{La}), we deduce that:
\begin{equation}
\|\frac{1}{\rho_n}\|_{L^\infty([0,T],L^\infty(\R^N))}\leq C(T),
\label{vide2}
\end{equation}
with $C$ a continuous fonction on $[0,T]$. We have seen previously that up to a subsequence $(\rho_n)_{n\in\mathbb{N}}$ converges to $\rho$ almost everywhere on $(0,T)\times\R^N$, it implies from (\ref{vide2}) that $(\frac{1}{\rho_n})_{n\in\mathbb{N}}$ converges to $\frac{1}{\rho}$ almost everywhere on $(0,T)\times\R^N$ and $\frac{1}{\rho}$ belongs to $L^\infty([0,T],L^\infty(\R^N))$.
Combining this last information and the fact that $m$ belongs to $L^2([0,T],L^2_{loc}(\R^N))$, we deduce that $u=\frac{m}{\rho}$ is in $L^2([0,T],L^2_{loc}(\R^N))$.
Now since  $(m_n)_{n\in\mathbb{N}}$ converges up to a subsequence to $m$ in $L_T^{\alpha_3(\e)}(L_{loc}^\infty)$ and $(\frac{1}{\rho_n})_{n\in\mathbb{N}}$ converges to $\frac{1}{\rho}$ in $L^p_{loc}([0,T]\times\R^N)$ for any $p\in[1,+\infty)$ (it suffices to use dominated convergence), we deduce that $(u_n)_{n\in\mathbb{N}}$ converges in the sense of distributions up to a subsequence to $\frac{m}{p}=u$.\\
\\
Similarly we verify that $(v_n)_{n\in\mathbb{N}}$ and $(u_n)_{n\in\mathbb{N}}$ are respectively uniformly bounded in $L^\infty([0,T],L^\infty(\R^N))$ and in $L^2([0,T],L^2_{loc}(\R^N))$, it implies since 
$\n\ln\rho_n=\frac{1}{2\mu}(v_n-u_n)$ that $(\n\ln\rho_n)_{n\in\mathbb{N}}$ is uniformly bounded in $L^2([0,T],L^2_{loc}(\R^N))$ and converges up to a subsequence weakly in $L^2_{loc}([0,T]\times\R^N)$ to $w$. We deduce in addition since $(\ln\rho_n)_{n\in\mathbb{N}}$ converges in $L^1_{loc}([0,T]\times\R^N)$ to $\ln\rho$ that $w=\n\ln\rho=\frac{1}{2\mu}(v-u)$. We have proved in particular that $\n\ln\rho\in L^2_{loc}([0,T]\times\R^N)$.
\\
From (\ref{solfaible1}), we recall that $(\rho,m,v)$ is a weak solution of (\ref{0.1m}) verifying for any  $\varphi\in C^\infty_c([0,T]\times\R^N)$:
\begin{equation}
\begin{aligned}
&\int_{\R^N}m(T,x)\cdot\varphi(T,x) dx-\int_{\R^N}\rho_0 v_0(x)\cdot\varphi(0,x) dx-2\mu\int_{\R^N}\rho_0(x){\rm div}\va(0,x)dx\\
&-\int^T_0\int_{\R^N}m(s,x)\cdot \p_s \varphi(s,x) ds dx-\sum_{i,j=1}^N \int^T_0\int_{\R^N}m_j v_i (s,x)\p_i \varphi_j(s,x) ds dx\\
&-2\mu\sum_{i,j=1}^N \int^T_0\int_{\R^N}m_j (s,x)\p_{ii}\varphi_j(s,x) ds dx-\int^T_0 \int_{\R^N} P(\rho)(s,x){\rm div} \varphi(s,x) ds dx=0.
\end{aligned}
\label{solfaible2}
\end{equation}
Since $u=\frac{m}{\rho}\in L^2_{loc}([0,T]\times\R^N)$ and $v=u+2\mu\n\ln\rho$, we have:
\begin{equation}
\begin{aligned}
&\sum_{i,j=1}^N \int^T_0\int_{\R^N}m_j v_i (s,x)\p_i \varphi_j(s,x) ds dx=\sum_{i,j=1}^N \int^T_0\int_{\R^N}\rho u_j (u_i+2\mu\p_i ln\rho) (s,x)\p_i \varphi_j(s,x) ds dx
\end{aligned}
\label{solfaible3}
\end{equation}
It is important to mention that the last integral is well defined since $\rho u_i u_j$ and $\rho u_j\p_i\ln\rho$ belongs to $L^1_{loc}([0,T]\times\R^N)$. In addition we observe that $\rho\p_i\ln\rho=\p_i\rho$. Indeed we have seen that $(\rho_n)_{n\in\mathbb{N}}$ converges
 up to a subsequence almost everywhere to $\rho$ on $[0,T]\times\R^N$ and is uniformly bounded in $L^\infty([0,T],L^\infty(\R^N))$. It implies that $(\rho_n)_{n\in\mathbb{N}}$ converges to $\rho$ in $L^p_{loc}([0,T]\times\R^N)$ for $p\in[1,+\infty)$. We recall now that $(\n\ln\rho_n)_{n\in\mathbb{N}}$ is uniformly bounded in $L^2([0,T],L^2_{loc}(\R^N))$ and converges up to a subsequence weakly in $L^2([0,T],L^2_{loc}(\R^N))$ to $\n\ln\rho$. It implies in particular that for any $\varphi\in (C^\infty_c(\R^N))^N$ we have:
\begin{equation}
\int^T_0\int_{\R^N}\rho_n\n\ln\rho_n\cdot\varphi dx ds\rightarrow_{n\rightarrow+\infty}\int^T_0\int_{\R^N}\rho\n\ln\rho\cdot\varphi dx ds.
\label{egal1}
\end{equation}
In addition we have similarly:
\begin{equation}
\int^T_0\int_{\R^N}\rho_n \n\ln\rho_n\cdot\varphi dx ds\rightarrow_{n\rightarrow+\infty}\int^T_0\int_{\R^N}\n\rho\cdot\varphi dx ds.
\label{egal2}
\end{equation}
Combining (\ref{egal1}) and (\ref{egal2}), we obtain that $\n\rho=\rho\n\ln\rho$.

From (\ref{solfaible2}) and (\ref{solfaible3}), we deduce that we have 
\begin{equation}
\begin{aligned}
&\int_{\R^N}m(T,x)\cdot\varphi(T,x) dx-\int_{\R^N}\rho_0 v_0(x)\cdot\varphi(0,x) dx-2\mu\int_{\R^N}\rho_0(x){\rm div}\va(0,x)dx\\
&-\int^T_0\int_{\R^N}\rho u(s,x)\cdot \p_s \varphi(s,x) ds dx-\sum_{i,j=1}^N \int^T_0\int_{\R^N}(\rho u_j  u_i +2\mu\p_i \rho)(s,x)\p_i \varphi_j(s,x) ds dx\\
&-2\mu\sum_{i,j=1}^N \int^T_0\int_{\R^N}\rho u_j (s,x)\p_{ii}\varphi_j(s,x) ds dx-\int^T_0\int_{\R^N}P(\rho)(s,x){\rm div}\varphi(s,x) ds dx=0
\end{aligned}
\label{solfaible4}
\end{equation}
We have then finally for any any $\varphi\in C^\infty_c([0,T)\times\R^N)$:
\begin{equation}
\begin{aligned}
&-\int_{\R^N}\rho_0 v_0(x)\cdot\varphi(0,x) dx-2\mu\int_{\R^N}\rho_0(x){\rm div}\va(0,x)dx\\
&-\int^{+\infty}_0\int_{\R^N}\rho u(s,x)\cdot \p_s \varphi(s,x) ds dx-\sum_{i,j=1}^N \int^{+\infty}_0\int_{\R^N}(\rho u_j  u_i +2\mu\p_i \rho)(s,x)\p_i \varphi_j(s,x) ds dx\\
&-2\mu\sum_{i,j=1}^N \int^{+\infty}_0\int_{\R^N}\rho u_j (s,x)\p_{ii}\varphi_j(s,x) ds dx-\int^{+\infty}_0\int_{\R^N}P(\rho)(s,x){\rm div}\varphi(s,x) ds dx=0
\end{aligned}
\label{solfaible5}
\end{equation}
Proceeding in a similar way for the first equation of (\ref{0.1}), we deduce that $(\rho,u)$ is a weak solution in finite time on $[0,T)$ of (\ref{0.1}).
\subsection*{Global weak solution for (\ref{0.1}) when $N=2$}
We are now interested in proving the existence of global weak solution $(\rho,u)$ for the system (\ref{0.1}) with the assumptions of the Theorem \ref{theo2}. We assume now in particular that  $(\rho_0-1)$, $v_0$ belong to $L^2(\R^N)$ and $u_0$, $v_0$ are in $L^1(\R^N)$. It implies in particular that $\n\rho_0$ is also in $L^1(\R^N)$. We consider the same approximate sequence $(\rho_n,m_n,v_n)_{n\in\mathbb{N}}$ as in the previous section except that we change the regularization on the initial density $\rho_0^n$ and we set:
\begin{equation}
\rho_0^n=\va_n j_n*(\rho_0-1)+1,
\label{edefinitu}
\end{equation}
with $j_n$ and $\va_n$ defined in (\ref{definitu}). Since it exists $c>0$ such that $\rho_0\geq c$, we have then $\rho_0^n\geq \va_n (c-1)+1\geq \min(c,1)>0. $ We define now $m_{0,1}^n$, $m_0^n$, $v_{0}^n$ and $u_{0}^n$ as in  (\ref{definitu}). We observe in particular that we have:
\begin{equation}
\begin{cases}
\begin{aligned}
&\|\rho_0^n-1\|_{L^2(\R^N)}\leq \|\rho_0-1\|_{L^2(\R^N)},\;\sup_{x\in\R^N}\rho_0^n(x)\leq\|\rho_0\|_{L^\infty}\\
&\|v_0^n\|_{L^2(\R^N)}\leq \frac{1}{\min(c,1)} \|\rho_0 v_0\|_{L^2(\R^N)},\;\|\rho_0^n v_0^n\|_{L^1(\R^N)}\leq  \|\rho_0 v_0\|_{L^1(\R^N)}\\
&\|\rho_0^n u_0^n\|_{L^1(\R^N)}\leq  \|\rho_0 v_0\|_{L^1(\R^N)}+\|\n\rho_0\|_{L^1(\R^N)}+\|\n\va_n\|_{L^N(\R^N)}\|\rho_0-1\|_{L^{N'}(\R^N)}\\
\end{aligned}
\end{cases}
\label{ininiy}
\end{equation}
We observe now that the approximate solution $(\rho_n,u_n,v_n)$ satisfies the BD entropy. 
In other words we have  far all $t>0$:
\begin{equation}
\begin{aligned}
&\frac{1}{2}\int_{\R^N} \rho_n(t,x)|v_n|^2(t,x) dx+\int_{\R^N} (\Pi_n(\rho_n(t,x))-\Pi_n(1)) dx+C_\gamma \int^t_0\int_{\R^N}|\n F_n(\rho_n)|^2( s,x) ds dx\\
&\hspace{3,5cm}=\frac{1}{2}\int_{\R^N} \rho_0^n(x)|v_0^n|^2(x) dx+\int^t_0\int_{\R^N} (\Pi_n(\rho^n_0(x))-\Pi_n(1)) dx
\end{aligned}
\label{aBD1}
\end{equation}
with $\Pi_n(s)=s(\int^s_1\frac{P_n(z)}{z^2} dz -P_n(1))$ and $F'_n(\rho_n)=\sqrt{\frac{P'_n(\rho_n)}{\rho_n}}$.
From (\ref{ininiy}), (\ref{aBD1}), (\ref{infofin}) and (\ref{vide2}), we deduce that $(\n\ln\rho_n)_{n\in\mathbb{N}}$
is uniformly bounded in $L^2([0,T],L^2(\R^N))$ with $T$ defined above for $n$ large enough. Indeed we have use the fact that for $n$ large enough we have $F'_n(\rho_n)=\sqrt{\frac{P'(\rho_n)}{\rho_n}}$.
 Similarly from (\ref{aBD1}) and (\ref{vide2}) we deduce that  $(v_n)_{n\in\mathbb{N}}$
is uniformly bounded in $L^\infty([0,T],L^2(\R^N))$ for $n$ large enough. Since $u_n=v_n-2\mu\n\ln\rho_n$, we deduce that $(u_n)_{n\in\mathbb{N}}$
is uniformly bounded in $L^2([0,T],L^2(\R^N))$ for $n$ large enough.
\\
Multiplying the momentum equation by $s u_n$ and integrating over $(0,t)\times\R^N$ with $0<t\leq T$, we have for $n$ large enough:
\begin{equation}
\begin{aligned}
&\frac{1}{2}\int_{\R^N} t\rho_n(t,x)|u_n|^2(t,x) dx-\frac{1}{2}\int^t_0\int_{\R^N}\rho_n(s,x)|u_n|^2(s,x) ds dx\\
&+2\mu\int^t_0\int_{\R^N}s \rho_n(s,x)|\n u_n|^2( s,x) ds dx+\int_{\R^N} t(\Pi_n(\rho_n(t,x))-\Pi_n(1)) dx\\
&\hspace{5cm}=\int^t_0\int_{\R^N} t(\Pi_n(\rho_n(t,x))-\Pi_n(1)) dx
\end{aligned}
\label{bBD1}
\end{equation}
From (\ref{aBD1}), (\ref{infofin}) and the fact that $(u_n)_{n\in\mathbb{N}}$
is uniformly bounded in $L^2([0,T],L^2(\R^N))$ for $n$ large enough, we deduce that $(\sqrt{t}\sqrt{\rho_n}u_n)_{n\in\mathbb{N}}$ is uniformly bounded in $L^\infty([0,T],L^2(\R^N))$ for $n$ large enough. From (\ref{vide2}) it implies that $(\sqrt{t}u_n)_{n\in\mathbb{N}}$ is uniformly bounded in $L^\infty([0,T],L^2(\R^N))$. We have then obtained from (\ref{aBD1}), (\ref{bBD1}), (\ref{infofin}),(\ref{vide2}) for $C(T)>0$ independent on $n$ (with $T$ defined as in (\ref{defM1}) and $n$ large enough) and any $t\in(0,T)$:
\begin{equation}
\begin{cases}
\begin{aligned}
&\|v_n\|_{ L^\infty([0,T],L^2(\R^N)\cap L^\infty(\R^N))}\leq C(T)\\
&\sqrt{t} \|u_n(t,\cdot)\|_{L^\infty(\R^N)\cap L^2(\R^N)}\leq C(T)\\
&\|\rho_n-1\|_{ L^\infty([0,T],L^2(\R^N))}\leq C(T)\\
&\sqrt{t}\|\n\rho_n(t,\cdot)\|_{L^2(\R^N)}\leq C(T)\\
&\sqrt{t}\|\rho_n(t,\cdot)\|_{W^{1,\infty}(\R^N)}\leq C(T).
\end{aligned}
\end{cases}
\label{superimptechab1}
\end{equation}
By interpolation we deduce that for any $p\in[2,+\infty]$ we have for any $t\in(0,T)$:
\begin{equation}
\begin{cases}
\begin{aligned}
&\sqrt{t} \|u_n(t,\cdot)\|_{L^p(\R^N)}\leq C(T),\\
&\sqrt{t} \| \rho_n(t,\cdot)-1\|_{W^{1,p}(\R^N)}\leq C(T).
\end{aligned}
\end{cases}
\label{superimptechab2}
\end{equation}
We are now going to estimate the $L^1$ norm of $(\rho_n u_n)_{n\in\mathbb{N}}$. We recall that we have since $u_n$ is irrotationnal:
$$
\begin{aligned}
\begin{cases}
&\p_t(\rho_n u_n)+{\rm div}(\rho_n v_n\otimes u_n)-2\mu\D(\rho_n u_n)+\frac{\rho_n P_n'(\rho_n)}{2\mu}(v_n-u_n)=0\\
&\p_t(\rho_n v_n)+{\rm div}(\rho_n u_n\otimes v_n)+\frac{\rho P_n'(\rho_n)}{2\mu}(v_n-u_n)=0.
\end{cases}
\end{aligned}
$$
Let $f$ a regular function from $\R$ to $\R$, we have then for any $j\in\{ 1,\cdots,N\}$:
$$
\begin{aligned}
\begin{cases}
&\p_t f(\rho_n u_n^j)+{\rm div}(v_n f(\rho_n u_n^j))+{\rm div}v_n (\rho_n u_n^j f'(\rho_n u_n^j)-f(\rho_n u^j_n))-2\mu\D f(\rho_n u_n^j)\\
&+2\mu f''(\rho_n u_n^j)|\n (\rho_n u_n^j)|^2+\frac{\rho_n P_n'(\rho_n)}{2\mu}(v_n^j-u_n^j) f'(\rho_n u_n^j)=0\\[2mm]
&\p_t(f(\rho_n v_n^j))+{\rm div}(u_n f(\rho_n v_n^j))+{\rm div}u_n (\rho_n v_n^j f'(\rho_n v_n^j)-f(\rho_n v_n^j))\\
&\hspace{6cm}+\frac{\rho_n P_n'(\rho_n)}{2\mu}(v_n^j-u_n^j) f'(\rho_n v_n^j)=0.
\end{cases}
\end{aligned}
$$
Taking in a classical way $f(x)=|x|$ (in fact we approximate $f$ by a regularizing sequence $(f_k)_{k\in\mathbb{N}}$ and we pass to the limit when $k$ goes to $+\infty$ with $f_k(x)=\sqrt{x^2+\frac{1}{k}}-\sqrt{\frac{1}{k}}$), we obtain by integrating over $(0,t)\times\R^N$ with $t\in(0,T)$:
\begin{equation}
\begin{aligned}
&\int_{\R^N}\rho_n|u_n^j|(t,x) dx\leq \int_{\R^N}\rho_n |u_n^j|(0,x) dx+\int^t_0\int_{\R^N} \frac{\rho_n P_n'(\rho_n)}{2\mu}(|v_n^j|+|u_n^j|) ds dx\\[2mm]
&\int_{\R^N}\rho_n|v_n^j|(t,x) dx+\int^t_0\int_{\R^N} \frac{\rho_n P_n'(\rho_n)}{2\mu}|v_n^j| dx ds \leq \int_{\R^N}\rho_n|v_n^j|(0,x) dx\\
&\hspace{7cm}+\int^t_0\int_{\R^N} \frac{\rho_n P_n'(\rho_n)}{2\mu}|u_n^j| ds dx
\end{aligned}
\label{NormeL1}
\end{equation}
For $t\in(0,T)$ we deduce from (\ref{ininiy}), (\ref{infofin}), (\ref{NormeL1}), (\ref{Pressiontech}) that we have using the Gronwall lemma for $n$ large enough:
\begin{equation}
\begin{aligned}
\|\rho_n u_n(t,\cdot)\|_{L^1(\R^N)}+\|\rho_n v_n(t,\cdot)\|_{L^1(\R^N)}\leq C(T),
\end{aligned}
\label{NormeL11}
\end{equation}
with $C(T)>0$ sufficiently large and independent on $n$.\\
Using Besov embedding and arguments of interpolation, we have from (\ref{superimptechab1}), (\ref{superimptechab2}), (\ref{vide2}),
(\ref{infofin}) for any $t\in(0,T)$ and $n$ large enough:
\begin{equation}
\begin{cases}
\begin{aligned}
&\|u_n(t,\cdot)\|_{\dot{B}^{\NN-1}_{p,1}(\R^N)\cap \dot{B}^{\NN-1+\e_p}_{p,1}(\R^N)}\leq \frac{C(T)}{\sqrt{t}},\\
&\| \rho_n(t,\cdot)-1\|_{\dot{B}^{\NN}_{p,1}(\R^N)\cap \dot{B}^{\NN+\e_p}_{p,1}(\R^N)}\leq \frac{C(T)}{\sqrt{t}}\\
&\|(\rho_n,\frac{1}{\rho_n})(t,\cdot)\|_{L^\infty(\R^N)}\leq C(T)\\
&\|\rho_n(t)-1\|_{L^2(\R^N)\cap L^\infty(\R^N)}\leq C(T)\\
&\|\n\sqrt{\rho_n}(t,\cdot)\|_{L^2(\R^N)}\leq\frac{C(T)}{\sqrt{t}}\\
&\|u_n(t,\cdot)\|_{L^1(\R^N)}+\|v_n(t,\cdot)\|_{L^1(\R^N)}\leq C(T)\\
&\|u_n(t,\cdot)\|_{L^2(\R^N)\cap L^\infty(\R^N)}\leq\frac{C(T)}{\sqrt{t}},\;\|v_n(t,\cdot)\|_{L^2(\R^N)\cap L^\infty(\R^N)}\leq C(T),
\end{aligned}
\end{cases}
\label{superimptechab3}
\end{equation}
with $p>N$, $\e_p>0$ sufficiently small such that $\NN-1+\e_p<0$ and $C(T)>0$ large enough depending only on $T$.
We are now going to prove global time estimates on $(\rho_n,u_n,v_n)_{n\in\mathbb{N}}$. 
From (\ref{superimptechab3}) we deduce that it exists $C(T)>0$ independent on $n$ such that for $p\in]N,2N[$ (and $\e_p>0$ small enough) we have with $T$ defined as above and $n$ large enough:
\begin{equation}
\begin{cases}
\begin{aligned}
&\|u_n(\frac{T}{2},\cdot)\|_{\dot{B}^{\NN-1}_{p,1}(\R^N)\cap \dot{B}^{\NN-1+\e_p}_{p,1}(\R^N)}\leq C(T),\\
&\| \rho_n(\frac{T}{2},\cdot)-1\|_{\dot{B}^{\NN}_{p,1}(\R^N)\cap \dot{B}^{\NN+\e_p}_{p,1}(\R^N)}\leq  C(T)\\
&\|(\rho_n,\frac{1}{\rho_n})(\frac{T}{2},\cdot)\|_{L^\infty(\R^N)}\leq C(T)\\
&\|\rho_n(\frac{T}{2})-1\|_{L^2(\R^N)\cap L^\infty(\R^N)}\leq C(T)\\
&\|\n\sqrt{\rho_n}(\frac{T}{2},\cdot)\|_{L^2(\R^N)}\leq  C(T)\\
&\|u_n(\frac{T}{2},\cdot)\|_{L^2(\R^N)\cap L^\infty(\R^N)}\leq C(T),\;\|v_n(t,\cdot)\|_{L^2(\R^N)\cap L^\infty(\R^N)}\leq C(T)\\
&\|u_n(\frac{T}{2},\cdot)\|_{L^1(\R^N)}+\|v_n(\frac{T}{2},\cdot)\|_{L^1(\R^N)}\leq C(T)\\
\end{aligned}
\end{cases}
\label{superimptechab5}
\end{equation}
In addition we have seen that $\rho_n(\frac{T}{2},\cdot)$ is radial and $u_n(\frac{T}{2},\cdot)$ is radially symmetric.
Using the Theorem \ref{theo4}, there exists a global strong solution $(\rho_{n,1},u_{n,1})_{n\in\mathbb{N}}$ on $(\frac{T}{2},+\infty)$ with initial data $(\rho_n(\frac{T}{2},\cdot),u_n(\frac{T}{2},\cdot))$ at the time $\frac{T}{2}$. Furthermore there exists a continuous function $C_p$ on $[0,+\infty)$ independent on $n$ such that for any $t>\frac{T}{2}$  we have:
\begin{equation}
\begin{cases}
\begin{aligned}
&\|u_{n,1}\|_{\widetilde{L}^\infty([\frac{T}{2},t],\dot{B}^{\NN-1}_{p,1}\cap \dot{B}^{\NN-1+\e_p}_{p,1})}+\|u_n\|_{\widetilde{L}^1([\frac{T}{2},t],\dot{B}^{\NN+1}_{p,1}\cap \dot{B}^{\NN+1+\e_p}_{p,1})}\leq C(t-\frac{T}{2}),\\
&\|\rho_{n,1}-1\|_{\widetilde{L}^\infty([\frac{T}{2},t],\dot{B}^{\NN}_{p,1}\cap \dot{B}^{\NN+\e_p}_{p,1})}\leq C(t-\frac{T}{2}),\\
&\|(\rho_{n,1},\frac{1}{\rho_{n,1}})\|_{L^\infty([\frac{T}{2},t],L^\infty)}\leq C(t-\frac{T}{2}).
\end{aligned}
\end{cases}
\label{crucialuni}
\end{equation}
Now since $(\rho_n,u_n)_{n\in\mathbb{N}}$ is a regular solution, we recall that $(\rho_n,u_n)_{n\in\mathbb{N}} $ verify also (\ref{estimo1}), and we deduce by uniqueness of the strong solution for compressible Navier-Stokes equation (see \cite{Fourier}) that:
\begin{equation}
\rho_n=\rho_{n,1},\;u_n=u_{n,1}\;\;\mbox{on}\;\;[\frac{T}{2},+\infty),
\label{mitech}
\end{equation}
for $n$ large enough.
From (\ref{mitech}) and (\ref{crucialuni}) we deduce that $(\rho_n,u_n)_{n\in\mathbb{N}}$ verifies uniformly in $n$ the estimates (\ref{crucialuni}) for $n$ large enough.  It is then easy now to prove the existence of global weak solution when $N=2$ . It suffices to combine the ideas that we use to prove the existence of weak solution in finite time (typically we focus on the interval $[0,\frac{T}{2}]$) and the classical arguments developed in \cite{Danchin} to prove the existence of strong solution (at this level, we consider the interval $(\frac{T}{2},+\infty)$). We note that the solution $(\rho,u)$ satisfies in addition:
\begin{equation}
\begin{cases}
\begin{aligned}
&(\rho-1,u)\in \widetilde{L}_{loc}^\infty([\frac{T}{2},+\infty),\dot{B}^{\NN}_{p,1}\cap \dot{B}^{\NN+\e_p}_{p,1})\\
&\hspace{0,1cm}\times (\widetilde{L}_{loc}^\infty([\frac{T}{2},+\infty),\dot{B}^{\NN-1}_{p,1}\cap \dot{B}^{\NN-1+\e_p}_{p,1})\cap \widetilde{L}_{loc}^1([\frac{T}{2},+\infty),\dot{B}^{\NN+1}_{p,1}\cap \dot{B}^{\NN+1+\e_p}_{p,1}))^N\\[2mm]
&(\rho,\frac{1}{\rho})\in L_{loc}^\infty((0,+\infty),L^\infty).
\end{aligned}
\end{cases}
\end{equation}

\section{Proof of the Theorem \ref{theo3}}
The proof follows exactly the same lines as the proof of the Theorem \ref{theo2}. Indeed in one dimension, we can rewrite the system (\ref{0.1}) as follows:
$$
\begin{cases}
\begin{aligned}
&\p_t\rho+\p_x(\rho u)=0\\
&\p_t(\rho u)+\p_x(\rho u^2)-2\mu\p_x(\rho\p_x u)+\p_x P(\rho)=0.
\end{aligned}
\end{cases}
$$
The crucial point is that in one dimension, the effective velocity $v=u+2\mu\p_x\ln\rho$ satisfies a damped transport equation:
$$\p_t v+u\p_x v+\p_x F(\rho)=0,$$
with $F'(\rho)=\frac{P'(\rho)}{\rho}$. In particular it does not require any geometric assumption in order to ensure that the solution remains irrotationnal (roughly speaking it is always the case in one dimansion). We can then apply the same estimates as in the previous proof what is sufficient to show the Theorem \ref{theo3}.
\label{section5}
\section*{Appendix}
In this Appendix, we wish to understand the form of the nonlinear terms in the basis $(e_{r,x},e_{\theta,x},e_{\phi,x})$ of $\R^3$ with $x\in\R^3\backslash\{x_1=x_2=0\}$. We assume here that $(\rho,u)$ is a regular solution of the system (\ref{0.1}) on $(0,T^*)$ verifying the geometric assumptions (\ref{densradial1}), (\ref{decomporad}):
\begin{equation}
\begin{cases}
\begin{aligned}
&\rho(t,x)=\rho_1(t,|x|),\\
&u(t,x)=u_{1}(t,|x|) e_{r,x}+u_{2}(t,|x'|,x_3)e_{\theta,x}+u_{3}(t,|x'|,x_3)e_{\phi,x}\\
&u^1(t,x)=u_{1}(t,|x|) e_{r,x},\,u^2(t,x)=u_{2}(t,|x'|,x_3)e_{\theta,x},\,u^3(t,x)=u_{3}(t,|x'|,x_3)e_{\phi,x}.
\end{aligned}
\end{cases}
\label{decomporadk}
\end{equation}
Setting $r=|x|$ and $r'=\sqrt{x_1^2+x_2^2}$ with $x\in\R^3\backslash \{x_1=x_2=0\}$, we will defined in the sequel the derivatives $\p_r u_1$, $\p_{rr}u_1$, $\p_{r'}u_2$, $\p_{r'r'}u_2$, $\p_3 u_2$, $\p_{3, r'}u_2$...First we can verify that these derivatives are well defined on $R^3\backslash \{x_1=x_2=0\}$, let us consider as an example $\p_r u_1$; we have then when
$|x|=r$, $\lambda(h)=1+\frac{h}{r}$:
$$
\begin{aligned}
\p_r u_1(t,r)&=\lim_{h\rightarrow 0}\frac{ u_1(t,r+h)-u_1(t,r)}{h}=\lim_{h\rightarrow 0}\frac{ u(t,\lambda(h) x)\cdot x-u(t,x)\cdot x}{h|x|} \\
&=\lim_{h\rightarrow 0}\frac{ (u(t,\lambda(h) x)-u(t,x))\cdot x}{h|x|}
\end{aligned}
$$
We have then for $|x|>0$:
$$
\begin{aligned}
&\p_r u_1(t,r)=Du(t,x)(x)\cdot\frac{x}{r^2}=\frac{(\n u(t,x)\cdot x)\cdot x}{|x|^2}.
\end{aligned}
$$
We can observe that this formula does not depend on the choice of $x$. Indeed if $A$  is an isometry, we have since $^t A u(t,Ax)=u(t,x)$ for any $x\in\R^3$:
$$Du(t,Ax)(Ax)\cdot\frac{Ax}{r^2}=Du(t,x)(x)\cdot\frac{x}{r^2}.$$
Similarly we can show that $u_1(t,r), u_2(t,r',x_3),u_3(t,r',x_3)$ are $C^\infty$ in $\R^3\backslash\{x_1=x_2=0\}$. 
Let us give now the form of the nonlinear terms in the momentum equation of (\ref{0.1}).
Few computations give for $x\in\R^3\backslash \{x=0\}$ and $t\in(0,T^*)$:
\begin{equation}
\begin{cases}
\begin{aligned}
&2\mu{\rm div}(\rho D u^1)(t,x)=2\mu{\rm div}(\rho \n u^1)(t,x)=B_1(t,|x|)\frac{x}{|x|}\\
&\rho u^1\cdot\n u^1(t,x)=B_2(t,|x|)\frac{x}{|x|}\\
&\n P(\rho)(t,x)= (P(\rho_1))_r (t,|x|)\frac{x}{|x|}
\end{aligned}
\end{cases}
\label{calcu1}
\end{equation}
with $B_1(t,|x|)=|x|^2(\frac{\rho_1}{r}(\frac{u_1}{r})_r)_r(t,|x|)+4\rho_1(t,|x|)(\frac{u_1}{r})((t,|x|)+(\frac{\rho_1 u_1}{r})_r(t,|x|)$, $B_2(t,|x|)=\rho_1 u_1\p_r u_1(t,|x|)$ regular functions on $\R^3\backslash \{x=0\}$.
Similarly tedious calculus give for $x\in \R^3\backslash \{x_1=x_2=0\}$ and $t\in(0,T^*)$ (with $r'=|x'|$):
\begin{equation}
\begin{cases}
\begin{aligned}
& u^2\cdot\n u^1(t,x)=\frac{u_2(t,x)u_1(t,|x|)}{|x|} e_{\theta,x},\\
&u^1\cdot\n u^2(t,x)=(\frac{|x'|}{|x|}u_1(t,|x|)\p_{r'}u_2(t,r',x_3)+\frac{x_3}{|x|}u_1(t,|x|)\p_3 u_2(t,r',x_3))e_{\theta,x}\\
&u^2\cdot\n u^2(t,x)=-e_{r,x}(\frac{u_2^2(t,r',x_3)}{|x|})+e_{\theta,x}\big(-\frac{|x'|}{|x|}u_2(t,r',x_3)\p_3 u_2(t,r',x_3)\\
&+\frac{x_3}{|x|}u_2(t,r',x_3)\p_{r'}u_2(t,r',x_3)\big).
\end{aligned}
\end{cases}
\label{calculs1a}
\end{equation}
We wish now to estimate $\D u^2$, we get in a first time for $x\in \R^3\backslash \{x_1=x_2=0\}$:
\begin{equation}
\begin{cases}
\begin{aligned}
& \D(\frac{x_1 x_3}{|x'||x|}u_2(t,|x'|,x_3))=\frac{x_1 x_3}{|x'||x|}\big(-u_2(t,r',x_3)(\frac{2}{|x|^2}+\frac{1}{|x'|^2})\\
&+\p_{r'}u_2(t,r',x_3) (\frac{1}{|x'|}-\frac{2|x'|}{|x|^2})+\p^2_{r' ,r'}u_2(t,r',x_3)+\p^2_{3 ,3}u_2(t,r',x_3)\big)\\
&+2\frac{x_1}{|x|}\p_3 u_2(t,r',x_3) \frac{|x'|}{|x|^2}\\
&\D(-\frac{|x'|}{|x|}u_2(t,|x'|,x_3))=-\frac{|x'|}{|x|}\big(-u_2(t,r',x_3)(\frac{2}{|x|^2}-\frac{1}{|x'|^2})+\p_{r'}u_2(t,r',x_3)(\frac{3}{|x'|}-2\frac{|x'|}{|x|^2})\\
&+\p^2_{r' ,r'}u_2(t,r',x_3)+\p^2_{3 ,3}u_2(t,r',x_3)\big)+2\frac{x_3}{|x|}\p_3 u_2 (t,r',x_3)\frac{|x'|}{|x|^2}\\
\end{aligned}
\end{cases}
\label{calculs2a}
\end{equation}
We deduce that we have $x\in \R^3\backslash \{x_1=x_2=0\}$ and $t\in(0,T^*)$:
\begin{equation}
\begin{aligned}
&\D u^2(t,x)=\big(-u_2(\frac{2}{|x|^2}+\frac{1}{|x'|^2})+\p_{r'}u_2(\frac{1}{|x'|}-2\frac{|x'|}{|x|^2})+\p_{r' ,r'}u_2+\p_{3 ,3}u_2\big)e_{\theta,x}\\
&+2 \p_3 u_2 \frac{|x'|}{|x|^2} e_{r,x}+2\begin{pmatrix} 
0 \\
0 \\
-\frac{u_2}{|x'|\,|x|}-\frac{\p_{r'}u_2}{|x|}
\end{pmatrix}
\end{aligned}
\end{equation}
From (\ref{algebre}), we deduce that:
\begin{equation}
\begin{aligned}
&\D u^2(t,x)=\big(-u_2(t,r',x_3)(\frac{2}{|x|^2}+\frac{1}{|x'|^2})+\p_{r'}u_2(t,r',x_3)(\frac{1}{|x'|}-2\frac{|x'|}{|x|^2})\\
&+\p^2_{r' ,r'}u_2(t,r',x_3)+\p^2_{3 ,3}u_2(t,r',x_3)\big)e_{\theta,x}+2 \p_3 u_2(t,r',x_3) \frac{|x'|}{|x|^2} e_{r,x}\\
&-2(\frac{u_2(t,r',x_3)}{|x'|\,|x|}+\frac{\p_{r'}u_2(t,r',x_3)}{|x|})(\frac{x_3}{|x|} e_{r,x}-\frac{|x'| }{|x|}e_{\theta,x}).
\end{aligned}
\label{Lapla}
\end{equation}
Similarly we have $x\in \R^3\backslash \{x_1=x_2=0\}$ and $t\in(0,T^*)$:
\begin{equation}
\begin{aligned}
&{\rm div}u^2(t,x)=\frac{x_3}{|x||x'|}u_2(t,r',x_3)+\frac{x_3}{|x|}\p_{r'}u_2(t,r',x_3)-\frac{|x'|}{|x|}\p_3 u_2(t,r',x_3)\\
&\n{\rm div} u^2(t,x)=B(t,|x'|,x_3)e_{\theta, x}+C(t,|x'|,x_3)e_{r,x},
\end{aligned}
\label{Lapla1}
\end{equation}
with $B$ and $C$ functions depending on $(t,|x'|,x_3)$ and regular on $ \R^3\backslash \{x_1=x_2=0\}$. Finally we obtain $x\in \R^3\backslash \{x_1=x_2=0\}$ and $t\in(0,T^*)$:
\begin{equation}
\n\rho\cdot ^t\n u^2(t,x)=-\frac{\p_r\rho(t,|x|)u_2(t,|x'|,x_3)}{|x|}e_{\theta,x}.
\label{geometrie5}
\end{equation}
We have finally:
$$\n\rho\cdot\n u^2=[\frac{|x'|}{|x|}\p_r\rho\p_{r'}u_2+\frac{x_3}{|x|}\p_3 u_2\p_r\rho]e_{\theta,x}.$$
\section*{Acknowledgements}
The author has been partially funded by the ANR project 
INFAMIE ANR-15-CE40-0011. This work was realized during the secondment of the author in the ANGE Inria team.

\end{document}